\documentclass[10pt]{article}
\usepackage{amsmath,amssymb}
\makeatletter\@addtoreset {equation}{section}\makeatother

\newtheorem{theo}{Theorem}[section]
\newtheorem{lem}[theo]{Lemma}

\newtheorem{prop}[theo]{Proposition}
\newtheorem{cor}[theo]{Corollary}
\newtheorem{rem}[theo]{Remark}
\newenvironment{Proof}
{\begin{trivlist} \item[]{\bf Proof. }}%
{\hspace*{\fill}$\rule{.3\baselineskip}{.35\baselineskip}$\end{trivlist}}

\newcommand{\R}{\mathbb{R}}

\newcommand{\N}{\mathbb{N}}
\newcommand{\Q}{\mathbb{Q}}

\renewcommand{\geq}{\geqslant}
\renewcommand{\leq}{\leqslant}

\renewcommand{\phi}{\varphi}
\newcommand{\be}{\begin{eqnarray}}
\newcommand{\ee}{\end{eqnarray}}

\newcommand{\eps}{\varepsilon}

\begin{document}
\title{\bf The ground state of two coupled Gross--Pitaevskii equations in
  the Thomas--Fermi limit}
\author{Cl\'ement Gallo
}
\date{}
\maketitle
\abstract{We prove existence and uniqueness of a positive solution to
  a system of two coupled Gross-Pitaevskii equations. We give a full
  asymptotic expansion of this solution into powers of the semi
  classical parameter $\eps$ in the Thomas--Fermi limit $\eps\to 0$.}

\section{Introduction}
Recent experiments with Bose--Einstein condensates \cite{PitStr}
have stimulated new interest in the Gross--Pitaevskii equation
with a harmonic potential. This equation can be written as
\begin{equation}
\label{GP0} i \eps u_t + \eps^2 \Delta u + (1 - |x|^2) u - |u|^2 u = 0,
\quad x \in \R^d, \quad t \in \R_+,
\end{equation}
where $u(t,x)$ denotes the complex valued wave function of the Bose
gas, and $\eps$ is a small parameter. The limit $\eps\to 0$ corresponds to the
Thomas--Fermi approximation of a nearly compact atomic cloud
\cite{Fermi}, \cite{Thomas}. At equilibrium and in the absence of
rotation, the condensate is described by the ground--state, which is a
positive, time independent solution $u(t,x)=\eta_\eps(x)$ to
(\ref{GP0}). The ground state minimizes the
Gross--Pitaevskii energy  
\be\label{GP-en}
E_\eps(u)=\int_{\R^d}\left(\eps^2|\nabla
u|^2+(|x|^2-1)|u|^2+\frac{1}{2}|u|^4\right)dx
\ee
among functions with finite energy. The understanding of the profile
of the ground state is particularily important \cite{A}. It is well
known (see for instance
\cite{IM}) that in the
Thomas--Fermi limit $\eps\to 0$, the ground state $\eta_\eps$
converges to the Thomas--Fermi's compactly supported function
\begin{equation}
\label{Thomas-Fermi}
\eta_0(x) = \left\{ \begin{array}{cl} (1 - |x|^2)^{1/2} & \text{for }  |x| < 1, \\
0 & \text{for }  |x| > 1. \end{array} \right.
\end{equation}
The function $\eta_0$ has a singularity at $|x|=1$, whereas for
$\eps>0$, $\eta_\eps$ is regular. The question of the description of the
behaviour of $\eta_\eps$ close to the turning point $|x|=1$ as
$\eps\to 0$ has been adressed by Dalfovo, Pitaevskii and Stringari
\cite{DPS} and by Fetter and Feder \cite{FF} on a formal level. Among other reasons, this question is relevant because an important part of the kinetic energy is concentrated in the region $|x|\approx 1$ (see also \cite{G}). In
particular, it is shown in \cite{DPS} and \cite{FF} that it is possible to describe $\eta_\eps$ close to $|x|=1$ as $\eps\to 0$
thanks to solutions of the Painlev\'e II equation. This analysis has
been made rigorous in \cite{GP}, where a full asymptotic expansion of
$\eta_\eps$ in terms of powers of $\eps^{2/3}$ is calculated. The
proof consists in introducing a new variable $y=(1-|x|^2)/\eps^{2/3}$
that blows up the solution close to the turning point
$|x|=1$, writing $\eta_\eps(x)=\eps^{1/3}\nu_\eps(y)$ and solving
the equation satisfied by $\nu_\eps$ in terms of the variable $y$. It
turns out that the variable $y$ makes it possible to describe the behaviour of
$\eta_\eps$ as $\eps\to 0$ not only close to the turning point, but also globally for all $x\in \R^d$. In \cite{KS}, Karali
and Sourdis have extended this result to more general potentials.

The purpose of this paper is to adapt the result obtained in \cite{GP}
to the case of a two--component Bose--Einstein condensate. As we shall
see, one of the new difficulties we are facing to is that the
ground state has now two turning points instead of one in the case of
a scalar Gross-Pitaevskii equation. As a matter of fact, it will be
necessary to use three different variables to describe the ground
state, instead of one for the scalar equation. Then,
denoting by $\eta_1$ and $\eta_2$ the wave functions of the two
components, $\eta_1$ and $\eta_2$ solve the following system of two
coupled Gross-Pitaevskii equations with quadratic potentials,
\be\label{sys0}
\left\{\begin{array}{l}\eps^2\Delta\eta_{1}+\left(\mu_1-|x|^2\right)\eta_{1}-2\alpha_1\eta_{1}^3-2\alpha_0\eta_{2}^2\eta_{1}=0\\
\eps^2\Delta\eta_{2}+\left(\mu_2-|x|^2\right)\eta_{2}-2\alpha_2\eta_{2}^3-2\alpha_0\eta_{1}^2\eta_{2}=0,\end{array}\right.
\ee
where $\alpha_0,\alpha_1,\alpha_2>0$, $\mu_1,\mu_2>0$ are chemical
potentials, $\eps$ is a small parameter and $x\in\R^d$ where the
dimension $d$ is 1,2 or 3. Ground states of this system have also been studied in the case $d=2$ and with different methods by Aftalion, Noris and Sourdis \cite{ANS}. They prove various estimates on the difference between the Ground state and the Thomas-Fermi limit, which can be recovered by using the full asymptotic expansion of the ground state we prove here. 

For convenience, we define
$$\Gamma_1=1-\frac{\alpha_0}{\alpha_1},\quad
\Gamma_2=1-\frac{\alpha_0}{\alpha_2},\quad
\Gamma_{12}=1-\frac{\alpha_0^2}{\alpha_1\alpha_2}.$$
We will consider here only values of the parameters such that the two
components of the Thomas--Fermi limit $(\eta_{10},\eta_{20})$ are
supported and do not vanish on disks centered at $x=0$, in opposition
with other cases where one component is supported in an annulus
and the other one in a disk. More specific conditions are given below.
One of the differences between this case and the one component case is
that, as we shall see in the next section, the Thomas--Fermi limit
$(\eta_{10},\eta_{20})$ has now two turning points. Thus, we have to
introduce two different new variables. We will still be able to
give a full asymptotic expansion of $(\eta_1,\eta_2)$ into powers of
$\eps$ in the limit $\eps\to 0$, but functions of each of these
two new variables will appear in every term of the expansion. 

\subsection{Calculation of the Thomas-Fermi limit}
We are interested in solutions of
(\ref{sys0}) which converge in the Thomas-Fermi limit $\eps\to 0$ to
functions $\eta_{10}$ and $\eta_{20}$ which are both supported in a
disk, with respective radii $R_1$ and $R_2$ (for $j=1,2$, $R_j=\inf\left\{R>0,{\rm Supp}\eta_{j0}\subset B(0,R)\right\}$), and such that
$(\eta_{10},\eta_{20})$ solves (\ref{sys0}) with $\eps=0$. Let us
recall the arguments leading to the expression of the Thomas--Fermi
profile $(\eta_{10},\eta_{20})$ of the ground state, as it has been done in \cite{AMW}. Up to a change of the
indices, we assume (see Remark \ref{R1R2} below for the case $R_1=R_2$)
$$R_1<R_2.$$
From our definition of $R_1$ and $R_2$, we have $\eta_{10}(x)=\eta_{20}(x)=0$ for $|x|\geq
R_2$. For $R_1<|x|<R_2$, $\eta_{10}(x)=0$, and the second equation in
(\ref{sys0}) implies
$$\eta_{20}(x)^2=\frac{\mu_2-|x|^2}{2\alpha_2}.$$
Thus, 
\begin{eqnarray}\label{mu2}
\mu_2=R_2^2,
\end{eqnarray}
and $\eta_{20}(x)>0$ for $|x|=R_1$, which implies that
$\eta_{10}(x)\neq 0$ and $\eta_{20}(x)\neq 0$ for $|x|\approx R_1$,
$|x|<R_1$. If $\eps=0$, $\eta_1\neq 0$ and $\eta_2\neq 0$, then (\ref{sys0}) can be
rewritten into a non-homogeneous linear system in the variables $\eta_1^2,\eta_2^2$. Solving this
system (for the peculiar case $\Gamma_{12}=0$, see Remark \ref{Gamma120} below), we get, for $|x|\approx R_1$ and $|x|<R_1$,
\begin{eqnarray}
\eta_{10}(x)^2=\frac{1}{2\alpha_1\Gamma_{12}}\left(\mu_1-\frac{\alpha_0}{\alpha_2}\mu_2-\Gamma_2|x|^2\right),\label{et10}\\
\eta_{20}(x)^2=\frac{1}{2\alpha_2\Gamma_{12}}\left(\mu_2-\frac{\alpha_0}{\alpha_1}\mu_1-\Gamma_1|x|^2\right).\label{et20}
\end{eqnarray}
In particular, since $\eta_{10}$ vanishes on the sphere $|x|=R_1$ (or
equivalently, using the continuity of $\eta_{20}$ on the same sphere), we deduce
\begin{eqnarray}\label{mu1}
\mu_1=\frac{\alpha_0}{\alpha_2}\mu_2+\Gamma_2R_1^2=\frac{\alpha_0}{\alpha_2}R_2^2+\Gamma_2R_1^2.
\end{eqnarray}
Moreover, we also infer from the positiveness of $\eta_{10}^2$ and
(\ref{et10}) that the condition
\begin{eqnarray}\label{posgamma}
\Gamma_2/\Gamma_{12}>0
\end{eqnarray}
has to be satisfied.
Finally, (\ref{et20}) and the assumption of positiveness of $\eta_2$ on the disk with
radius $R_2$ (and not on an annulus) yields
$$\frac{1}{\Gamma_{12}}\left(\mu_2-\frac{\alpha_0}{\alpha_1}\mu_1\right)>0,$$
which can be rewritten in terms of $R_1$ and $R_2$ as
\begin{eqnarray}\label{cond-disk}
R_2^2>\frac{\alpha_0}{\alpha_1}\frac{\Gamma_2}{\Gamma_{12}}R_1^2.
\end{eqnarray}
As a result, provided that the parameters satisfy conditions (\ref{posgamma}) and
(\ref{cond-disk}), 
\begin{eqnarray}\label{eta10}
\eta_{10}(x)=\left\{\begin{array}{ccc}\left(\frac{\Gamma_2}{2\alpha_1\Gamma_{12}}\right)^{1/2}(R_1^2-|x|^2)^{1/2}&
    {\rm if }& |x|\leq R_1\\0& {\rm if }& |x|\geq R_1\end{array}\right.
\end{eqnarray}
and
\begin{eqnarray}\label{eta20}
\eta_{20}(x)=\left\{\begin{array}{ccc}\left(\frac{R_2^2-R_1^2}{2\alpha_2}+\frac{\Gamma_1}{2\alpha_2\Gamma_{12}}(R_1^2-|x|^2)\right)^{1/2}&
    {\rm if }&
    |x|\leq R_1\\\left(\frac{R_2^2-|x|^2}{2\alpha_2}\right)^{1/2}& {\rm if
      }&R_1\leq |x|\leq R_2\\0& {\rm if }& |x|\geq R_2\end{array}\right.
\end{eqnarray}
define a solution for $\eps=0$ to the system (\ref{sys0}), which, taking into account (\ref{mu2}) and (\ref{mu1}), can
be rewritten as
\begin{eqnarray}\label{sys}
\left\{\begin{array}{l}\eps^2\Delta\eta_{1}+\left(\frac{\alpha_0}{\alpha_2}(R_2^2-R_1^2)+R_1^2-|x|^2\right)\eta_{1}-2\alpha_1\eta_{1}^3-2\alpha_0\eta_{2}^2\eta_{1}=0\\
\eps^2\Delta\eta_{2}+\left(R_2^2-|x|^2\right)\eta_{2}-2\alpha_2\eta_{2}^3-2\alpha_0\eta_{1}^2\eta_{2}=0.\end{array}\right.
\end{eqnarray}
\begin{rem}
From (\ref{mu2}) and (\ref{mu1}), $\mu_2-\mu_1=\Gamma_2(R_2^2-R_1^2)$.
Thus, under the extra assumption $\Gamma_2>0$ (an assumption which
will be made later), the assumption $R_1<R_2$ implies $\mu_1<\mu_2$.
\end{rem}
\begin{rem}\label{R1R2}
If $R_1=R_2=R$ and $\Gamma_{12}\neq 0$, $\eta_{10}(x)$ and $\eta_{20}(x)$ are given by
(\ref{et10}) and (\ref{et20}) for $|x|\leq R$, and they both vanish at
$|x|=R$. We infer $\Gamma_2/\Gamma_{12}>0$, $\Gamma_1/\Gamma_{12}>0$
(if $\Gamma_1=0$ or $\Gamma_2=0$, then one of the two components is
identically equal to 0, and therefore we are brought back to the study
of one simple equation, like the one which was studuied in \cite{GP}) and
$$
\frac{1}{\Gamma_2}\left(\mu_1-\frac{\alpha_0}{\alpha_2}\mu_2\right)=R^2=\frac{1}{\Gamma_1}\left(\mu_2-\frac{\alpha_0}{\alpha_1}\mu_1\right),$$
which implies $\mu_1=\mu_2=\mu$. Then,  for $\eps>0$, if $\eta$ denotes the ground
state of
$$\eps^2\Delta\eta+(\mu-|x|^2)\eta-2|\Gamma_{12}|\eta^3=0$$
(which, up to a rescaling, is the one which is described in \cite{GP}), then
$$(\eta_1,\eta_2)=\left((|\Gamma_{2}|/\alpha_1)^{1/2}\eta,(|\Gamma_{1}|/\alpha_2)^{1/2}\eta\right)$$
solves (\ref{sys0}).
\end{rem}
\begin{rem}\label{Gamma120}
If $\Gamma_{12}=0$, then an analysis similar to the one which is done above implies $\alpha_1=\alpha_2=\alpha_0=\alpha$ and $\mu_1=\mu_2=\mu$. Then,
$$(\eta_1,\eta_2)=(\eta,\eta)$$
solves (\ref{sys0}), where $\eta$ is the ground state solution of 
$$\eps^2\Delta\eta+(\mu-|x|^2)\eta-4\alpha\eta^3=0,$$
which is described in \cite{GP} (up to a rescaling).
\end{rem}

\subsection{Goal and strategy}
Our goal is to construct a solution $(\eta_{1},\eta_{2})$ of (\ref{sys}) for $\eps>0$ sufficiently small, and to
describe its convergence to $(\eta_{10},\eta_{20})$ as $\eps\to
0$.  The first step consists in constructing aproximate solutions of
(\ref{sys}). Because of the singularities of $\eta_{10}$ and
$\eta_{20}$ at $|x|=R_1$ and $|x|=R_2$, $(\eta_1(x),\eta_2(x))$ will
be described by functions of different variables, depending on the
region of $\R^d$ $x$ belongs to. We write $\R^d=D_0\cup D_1\cup D_2$, where
\begin{eqnarray}
D_0&=&\left\{x\in\R^d\big| |x|^2\leq
R_1^2-\eps^\beta\right\},\nonumber
\end{eqnarray}
\begin{eqnarray}
D_1&=&\left\{x\in\R^d\big| R_1^2-2\eps^\beta\leq|x|^2\leq
R_1^2+2\eps^\beta\right\}\nonumber
\end{eqnarray}
and
\begin{eqnarray}
D_2=\left\{x\in\R^d\big| |x|^2\geq
R_1^2+\eps^\beta\right\},\nonumber
\end{eqnarray}
where $\beta\in (0,2/3)$ is some number that will be fixed later (note that $D_0\cap D_1$ and $D_1\cap D_2$ are not empty). Then, for $x\in D_0$,
$(\eta_1(x),\eta_2(x))$ will be described as a function of the
variable $z=R_1^2-|x|^2$, whereas for $j=1,2$ and $x\in D_j$, it will
be described as a function of the real variables $y_j$ given by
\begin{eqnarray}\label{var}
y_j=\frac{R_j^2-|x|^2}{\eps^{2/3}}.
\end{eqnarray}
In order to be more specific, let us introduce the following
truncation functions. Let $\phi$ be a $\mathcal{C}^\infty$ function on
$\R$ wich is identically equal to 0 on $\R_-$ and identically equal to 1 on
$[1,+\infty)$. Then, let us define
$$\Phi_\eps(z)=\phi\left(\frac{z-\eps^\beta}{2\eps^\beta-\eps^\beta}\right),$$
such that $\Phi_\eps(z)\equiv 0$ for $z\leq \eps^\beta$ and
$\Phi_\eps(z)\equiv 1$ for $z\geq 2\eps^\beta,$
which means (if $\Phi_\eps(z)=\Phi_\eps(R_1^2-|x|^2)$ is considered as a function of the variable $x$, also denoted $\Phi_\eps$ for convenience) that ${\rm Supp} \Phi_\eps\subset D_0$ and
$\Phi_\eps\equiv 1$ for $x\in D_0\string\ D_1$.
Similarly, we set
$$\chi_\eps(y_1)=\left(1-\phi\left(\frac{\eps^{2/3}y_1-\eps^\beta}{2\eps^\beta-\eps^\beta}\right)\right)\phi\left(\frac{\eps^{2/3}y_1+2\eps^\beta}{-\eps^\beta+2\eps^\beta}\right),$$
such that $\chi_\eps(y_1)\equiv 0$ for $y_1\geq
  2\eps^{\beta-2/3}$ and $y_1\leq -2\eps^{\beta-2/3}$,
  whereas $\chi_\eps(y_1)\equiv 1$ for
  $-\eps^{\beta-2/3}\leq y_1\leq \eps^{\beta-2/3}$, which means (if $\chi_\eps(y_1)=\chi_\eps((R_1^2-|x|^2)/\eps^{2/3})$ is considered as a function of the variable $x$, also denoted $\chi_\eps$)  that ${\rm Supp} \chi_\eps\subset D_1$ and
$\chi_\eps\equiv 1$ for $x\in D_1\string\ (D_0\cup D_2)$.
We also define
$$\Psi_\eps(y_2)=1-\phi\left(\frac{z}{\eps^\beta}+2\right)=1-\phi\left(-\frac{R_2^2-R_1^2}{\eps^{\beta}}+\eps^{2/3-\beta}y_2+2\right),$$
such that $\Psi_\eps(y_2)\equiv 0$ for $y_2\geq
\frac{R_2^2-R_1^2}{\eps^{2/3}}-\eps^{\beta-2/3}$ and $\Psi_\eps(y_2)\equiv 1$ for
$y_2\leq\frac{R_2^2-R_1^2}{\eps^{2/3}}-2\eps^{\beta-2/3}$, which means (if $\Psi_\eps(y_2)=\Psi_\eps((R_2^2-|x|^2)/\eps^{2/3})$ is considered as a function of $x$, also denoted $\Psi_\eps$) that ${\rm Supp} \Psi_\eps\subset D_2$ and
$\Psi_\eps\equiv 1$ for $x\in D_2\string\ D_1$.
Formally, we look for $(\eta_{1},\eta_{2})$ under the form 
\begin{eqnarray}\label{ansatznorest}
\left\{\begin{array}{l}\eta_{1}(x)=\Phi_\eps\omega(z)+\eps^{1/3}\chi_\eps\nu(y_1)\\
\eta_{2}(x) =\Phi_\eps\tau(z)+\eps^{1/3}\chi_\eps\lambda(y_1)^{1/2}+\eps^{1/3}\Psi_\eps\mu(y_2),
\end{array}\right.
\end{eqnarray}
in such a way that
\be\label{apD0}
{\rm for}\ x\in D_0,\quad (\eta_1,\eta_2)(x)\approx\left(\omega(z),\tau(z)\right),
\ee
\be\label{apD1}
{\rm for}\ x\in D_1,\quad (\eta_1,\eta_2)(x)\approx\eps^{1/3}\left(\nu(y_1),\lambda(y_1)^{1/2}\right)
\ee
and
\be\label{apD2}
{\rm for}\ x\in D_2,\quad (\eta_1,\eta_2)(x)\approx\left(0,\eps^{1/3}\mu(y_2)\right).
\ee
We look for approximate values of the functions
$\omega,\ \tau,\ \nu,\ \lambda$ and $\mu$ by using a multi-scale
analysis. Namely, we write
\begin{eqnarray}\label{multiscale}
\begin{array}{l}
\omega=\omega_0+\eps^2\omega_1+\eps^{4}\omega_2+\cdots\\
\tau=\tau_0+\eps^2\tau_1+\eps^{4}\tau_2+\cdots\\
\nu=\nu_0+\eps^{2/3}\nu_1+\eps^{4/3}\nu_2+\cdots\\
\lambda=\lambda_{-1}\eps^{-2/3}+\lambda_0+\eps^{2/3}\lambda_1+\eps^{4/3}\lambda_2+\cdots\\
\mu=\mu_0+\eps^{2/3}\mu_1+\eps^{4/3}\mu_2+\cdots
\end{array}
\end{eqnarray}
The $\omega_j$'s, $\tau_j$'s, $\nu_j$'s, $\lambda_j$'s and $\mu_j$'s in this
expansions are $\eps$-independent functions, which are chosen in such
a way that (\ref{apD0}), (\ref{apD1}) and (\ref{apD2}) provide at
least formally solutions to (\ref{sys}) at any order.\\
Then, we prove rigorously that the truncation of the formal asymptotic
expansions we have obtained are indeed approximations of positive
solutions to (\ref{sys}) which converge to $(\eta_{10},\eta_{20})$ as
$\eps\to 0$. For this purpose, we use the ansatz
\begin{eqnarray}\label{ansatz}
\left\{\begin{array}{l}\eta_{1}(x)=\Phi_\eps\omega(z)+\eps^{1/3}\left(\chi_\eps\nu(y_1)+\eps^{2(N+1)/3}P(x)\right)\\
\eta_{2}(x)=\Phi_\eps\tau(z)+\eps^{1/3}\left(\chi_\eps\lambda(y_1)^{1/2}+\Psi_\eps\mu(y_2)+\eps^{2(N+1)/3}Q(x)\right).
\end{array}\right.
\end{eqnarray}
where $\omega$, $\tau$, $\nu$, $\lambda$  and $\mu$  are now
truncations up to some finite order ($N\in\mathbb{N}$ for $\nu$,
$\lambda$, $M=M(N)$ for
$\omega$, $\tau$ and $L=L(N)$ for $\mu$) of the formal series (\ref{multiscale}), and $P$, $Q$ are remainder terms. A fixed point theorem provides the existence of $P,Q$ as
well as estimates which ensure that the remainder terms in (\ref{ansatz}) are indeed
small. The better $\omega$, $\tau$, $\nu$, $\lambda$  and $\mu$ are
chosen (that is, the larger is N), the smaller is
$\eps^{2(N+1)/3}(P,Q)$. The functional space in which $(P,Q)$ is obtained is 
 $H^1_w(\R^d)^2$, where 
 $$H^1_w(\R^d)=\left\{f\in H^1(\R^d)\ |\  \min(|y_1|,|y_2|)^{1/2}f\in L^2(\R^d)\right\}.$$
$H^1_w(\R^d)^2$ is endowed with the norm
$$\left\|(P,Q)\right\|_{H^1_w(\R^d)^2}=\left(\int_{\R^d}\left(|\nabla
  P|^2+|\nabla Q|^2\right)dx+\int_{\R^d}\max(1,\min(|y_1|,|y_2|))(|P|^2+|Q|^2)dx\right)^{1/2}.$$
  
\begin{rem}
Note that the set  $H^1_w(\R^d)^2$ does not depend on $\eps$, even though it's norm does. However, this norm has been chosen in such a way that the norm of the continuous embedding of $H^1_w(\R^d)^2$ into $H^1(\R^d)^2$ is uniformly bounded in $\eps$.
\end{rem}
 
Once we have constructed $(\eta_1,\eta_2)$, we would like to estimate
in different norms the difference between the exact solution
$(\eta_1,\eta_2)$ and its approximation
\begin{eqnarray}\label{approximation}
\left\{\begin{array}{l}\eta_{1app}(x)=\Phi_\eps\omega(z)+\eps^{1/3}\chi_\eps\nu(y_1)\\
\eta_{2app}(x)=\Phi_\eps\tau(z)+\eps^{1/3}\chi_\eps\lambda(y_1)^{1/2}+\eps^{1/3}\Psi_\eps\mu(y_2),
\end{array}\right.
\end{eqnarray}
where $\omega$, $\tau$, $\nu$, $\lambda$  and $\mu$  are truncations
of the formal series (\ref{multiscale}) up to some fixed orders $M_0$
for $\omega$ and $\tau$, $N_0$ for $\nu$ and $\lambda$ and $L_0$ for $\mu$. However, the estimates on
$P$ and $Q$ provided by the fixed point argument are not very good. So
in order to get better estimates on $\eta_j-\eta_{japp}$ ($j=1,2$), we
proceed as follows. We choose three large integers $M>M_0$, $N>N_0$
and $L>L_0$ and write $\eta_j$ as in (\ref{ansatz}) (with truncations
of the formal power series at orders $M$, $N$ and $L$ instead of
$M_0$, $N_0$ and $L_0$). Then, the estimate on $\eta_j-\eta_{japp}$ is
obtained thanks to estimates on $\eps^{2m}\omega_m$ and
$\eps^{2m}\tau_m$ for $M_0+1\leq m\leq M$, on $\eps^{2n/3}\nu_n$ and
$\eps^{2n/3}\lambda_n$ for $N_0+1\leq n\leq N$ and on
$\eps^{2n/3}\mu_n$ for $L_0+1\leq n\leq L$. The estimates on $P$ and
$Q$ provided by the fixed point argument are good enough to ensure
that $\eps^{2N/3+1}P$ and $\eps^{2N/3+1}Q$ are negligible in
comparison with the other terms in the expression of
$\eta_j-\eta_{japp}$.

In our main result below, we give estimates on the $L^p(\R^d)$ and
$H^1(\R^d)$ norms of $\eta_j-\eta_{japp}$ for $j=1,2$. Note however
that depending on the need of the reader, our strategy can give many
other informations on  $\eta_j-\eta_{japp}$ (see Remark
\ref{rem-post-theo} below).

\begin{theo}\label{main}
Let $d\in \{1,2,3\}$, and $\alpha_0,\alpha_1,\alpha_2>0$, $R_2>R_1>0$
such that $\Gamma_2,\Gamma_{12}>0$ and such that (\ref{cond-disk}) is
satisfied. Then, for $\eps>0$ sufficiently small, (\ref{sys}) has a
unique solution $(\eta_1,\eta_2)\in\mathcal{C}_0(\R^d)^2$ such that
the two components $\eta_1$ and $\eta_2$ are both positive. Moreover,
if $M_0,N_0,L_0\in\N$, if $\beta\in(0,\eps^{2/3})\backslash \Q$, if $\omega_m=\omega_m(z),\tau_m=\tau_m(z),\nu_n=\nu_n(y_1),\lambda_n=\lambda_n(y_1),\mu_n=\mu_n(y_2)$ are the functions given by (\ref{eq-omega0}), (\ref{eq-tau0}), (\ref{rec-omegan-taun}), (\ref{mun}), (\ref{I-1}), (\ref{nu00}), (\ref{II0}), (\ref{In}) and (\ref{IIn}), then
\be\label{est-main}
\left\|\eta_1-\eta_{1app}\right\|_{E}=\mathcal{O}\left(\eps^{\gamma_1(E)}\right)
\quad
{\rm and}
\quad
\left\|\eta_2-\eta_{2app}\right\|_{E}=\mathcal{O}\left(\eps^{\gamma_2(E)}\right),
\ee
where $E$ can be either $L^p(\R^d)$ for any $p\in [2,+\infty]$ or $H^1(\R^d)$,
$$\eta_{1app}=\Phi_\eps\sum_{m=0}^{M_0}\eps^{2m}\omega_m
+\eps^{1/3}\chi_\eps\sum_{n=0}^{N_0}\eps^{2n/3}\nu_n,$$
$$\eta_{2app}=\Phi_\eps\sum_{m=0}^{M_0}\eps^{2m}\tau_m +\eps^{1/3}\chi_\eps\left(\sum_{n=-1}^{N_0}\eps^{2n/3}\lambda_n\right)^{1/2} +\eps^{1/3}\Psi_\eps\sum_{n=0}^{L_0}\eps^{2n/3}\mu_n$$
and for $p\in [2,+\infty]$, 
$$\gamma_1(L^p(\R^d))=\left\{\begin{array}{ll}
\min\big((2-3\beta)M_0+2-5\beta/2+\beta/p\ ,\ 1+2/(3p)\big) & {\rm if}\ N_0=0\\
\min\big((2-3\beta)M_0+2-5\beta/2+\beta/p\ ,\ 5/3+2/(3p)\big) & {\rm
  if}\ N_0=1 \ {\rm and}\ p>2\\
\min\big((2-3\beta)M_0+2-5\beta/2+\beta/p\ ,\ 2-\delta\big) & {\rm
  if}\ N_0=1 \ {\rm and}\ p=2\\
\min\big((2-3\beta)M_0+2-5\beta/2+\beta/p\ ,\ \beta N_0+2-3\beta/2+\beta/p\big)& {\rm
  if}\ N_0\geq 2\end{array}\right.$$
where $\delta>0$ is arbitrarily small, and
$$\gamma_2(L^p(\R^d))=\left\{\begin{array}{ll}
\min\big((2-3\beta)M_0+2-2\beta+\beta/p\ ,\ 4/3+2/(3p)\ ,\ 2L_0/3+1+2/(3p)\big)&
    {\rm if\ } N_0=0\\
\min\big((2-3\beta)M_0+2-2\beta+\beta/p\ ,\ \beta N_0+2-\beta+\beta/p\ ,\ 2L_0/3+1+2/(3p)\big)&  {\rm if\ } N_0\geq 1,\end{array}\right.$$
whereas
$$\gamma_1(H^1(\R^d))=\left\{\begin{array}{ll}
\min\big((2-3\beta)(M_0+1)\ ,\ 2/3\big)&
        {\rm if\ } N_0=0\\
\min\big((2-3\beta)(M_0+1)\ ,\ 4/3\big)&
        {\rm if\ } N_0=1\\
\min\big((2-3\beta)(M_0+1)\ ,\ 2-\delta\big)&
        {\rm if\ } N_0=2\\
\min\big((2-3\beta)(M_0+1)\ ,\ \beta(N_0-2)+2\big)&
        {\rm if\ } N_0\geq 3
\end{array}\right.$$
and 
$$
\gamma_2(H^1(\R^d))=\left\{\begin{array}{ll}
\min\big((2-3\beta)(M_0+1)+\beta/2\ ,\ 1\ ,\ 2L_0/3+2/3\big)&
        {\rm if\ } N_0=0\\
\min\big((2-3\beta)(M_0+1)+\beta/2\ ,\ 5/3\ ,\ 2L_0/3+2/3\big)&
        {\rm if\ } N_0=1\\
\min\big((2-3\beta)(M_0+1)+\beta/2\ ,\ \beta(N_0-3/2)+2\ ,\ 2L_0/3+2/3\big)&
        {\rm if\ } N_0\geq 2.
\end{array}\right.$$
\end{theo}

\begin{rem}\label{rem-post-theo}
Depending on the value of $M_0$, $N_0$, $L_0$ and $p$, the value of the
parameter $\beta\in (0,2/3)$ can be adjusted in such a way that the
values of $\gamma_1$ and $\gamma_2$ are as large as possible. If we
are only interested in the approximation of one of the two components
$\eta_j$, one can
even choose $\beta\in(0,2/3)$ to optimize $\gamma_j$ without
considering the other component. In some cases, one can be
interested in estimations on the norms of $\eta_1-\eta_{1app}$ and
$\eta_2-\eta_{2app}$, not on $\R^d$ as a whole, but only on a
subdomain like $D_0$, $D_1$ or $D_2$. In each minimum in the
expressions of $\gamma_1$ and $\gamma_2$ in the statement of the
theorem, the first argument 
corresponds to the rate of convergence of the norm in $D_0$, the
second one to the rate of convergence of the norm in $D_1$ and the
third one (for $\eta_2$) to the rate of convergence of the norm in
$D_2$. The $L^p$ and $H^1$ norms of the restriction of
$\eta_1-\eta_{1app}$ to $D_2$ converge to 0 faster than any power of
$\eps$ as $\eps\to 0$.
\end{rem}

In the following corollary, we write more expicitely upper bounds on
the rates of convergence of $\eta_1-\eta_{1app}$ and
$\eta_2-\eta_{2app}$ to 0 in the particular 
and important case where $M_0=N_0=L_0=0$ and $E=L^2,L^\infty{\rm\ or\ }H^1$.

\begin{cor}If $\beta\in (0,1/3)$, we have
$$\eta_1=\Phi_\eps\omega_0+\eps^{1/3}\chi_\eps\nu_0+\left\{\begin{array}{l}\mathcal{O}_{L^2(\R^d)}(\eps^{4/3})\\ \mathcal{O}_{L^\infty(\R^d)}(\eps)\\\mathcal{O}_{H^1(\R^d)}(\eps^{2/3})\end{array}\right.$$
and
$$\eta_2=\Phi_\eps\tau_0+\eps^{1/3}\chi_\eps\left(\frac{\lambda_{-1}}{\eps^{2/3}}+\lambda_0\right)+\eps^{1/3}\Psi_\eps\mu_0+\left\{\begin{array}{l}\mathcal{O}_{L^2(\R^d)}(\eps^{4/3})\\ \mathcal{O}_{L^\infty(\R^d)}(\eps)\\\mathcal{O}_{H^1(\R^d)}(\eps^{2/3}).\end{array}\right.$$
\end{cor}


\subsection{Organization of the paper}
In Section \ref{formal}, we calculate formally all the functions
$\omega_j$'s, $\tau_j$'s, $\nu_j$'s, $\lambda_j$'s and $\mu_j$'s
appearing in the formal series (\ref{multiscale}), in such a way that
truncations of these series provide at least formally, through the
ansatz (\ref{apD0}), (\ref{apD1}) and (\ref{apD2}), solutions to
(\ref{sys}) at any order. We also study  asymptotic behaviours of
these functions. In Section \ref{truncation}, we study the
functions obtained by truncations of the formal series. In particular,
if $\omega$, $\tau$, $\nu$, $\lambda$, $\mu$ denote these truncations,
we estimate the order at which (\ref{apD0}), (\ref{apD1}) and
(\ref{apD2}) solve (\ref{sys}), respectively on $D_0$, $D_1$ and
$D_2$. We also check that (\ref{apD0}) and (\ref{apD1}) are close one
from another on $D_0\cap D_1$ and that (\ref{apD1}) and (\ref{apD2})
are close one from another on $D_1\cap D_2$. Section \ref{rigorous} is
devoted to the proof of the main result.

\paragraph{Notations.}
\begin{itemize}
\item
If $A$ and $B$ are two quantities depending on a parameter $x$
belonging to some set $D$, the claim ``for $x\in D$, $A(x) \lesssim
B(x)$'' means ``there exists $C>0$ such that for
every $x\in D$, $A(x) \leq C B(x)$''.
\item
Let $F(x)$ be a function defined in a neighborhood of $\infty$.
Given $\alpha \in \R$, $\{ f_m \}_{m \in \N} \in \R$, and $\gamma > 0$, the
notation
\be\label{defasympseries}
F(x) \underset{x \to \infty}{\approx} x^\alpha\sum_{m=0}^\infty f_m x^{-\gamma m}\nonumber
\ee
means that for every $M\in\N$,
$$
F(x)-x^\alpha\sum_{m=0}^M f_m x^{-\gamma
  m}=\mathcal{O}(x^{\alpha-\gamma(M+1)})\quad \text{as }x\to\infty,
$$
and, moreover, that the asymptotic series can be differentiated term
by term. We use the same notation if $\gamma<0$ and if $F$ is defined
in a neighborhood of $0$.
\item
$\mathcal{C}_0(\R^d) $ denotes the space of continuous functions on $\R^d$ that converge to 0 at infinity.
\item
If $(f_\eps)_{0<\eps<\eps_0}$ is a sequence of functions such that for every $\eps$, $f_\eps$ belongs to some Banach space $E_\eps$ that may depend on $\eps$, if $\alpha\in \R$, $f_\eps=\mathcal{O}_{E_\eps}(\eps^\alpha)$  (respectively  $f_\eps=o_{E_\eps}(\eps^\alpha)$)  means that $\|f_\eps\|_{E_\eps}/\eps^\alpha$ remains bounded (respectively converges to 0) as $\eps\to 0$.
\end{itemize}

\section{Formal asymptotic expansions}\label{formal}
\subsection{Asymptotic behaviour of
  $\nu_0,\mu_0,\lambda_{-1},\lambda_0$.}\label{asymp}
We are looking for a solution $(\eta_{1},\eta_{2})$ to
(\ref{sys}) which converges to the Thomas-Fermi approximation as
$\eps\to 0$. Namely, for every $x\in\R^d$,
\begin{eqnarray}\label{limeta}
\eta_{1}(x)\underset{\eps\to
  0}{\longrightarrow}\eta_{10}(x),\quad
\eta_{2}(x)\underset{\eps\to
  0}{\longrightarrow}\eta_{20}(x).
\end{eqnarray}
The convergence of $(\eta_{1},\eta_{2})$ (expressed using the
ansatz (\ref{ansatz})) to the Thomas-Fermi
limit determines the asymptotic behaviour of
$\nu(y_1),\mu(y_2),\lambda(y_1)$ as
$y_1,y_2\to\pm\infty$. We will construct the functions
$\nu_0,\mu_0,\lambda_{-1}$ and $\lambda_0$ in such a way that they capture entirely this
asymptotic behaviour. More precisely,
$$\begin{array}{llll}{\rm for}\ |x|>R_2,&\eps^{1/3}\mu_0(y_2)\underset{\eps\to
  0}{\longrightarrow}0 & {\rm
  yields} & \mu_0(y_2)\underset{y_2\to-\infty}{\longrightarrow}0,\\
{\rm for}\ R_1<|x|<R_2, & \eps^{1/3}\mu_0(y_2)\underset{\eps\to
  0}{\longrightarrow}\left(\frac{R_2^2-|x|^2}{2\alpha_2}\right)^{1/2} &{\rm
  yields} &
\mu_0(y_2)\underset{y_2\to+\infty}{\sim}\left(\frac{y_2}{2\alpha_2}\right)^{1/2},\\
{\rm for}\ R_1<|x|<R_2,&\eps^{1/3}\nu_0(y_1)\underset{\eps\to
  0}{\longrightarrow}0 &{\rm
  yields}&\nu_0(y_1)\underset{y_1\to-\infty}{\longrightarrow}0,\\
{\rm for}\ |x|<R_1,&\eps^{1/3}\nu_0(y_1)\underset{\eps\to
  0}{\longrightarrow}\left(\frac{\Gamma_2}{2\alpha_1\Gamma_{12}}(R_1^2-|x|^2)\right)^{1/2}&{\rm
  yields}& \nu_0(y_1)\underset{y_1\to+\infty}{\sim}\left(\frac{\Gamma_2y_1}{2\alpha_1\Gamma_{12}}\right)^{1/2}\\
{\rm for}\ R_1<|x|<R_2,&\multicolumn{3}{l}{\eps^{1/3}\left(\frac{\lambda_{-1}(y_1)}{\eps^{2/3}}+\lambda_0(y_1)\right)^{1/2}\underset{\eps\to
  0}{\longrightarrow}\left(\frac{R_2^2-|x|^2}{2\alpha_2}\right)^{1/2}=\left(\frac{R_2^2-R_1^2}{2\alpha_2}+\frac{R_1^2-|x|^2}{2\alpha_2}\right)^{1/2}}\\
&& {\rm yields}&
\lambda_{-1}(y_1)\underset{y_1\to-\infty}{\longrightarrow}\frac{R_2^2-R_1^2}{2\alpha_2},\\
&&&\lambda_0(y_1)\underset{y_1\to-\infty}{\sim}\frac{y_1}{2\alpha_2},\\
{\rm for}\ |x|<R_1,&
\multicolumn{3}{l}{\eps^{1/3}\left(\frac{\lambda_{-1}(y_1)}{\eps^{2/3}}+\lambda_0(y_1)\right)^{1/2}\underset{\eps\to
  0}{\longrightarrow}\left(\frac{R_2^2-R_1^2}{2\alpha_2}+\frac{\Gamma_1}{2\alpha_2\Gamma_{12}}(R_1^2-|x|^2)\right)^{1/2}}\\
&& {\rm
  yields}&
\lambda_{-1}(y_1)\underset{y_1\to+\infty}{\longrightarrow}\frac{R_2^2-R_1^2}{2\alpha_2},\\
&&&\lambda_0(y_1)\underset{y_1\to+\infty}{\sim}\frac{\Gamma_1y_1}{2\alpha_2\Gamma_{12}}.
\end{array}$$

\subsection{Expansions of $\omega$ and $\tau$ in $D_0$}\label{secomtau}
In the domain $D_0$, we look for $(\eta_1,\eta_2)$ solution of (\ref{sys}) under the form
(\ref{apD0}). It follows that $\omega(z)$ and $\tau(z)$ have to solve
for $z\in (0,R_1^2)$ the following system of differential equations
\be\label{eq-omega}
-2d\eps^2\omega'+4(R_1^2-z)\eps^2\omega''+\left(\frac{\alpha_0}{\alpha_2}(R_2^2-R_1^2)+z\right)\omega-2\alpha_1\omega^3-2\alpha_0\tau^2\omega
&=&0\\
\label{eq-tau}
-2d\eps^2\tau'+4(R_1^2-z)\eps^2\tau''+\left(R_2^2-R_1^2+z\right)\tau-2\alpha_2\tau^3-2\alpha_0\omega^2\tau&=&0.
\ee
Then, we look for $\omega$ and $\tau$ under the form of formal power
series in the parameter $\eps^2$:
$$\omega=\sum_{m=0}^\infty\eps^{2m}\omega_m,\quad
\tau=\sum_{m=0}^\infty\eps^{2m}\tau_m.$$
Plugging these expansions into (\ref{eq-omega}), we get
\be\label{eq-I-omegan-taun}
-2d\sum_{m=1}^\infty\eps^{2m}\omega_{m-1}'+4(R_1^2-z)\sum_{m=1}^\infty\eps^{2m}\omega_{m-1}''+\left(\frac{\alpha_0}{\alpha_2}(R_2^2-R_1^2)+z\right)\sum_{m=0}^\infty\eps^{2m}\omega_m\nonumber\\
-2\alpha_1\sum_{m=0}^\infty\eps^{2m}\sum_{m_1+m_2+m_3=m}\omega_{m_1}\omega_{m_2}\omega_{m_3}-2\alpha_0\sum_{m=0}^\infty\eps^{2m}\sum_{m_1+m_2+m_3=m}\omega_{m_1}\tau_{m_2}\tau_{m_3}&=&0,
\ee
whereas (\ref{eq-tau}) yields
\be\label{eq-II-omegan-taun}
-2d\sum_{m=1}^\infty\eps^{2m}\tau_{m-1}'+4(R_1^2-z)\sum_{m=1}^\infty\eps^{2m}\tau_{m-1}''+\left(R_2^2-R_1^2+z\right)\sum_{m=0}^\infty\eps^{2m}\tau_m\nonumber\\
-2\alpha_2\sum_{m=0}^\infty\eps^{2m}\sum_{m_1+m_2+m_3=m}\tau_{m_1}\tau_{m_2}\tau_{m_3}
-2\alpha_0\sum_{m=0}^\infty\eps^{2m}\sum_{m_1+m_2+m_3=m}\omega_{m_1}\omega_{m_2}\tau_{m_3}=0.
\ee
At order $m=0$, we deduce that $\omega_0^2$, $\tau_0^2$ have to solve
in the domain $z\in(0,R_1^2)$ (a range of values of $z$ for which
they are expected not to vanish) the linear system
\be
\alpha_1\omega_0^2+\alpha_0\tau_0^2 &=
&\frac{1}{2}\left(\frac{\alpha_0}{\alpha_2}(R_2^2-R_1^2)+z\right)\nonumber\\
\alpha_0\omega_0^2+\alpha_2\tau_0^2 &=
&\frac{1}{2}\left(R_2^2-R_1^2+z\right).\nonumber
\ee
As already mentioned in (\ref{eta10}) and (\ref{eta20}), it follows
that 
\be\label{eq-omega0}
\omega_0^2&=&\frac{\Gamma_2}{2\alpha_1\Gamma_{12}}z
\ee
and
\be\label{eq-tau0}
\tau_0^2=\frac{R_2^2-R_1^2}{2\alpha_2}+\frac{\Gamma_1}{2\alpha_2\Gamma_{12}}z.
\ee
For $m\geq 1$, (\ref{eq-I-omegan-taun}) and (\ref{eq-II-omegan-taun})
imply that
\be\label{47}
M\left[\begin{array}{c}\omega_m\\ \tau_m\end{array}\right]=\left[\begin{array}{c}2d\omega_{m-1}'+4(z-R_1^2)\omega_{m-1}''+2\alpha_1\!\!\!\!\!\underset{\tiny\begin{array}{c}m_1+m_2+m_3=m\\m_1,m_2,m_3<m\end{array}}{\sum}\!\!\!\!\!\omega_{m_1}\omega_{m_2}\omega_{m_3}+2\alpha_0\!\!\!\!\!\underset{\tiny\begin{array}{c}m_1+m_2+m_3=m\\m_1,m_2,m_3<m\end{array}}{\sum}\!\!\!\!\!\omega_{m_1}\tau_{m_2}\tau_{m_3}\\
2d\tau_{m-1}'+4(z-R_1^2)\tau_{m-1}''+2\alpha_2\!\!\!\!\!\underset{\tiny\begin{array}{c}m_1+m_2+m_3=m\\m_1,m_2,m_3<m\end{array}}{\sum}\!\!\!\!\!\tau_{m_1}\tau_{m_2}\tau_{m_3}+2\alpha_0\!\!\!\!\!\underset{\tiny\begin{array}{c}m_1+m_2+m_3=m\\m_1,m_2,m_3<m\end{array}}{\sum}\!\!\!\!\!\omega_{m_1}\omega_{m_2}\tau_{m_3}\end{array}\right],
\ee
where
\be\label{def-M}
M&=&-4\left[\begin{array}{cc}\alpha_1\omega_0^2 &
      \alpha_0\omega_0\tau_0\\
\alpha_0\omega_0\tau_0 &\alpha_2\tau_0^2\end{array}\right].
\ee
Thus, the functions $\omega_m$, $\tau_m$ for $m\geq 1$ can be
calculated thanks to the recursion relation
\be\label{rec-omegan-taun}
\lefteqn{\left[\begin{array}{c}\omega_m\\ \tau_m\end{array}\right]=\frac{1}{\alpha_1\alpha_2\Gamma_{12}\omega_0^2\tau_0^2}\left[\begin{array}{cc} \alpha_2\tau_0^2&
      -\alpha_0\omega_0\tau_0\\
-\alpha_0\omega_0\tau_0 &\alpha_1\omega_0^2\end{array}\right]\times}\\
&&\left[\begin{array}{c}-\frac{d}{2}\omega_{m-1}'-(z-R_1^2)\omega_{m-1}''-\frac{\alpha_1}{2}\underset{\tiny\begin{array}{c}m_1+m_2+m_3=m\\m_1,m_2,m_3<m\end{array}}{\sum}\omega_{m_1}\omega_{m_2}\omega_{m_3}-\frac{\alpha_0}{2}\underset{\tiny\begin{array}{c}m_1+m_2+m_3=m\\m_1,m_2,m_3<m\end{array}}{\sum}\omega_{m_1}\tau_{m_2}\tau_{m_3}\\
-\frac{d}{2}\tau_{m-1}'-(z-R_1^2)\tau_{m-1}''-\frac{\alpha_2}{2}\underset{\tiny\begin{array}{c}m_1+m_2+m_3=m\\m_1,m_2,m_3<m\end{array}}{\sum}\tau_{m_1}\tau_{m_2}\tau_{m_3}-\frac{\alpha_0}{2}\underset{\tiny\begin{array}{c}m_1+m_2+m_3=m\\m_1,m_2,m_3<m\end{array}}{\sum}\omega_{m_1}\omega_{m_2}\tau_{m_3}\end{array}\right].\nonumber
\ee
From this relation, we deduce useful informations about the behaviour of $\omega_m$ and
$\tau_m$ for $z\in (0,R_1^2]$. 
\begin{lem}\label{est-omega-tau}
For every $m\geq 1$, there exists $(w_{m,n})_{n\geq 0},(t_{m,n})_{n\geq 0}\in\R^{\N}$ such that
\be\label{asymp-omegan}
\omega_m(z)\underset{z\to 0}{\approx}z^{1/2-3m}\sum_{n=0}^\infty
w_{m,n}z^n
\ee
and
\be\label{asymp-taun}
\tau_m(z)\underset{z\to 0}{\approx}z^{1-3m}\sum_{n=0}^\infty
t_{m,n}z^n.
\ee
In particular, there is a constant $c_m>0$ such that
$$\forall z\in (0,R_1^2], \quad |\omega_m(z)|\leq c_m
z^{1/2-3m}\quad {\rm and}\quad |\tau_m(z)|\leq c_mz^{1-3m}.$$
\end{lem}

\begin{rem}\label{rem-22}
Note that for $m=0$, (\ref{asymp-omegan}) is also true (with
$w_{0,n}=0$ for $n\geq 1$), whereas (\ref{asymp-taun}) has to be
replaced by the Taylor expansion of $\tau_0$ at $z=0$, which can be
written as
\be\label{asymp-tau0}
\tau_0(z)\underset{z\to 0}{\approx}\left(\frac{R_2^2-R_1^2}{2\alpha_2}\right)^{1/2}+\sum_{n=0}^\infty
t_{0,n}z^{1+n}
\ee
for some $(t_{0,n})_{n\geq 0}\in \R^\N$.
\end{rem}

\begin{Proof}
From (\ref{rec-omegan-taun}), $\omega_1$ and $\tau_1$ can be explicitely expressed by
\be\label{eq-omega1-tau1}
\left[\begin{array}{c}\omega_1\\\tau_1\end{array}\right]=\left[\begin{array}{c}\frac{-d\omega_0'/2-(z-R_1^2)\omega_0''}{\alpha_1\Gamma_{12}\omega_0^2}+\frac{\alpha_0(d\tau_0'/2+(z-R_1^2)\tau_0'')}{\alpha_1\alpha_2\Gamma_{12}\omega_0\tau_0}\\\frac{\alpha_0(d\omega_0'/2+(z-R_1^2)\omega_0'')}{\alpha_1\alpha_2\Gamma_{12}\omega_0\tau_0}-\frac{d\tau_0'/2+(z-R_1^2)\tau_0''}{\alpha_2\Gamma_{12}\tau_0^2}\end{array}\right].
\ee
Then, it follows from (\ref{eq-omega0}), (\ref{eq-tau0}) and the
expansions as $\eps\to 0$ of $\tau_0'$, $\tau_0''$, $1/\tau_0$ and
$1/\tau_0^2$ that (\ref{asymp-omegan})-(\ref{asymp-taun}) hold for
$m=1$. 
Let $m\geq 2$ and assume that (\ref{asymp-omegan}) and (\ref{asymp-taun}) are true
for $m$ replaced by any integer between 1 and $m-1$. Then
(\ref{asymp-omegan})-(\ref{asymp-taun}) also hold at order $m$ thanks
to (\ref{rec-omegan-taun}), the recursion assumption,
(\ref{eq-omega0}) and (\ref{asymp-tau0}).
\end{Proof}

\begin{rem}\label{23} A consequence of Lemma \ref{est-omega-tau} is that for
  every $x\in D_0$ (which in terms of the variable $z$, means
  $R_1^2\geq z\geq\eps^\beta$), for every $m\geq 1$, 
$$|\eps^{2m}\omega_m(z)|\leq c_m\eps^{\beta/2+m(2-3\beta)},\quad |\eps^{2m}\tau_m(z)|\leq c_m\eps^{\beta+m(2-3\beta)}.$$
In particular, since we have chosen $\beta\in (0,2/3)$, for every
$M\geq 1$, 
$$\left\|\sum_{m=0}^M\eps^{2m}\omega_m-\omega_0\right\|_{L^\infty(D_0)}\underset{\eps\to
  0}{\longrightarrow} 0\quad{\rm and}\quad \left\|\sum_{m=0}^M\eps^{2m}\tau_m-\tau_0\right\|_{L^\infty(D_0)}\underset{\eps\to
  0}{\longrightarrow} 0,$$
and for a fixed value of $M$, the larger is $m\in\{0,\cdots,M\}$, the
smaller are the $L^\infty(D_0)$ norms of $\eps^{2m}\omega_m$ and
$\eps^{2m}\tau_m$ in the limit $\eps\to 0$.
\end{rem}

\subsection{Expansion of $\mu$ in $D_2$}\label{secmu}
For $x\in D_2$, we look for a solution $(\eta_1,\eta_2)$ to (\ref{sys}) under the form (\ref{apD2}). Thus, $\mu$ is
constructed in such a way that
$\eta_{2}(x)=\eps^{1/3}\mu(y_2)$ solves, for $|x|>R_1$,
\begin{eqnarray}\label{eq-tildeeta}
\eps^2\Delta\eta_{2}+\left(R_2^2-|x|^2\right)\eta_{2}-2\alpha_2\eta_{2}^3=0,
\end{eqnarray}
which means that for $y_2<(R_2^2-R_1^2)/\eps^{2/3}$,
\begin{eqnarray}\label{eq-mu}
4|x|^2\mu''(y_2)-2d\eps^{2/3}\mu'(y_2)+y_2\mu(y_2)-2\alpha_2\mu(y_2)^3=0.
\end{eqnarray}
Moreover, we are looking for a solution $\eta_2$ that converges to
$\eta_{20}$ for $|x|>R_1$. Thus, as already discussed in Section \ref{asymp}, $\mu$ has to satisfy the following
asymptotics:
$$
 \mu(y_2)\underset{y_2\to-\infty}{\longrightarrow}0,\quad\quad
\mu(y_2)\underset{y_2\to+\infty}{\sim}\left(\frac{y_2}{2\alpha_2}\right)^{1/2}.
$$
We rescale to change the unknown function $\mu$ into
$\gamma$, defined by
$$\mu(y_2)=\frac{R_2^{1/3}}{(2\alpha_2)^{1/2}}\gamma\left(\frac{y_2}{R_2^{2/3}}\right).$$
Then, it turns out that $\mu$ solves (\ref{eq-mu}) if and only if $\gamma$ solves the differential equation
\be\label{eq-gamma}
4(1-\tilde{\eps}^{2/3}y)\gamma''(y)-2d\tilde{\eps}^{2/3}\gamma'(y)+y\gamma(y)-\gamma(y)^3=0,\quad
-\infty<y\leq \frac{R_2^2-R_1^2}{R_2^2}\tilde{\eps}^{-2/3},
\ee
where $\tilde{\eps}=\eps/R_2^2$. In \cite{GP}, we have constructed a
solution $\gamma$ of this equation for $y\in
(-\infty,\tilde{\eps}^{-2/3}]$ ($\gamma$ was denoted
$\nu_{\tilde{\eps}}$ in that paper). Moreover, this solution, for any
  $N\in\N$, can be expressed
under the form (see below for an explanation of the notations)
$$\gamma(y)=\sum_{n=0}^N\tilde{\eps}^{2n/3}\gamma_n(y)+\tilde{\eps}^{2(N+1)/3}R_{N,\tilde{\eps}}(y).$$
Thus, 
\begin{eqnarray}\label{defmu}
\mu(y_2)=\frac{R_2^{1/3}}{(2\alpha_2)^{1/2}}\sum_{n=0}^N\tilde{\eps}^{2n/3}\gamma_n\left(\frac{y_2}{R_2^{2/3}}\right)+\tilde{\eps}^{2(N+1)/3}\frac{R_2^{1/3}}{(2\alpha_2)^{1/2}}R_{N,\tilde{\eps}}\left(\frac{y_2}{R_2^{2/3}}\right).
\end{eqnarray}
In particular, the functions $\mu_n$ introduced in
(\ref{multiscale}) are given for every $n\geq 0$ by 
\begin{eqnarray}\label{mun}
\mu_n(y_2)=\frac{R_2^{1/3}}{(2\alpha_2)^{1/2}}R_2^{-4n/3}\gamma_n\left(\frac{y_2}{R_2^{2/3}}\right).
\end{eqnarray}
The functions $\gamma_n$ and $R_{N,\eps}$ mentioned above have been defined as follows in \cite{GP}.
\begin{itemize}
\item $\gamma_0$ is the Hastings-McLeod solution of the Painlev\'e-II
  equation, that is the unique solution of
\begin{eqnarray}
\label{P2-eq}
4 \gamma_0''(y) + y \gamma_0(y) - \gamma_0(y)^3 = 0, \quad y \in \mathbb{R},
\end{eqnarray}
with the asymptotic behaviour 
$$\gamma_0(y)\underset{y\to +\infty}{\sim}y^{1/2},\quad
\gamma_0(y)\underset{y\to -\infty}{\longrightarrow} 0.$$
\item for $1\leq n\leq N$, $\gamma_n$ is the unique solution of
\begin{eqnarray}\label{nun-eq}
-4\gamma_n''(y) + W_0(y)\gamma_n(y) =F_n(y),\quad y\in \R,
\end{eqnarray}
which goes to 0 as $y\to\pm\infty$, where
\begin{eqnarray}\label{W00}
W_0(y) = 3 \gamma_0^2(y) - y
\end{eqnarray}
and
$$
F_n(y)=-\!\!\!\!\!\!\!\!\underset{\tiny{\begin{array}{c}n_1,n_2,n_3<n\\n_1+n_2+n_3=n\end{array}}}{\sum}\!\!\!\!\gamma_{n_1}(y)\gamma_{n_2}(y)\gamma_{n_3}(y)-2d\gamma_{n-1}'(y)-4y\gamma_{n-1}''(y),$$
\item $R_{N,\tilde{\eps}}$ solves 
\begin{eqnarray} \label{r-eq}
-4(1-\tilde{\eps}^{2/3}y)R_{N,\tilde{\eps}}''+2\tilde{\eps}^{2/3}dR_{N,\tilde{\eps}}' + W_0 R_{N,\tilde{\eps}}=F_{N,\tilde{\eps}}(y,R_{N,\tilde{\eps}}),
\ y \in (-\infty,\tilde{\eps}^{-2/3}],
\end{eqnarray}
where 
\begin{eqnarray}
\lefteqn{F_{N,\tilde{\eps}}(y,R)\ =\
-(4y\nu_N''+2d\nu_N')-\sum_{n=0}^{2N-1}\tilde{\eps}^{2n/3}\underset{\tiny{\begin{array}{c}n_1+n_2+n_3=n+N+1\\0\leq
      n_1,n_2,n_3\leq
      N\end{array}}}{\sum}\gamma_{n_1}\gamma_{n_2}\gamma_{n_3}}\nonumber\\
&&-\left(3\sum_{n=1}^{2N}\tilde{\eps}^{2n/3}\underset{\tiny{\begin{array}{c}n_1+n_2=n\\0\leq
      n_1,n_2\leq N\end{array}}}{\sum}\gamma_{n_1}\gamma_{n_2}\right)R-
      \left(3\sum_{n=N+1}^{2N+1}\tilde{\eps}^{2n/3}\gamma_{n-(N+1)}\right)R^2-\tilde{\eps}^{4(N+1)/3}R^3.\nonumber
\end{eqnarray}
\end{itemize}
The analysis below requires the precise knowledge of the behaviour of $\gamma_n(y)$ as
$y\to \pm\infty$. This behaviour was already described
in \cite{GP}, and it is summarized in the next two
propositions:

\begin{prop}\label{proposition-Painleve}
$\ $
The behaviour of
$\gamma_0$ as $y\to-\infty$ is described by
\begin{eqnarray}\label{asympnu0-}
\gamma_0(y)=\frac{1}{\sqrt{\pi}(-y)^{1/4}}\exp\left(-\frac{1}{3}(-y)^{3/2}\right)
\left(1+\mathcal{O}(|y|^{-3/4})\right)\underset{y\to-\infty}{\approx}0,
\end{eqnarray} whereas as $y\to +\infty$,
\begin{eqnarray}\label{asympnu0+}
\gamma_0(y)\underset{y\to+\infty}{\approx}y^{1/2}\sum_{n=0}^{\infty}
a_ny^{-3n}, \end{eqnarray} where $a_0=1$, and for $n
\geq 0$,
$$
a_{n+1}=2\left(9n^2-\frac{1}{4}\right)a_n-\frac{1}{2}\sum_{\substack{
  n_1+n_2+n_3=n+1\\n_1,n_2,n_3\leq n}}a_{n_1}a_{n_2}a_{n_3}.
$$
\end{prop}

\begin{rem}
The calculation of the first terms in (\ref{asympnu0+}) gives
\begin{eqnarray}\label{DAHM}
\gamma_0(y)=y^{1/2}-\frac{1}{2}y^{-5/2}-\frac{73}{8}y^{-11/2}+O(y^{-17/2}).
\end{eqnarray}
\end{rem}

\begin{prop}\label{proposition-mun}
$\ $ 
For every $n\geq 1$, 
$$\gamma_n(y) \underset{y\to +\infty}{\approx}
y^{1/2-2n}\sum_{m=0}^{\infty} g_{n,m} y^{-3m} \text{ for some } \{ g_{n,m}\}_{m \in \mathbb{N}},$$
and\hfill $\gamma_n(y) \underset{y\to -\infty}{\approx} 0.\quad$\hfill\null\\
Moreover, if $d=1$, for every $n\geq 1$, $g_{n,0}=0$.\\
For instance,
$\gamma_1(y)\underset{y\to +\infty}{\sim}\frac{5(7-d)}{4}y^{-9/2}$ if $d=1$,
whereas $\gamma_1(y)\underset{y\to +\infty}{\sim}=\frac{1-d}{2}y^{-3/2}$
if $d=2,3$.
\end{prop}

\subsection{Expansions of $\nu$ and $\lambda$ in
  $D_1$}\label{sec-calcnul}
For $x\in D_1$, we formally look for a solution $(\eta_1,\eta_2)$ to
(\ref{sys}) under the form given in (\ref{apD1}). Then, it turns out that $\nu$ and $\lambda$ have to solve
\begin{eqnarray}
-2d\eps^{2/3}\nu'+4R_1^2\nu''-4\eps^{2/3}y_1\nu''
+\frac{\alpha_0}{\alpha_2}\frac{R_2^2-R_1^2}{\eps^{2/3}}\nu+y_1\nu
-2\alpha_1\nu^3-2\alpha_0\lambda\nu&=&0\label{eq-nu}\\
-d\eps^{2/3}\lambda\lambda'-(R_1^2-\eps^{2/3}y_1){\lambda'}^2+2(R_1^2-\eps^{2/3}y_1)\lambda\lambda''+y_2\lambda^{2}-2\alpha_2\lambda^{3}-2\alpha_0\nu^2\lambda^{2}&=&0\label{eq-lambda}
\end{eqnarray}
Moreover, we are looking for solutions $(\eta_1,\eta_2)$ that converge
to the Thomas-Fermi limit $(\eta_{10},\eta_{20})$ as $\eps\to 0$. As a
result, according to Section \ref{asymp}, $\nu$ and $\lambda$ have to satisfy the following asymptotics. On the one side, if $R_1<|x|<R_2$ is fixed, $\eps\to 0$ if and only if $y_1\to-\infty$, and
$$\nu(y_1)\underset{y_1\to-\infty}{\longrightarrow}0,\quad\quad \lambda(y_1)\underset{y_1\to-\infty}{\sim}\left(\frac{R_2^2-R_1^2}{2\alpha_2\eps^{2/3}}+\frac{y_1}{2\alpha_2}\right).$$
On the other side, if $|x|<R_1$ is fixed, $\eps\to 0$ if and only if $y_1\to+\infty$, and
$$\nu(y_1)\underset{y_1\to+\infty}{\sim}\left(\frac{\Gamma_2y_1}{2\alpha_1\Gamma_{12}}\right)^{1/2},\quad\quad\lambda(y_1)\underset{y_1\to+\infty}{\sim}\frac{R_2^2-R_1^2}{2\alpha_2\eps^{2/3}}+\frac{\Gamma_1y_1}{2\alpha_2\Gamma_{12}}.$$
We formally develop $\nu$ and $\lambda$ into powers of $\eps^{2/3}$:
\begin{eqnarray}\label{DAlnu}
\nu(y_1)=\sum_{n=0}^\infty\eps^{2n/3}\nu_n(y_1),\quad
\lambda(y_1)=\sum_{n=-1}^\infty\eps^{2n/3}\lambda_n(y_1),
\end{eqnarray}
and we plug these expansions of $\nu$ and $\lambda$ into (\ref{eq-nu}). We obtain
\begin{eqnarray*}
-2d\sum_{n=1}^\infty\eps^{2n/3}\nu_{n-1}'+4R_1^2\sum_{n=0}^\infty\eps^{2n/3}\nu_{n}''-4y_1\sum_{n=1}^\infty\eps^{2n/3}\nu_{n-1}''+\frac{\alpha_0}{\alpha_2}(R_2^2-R_1^2)\sum_{n=-1}^\infty\eps^{2n/3}\nu_{n+1}+y_1\sum_{n=0}^\infty\eps^{2n/3}\nu_{n}&&\\
-2\alpha_1\sum_{n=0}^\infty\eps^{2n/3}\sum_{n_1+n_2+n_3=n}\nu_{n_1}\nu_{n_2}\nu_{n_3}-2\alpha_0\sum_{n=-1}^\infty\eps^{2n/3}\sum_{\tiny \begin{array}{c}n_1+n_2=n,\\ n_1\geq
-1,n_2\geq 0\end{array}}\lambda_{n_1}\nu_{n_2}=0.&&
\end{eqnarray*}
At order $n=-1$, we get, in agreement with the asymptotics of
$\lambda_{-1}$ given in Section \ref{asymp},
\begin{eqnarray}\label{I-1}
\lambda_{-1}(y_1)&=&\frac{R_2^2-R_1^2}{2\alpha_2},
\end{eqnarray}
and therefore the equation can be simplified into
\begin{eqnarray}\label{eq-DAnun}
-2d\sum_{n=1}^\infty\eps^{2n/3}\nu_{n-1}'+4R_1^2\sum_{n=0}^\infty\eps^{2n/3}\nu_{n}''-4y_1\sum_{n=1}^\infty\eps^{2n/3}\nu_{n-1}''+y_1\sum_{n=0}^\infty\eps^{2n/3}\nu_{n}&&\nonumber\\
-2\alpha_1\sum_{n=0}^\infty\eps^{2n/3}\sum_{n_1+n_2+n_3=n}\nu_{n_1}\nu_{n_2}\nu_{n_3}-2\alpha_0\sum_{n=0}^\infty\eps^{2n/3}\sum_{\tiny \begin{array}{c}n_1+n_2=n,\\ n_1\geq
0,n_2\geq 0\end{array}}\lambda_{n_1}\nu_{n_2}=0&&
\end{eqnarray}
At order $n=0$, we obtain
\begin{eqnarray}\label{I0}
4R_1^2\nu_0''+y_1\nu_0-2\alpha_1\nu_{0}^3-2\alpha_0\lambda_{0}\nu_{0}=0.
\end{eqnarray}
Similarly, plugging (\ref{DAlnu}) into (\ref{eq-lambda}) and multiplying by $\eps^{4/3}$, we get
\begin{eqnarray}\label{eq-DAln}
\lefteqn{-d\sum_{n=1}^\infty\eps^{2n/3}\sum_{\tiny\begin{array}{c}n_1+n_2=n-3\\n_1,n_2\geq
    -1\end{array}}\lambda_{n_1}'\lambda_{n_2}
+2R_1^2\sum_{n=0}^\infty\eps^{2n/3}\sum_{\tiny\begin{array}{c}n_1+n_2=n-2\\n_1,n_2\geq
    -1\end{array}}\lambda_{n_1}''\lambda_{n_2}}\nonumber\\
&&-2y_1\sum_{n=1}^\infty\eps^{2n/3}\sum_{\tiny\begin{array}{c}n_1+n_2=n-3\\n_1,n_2\geq
    -1\end{array}}\lambda_{n_1}''\lambda_{n_2}
-R_1^2\sum_{n=0}^\infty\eps^{2n/3}\sum_{\tiny\begin{array}{c}n_1+n_2=n-2\\n_1,n_2\geq
    -1\end{array}}\lambda_{n_1}'\lambda_{n_2}'\nonumber\\
&&+y_1\sum_{n=1}^\infty\eps^{2n/3}\sum_{\tiny\begin{array}{c}n_1+n_2=n-3\\n_1,n_2\geq
    -1\end{array}}\lambda_{n_1}'\lambda_{n_2}'
+(R_2^2-R_1^2)\sum_{n=-1}^\infty\eps^{2n/3}\sum_{\tiny\begin{array}{c}n_1+n_2=n-1\\n_1,n_2\geq
    -1\end{array}}\lambda_{n_1}\lambda_{n_2}\nonumber\\
&&+y_1\sum_{n=0}^\infty\eps^{2n/3}\sum_{\tiny\begin{array}{c}n_1+n_2=n-2\\n_1,n_2\geq
    -1\end{array}}\lambda_{n_1}\lambda_{n_2}
-2\alpha_2\sum_{n=-1}^\infty\eps^{2n/3}\sum_{\tiny\begin{array}{c}n_1+n_2+n_3=n-2\\n_1,n_2,n_3\geq
    -1\end{array}}\lambda_{n_1}\lambda_{n_2}\lambda_{n_3}\nonumber\\
&&-2\alpha_0\sum_{n=0}^\infty\eps^{2n/3}\sum_{\tiny\begin{array}{c}n_1+n_2+n_3+n_4=n-2\\n_1,n_2\geq
    -1,\ n_3,n_4\geq 0\end{array}}\lambda_{n_1}\lambda_{n_2}\nu_{n_3}\nu_{n_4}=0
\end{eqnarray}
This equation at order $n=-1$ is satisfied thanks to (\ref{I-1}). At
order $n=0$, we obtain 
\begin{eqnarray}\label{II0}
\lambda_{0}(y_1)=\frac{y_1}{2\alpha_2}-\frac{\alpha_0}{\alpha_2}\nu_{0}(y_1)^2.
\end{eqnarray}
From (\ref{I0}) and (\ref{II0}), we infer the equation satisfied by
$\nu_0$:
\begin{eqnarray}\label{eq-nu0}
4R_1^2\nu_0''+\Gamma_2y_1\nu_0-2\alpha_1\Gamma_{12}\nu_{0}^3=0.
\end{eqnarray}
Moreover, according to Section \ref{asymp}, the asymptotic behaviour we need for $\nu_0$ is
$$\nu_0(y_1)\underset{y_1\to+\infty}{\sim}\left(\frac{\Gamma_2
  y_1}{2\alpha_1\Gamma_{12}}\right)^{1/2},\quad
\nu_0(y_1)\underset{y_1\to-\infty}{\to}0.$$
Looking for $\nu_0$ under the form 
$$\nu_0(y_1)=\frac{R_1^{1/3}|\Gamma_2|^{1/3}}{(2\alpha_1)^{1/2}|\Gamma_{12}|^{1/2}}\gamma\left(\frac{|\Gamma_2|^{1/3}y_1}{R_1^{2/3}}\right),$$
$\nu_0$ solves (\ref{eq-nu0}) if and only if $\gamma$ solves
\begin{eqnarray}\label{P2-eq-sign}
4{\rm sign}(\Gamma_2)\gamma''(y)+y\gamma(y)-\gamma(y)^3=0,\quad
y\in\R,
\end{eqnarray}
with the boundary conditions
\begin{eqnarray}\label{BCgamma}
\gamma(y)\underset{y\to+\infty}{\sim}\sqrt{y},\quad
\gamma(y)\underset{y\to-\infty}{\to}0.
\end{eqnarray}
If the sign of $\Gamma_2$ (which is the same as the sign of
$\Gamma_{12}$ according to (\ref{posgamma})) is negative, it can be
easily seen that (\ref{P2-eq-sign}) has no non-trivial solution with
fast decay to 0 as $y\to -\infty$. Indeed, if $\gamma$ solves (\ref{P2-eq-sign})
with $\gamma'(y)\to 0$ and $y\gamma(y)^2\to 0$ as $y\to -\infty$, then
by integration between $-\infty$ and $y$, we get 
$$2{\rm sign}(\Gamma_2)\gamma'(y)^2=-y\frac{\gamma(y)^2}{2}+\int_{-\infty}^y\frac{\gamma(t)^2}{2}dt+\frac{\gamma(y)^4}{4},$$
which implies $\gamma\equiv 0$ if $\Gamma_2<0$ and $y<0$. Also, from
now on, we assume 
\begin{eqnarray}\label{gamma2pos}
\Gamma_2>0,&& \Gamma_{12}>0.
\end{eqnarray}
Under this condition, $\gamma$ has to be the Hastings-McLeod solution
$\gamma_0$ of
the Painlev\'e II equation (\ref{P2-eq}), and
\begin{eqnarray}\label{nu00}
\nu_0(y_1)=\frac{R_1^{1/3}\Gamma_2^{1/3}}{(2\alpha_1)^{1/2}\Gamma_{12}^{1/2}}\gamma_0\left(\frac{\Gamma_2^{1/3}y_1}{R_1^{2/3}}\right).
\end{eqnarray}
Thanks to (\ref{I-1}), equation (\ref{eq-DAnun}) at order $n\geq 1$ gives  
\begin{eqnarray}\label{preIn}
\lefteqn{4R_1^2\nu_{n}''+y_1\nu_{n}-6\alpha_1\nu_0^2\nu_n-2\alpha_0\lambda_0\nu_n-2\alpha_0\lambda_n\nu_0}\\
&=&2d\nu_{n-1}'+4y_1\nu_{n-1}''+2\alpha_1\sum_{\tiny\begin{array}{c}n_1+n_2+n_3=n\\0\leq
  n_1,n_2,n_3\leq
  n-1\end{array}}\nu_{n_1}\nu_{n_2}\nu_{n_3}+2\alpha_0\sum_{\tiny\begin{array}{c}n_1+n_2=n,\\1\leq
  n_1,n_2\leq
n-1\end{array}}\lambda_{n_1}\nu_{n_2}.\nonumber
\end{eqnarray}
On the other side, equation (\ref{eq-DAln}) at order $n\geq 1$ yields
\begin{eqnarray}\label{star}
\lefteqn{-d\sum_{\tiny\begin{array}{c}n_1+n_2=n-3\\n_1,n_2\geq
    -1\end{array}}\lambda_{n_1}'\lambda_{n_2}
+2R_1^2\sum_{\tiny\begin{array}{c}n_1+n_2=n-2\\n_1,n_2\geq
    -1\end{array}}\lambda_{n_1}''\lambda_{n_2}}\nonumber\\
&&-2y_1\sum_{\tiny\begin{array}{c}n_1+n_2=n-3\\n_1,n_2\geq
    -1\end{array}}\lambda_{n_1}''\lambda_{n_2}
-R_1^2\sum_{\tiny\begin{array}{c}n_1+n_2=n-2\\n_1,n_2\geq
    -1\end{array}}\lambda_{n_1}'\lambda_{n_2}'\nonumber\\
&&+y_1\sum_{\tiny\begin{array}{c}n_1+n_2=n-3\\n_1,n_2\geq
    -1\end{array}}\lambda_{n_1}'\lambda_{n_2}'
+(R_2^2-R_1^2)\sum_{\tiny\begin{array}{c}n_1+n_2=n-1\\n_1,n_2\geq
    0\end{array}}\lambda_{n_1}\lambda_{n_2}+2(R_2^2-R_1^2)\lambda_{-1}\lambda_n\nonumber\\
&&+y_1\sum_{\tiny\begin{array}{c}n_1+n_2=n-2\\n_1,n_2\geq
    -1\end{array}}\lambda_{n_1}\lambda_{n_2}
-2\alpha_2\sum_{\tiny\begin{array}{c}n_1+n_2+n_3=n-2\\n-1\geq n_1,n_2,n_3\geq
    -1\end{array}}\lambda_{n_1}\lambda_{n_2}\lambda_{n_3}-6\alpha_2\lambda_{-1}^2\lambda_n\nonumber\\
&&-2\alpha_0\sum_{\tiny\begin{array}{c}n_1+n_2+n_3+n_4=n-2\\n_1,n_2\geq
    -1,\ n-1\geq n_3,n_4\geq 0\end{array}}\lambda_{n_1}\lambda_{n_2}\nu_{n_3}\nu_{n_4}-4\alpha_0\lambda_{-1}^2\nu_0\nu_n=0,
\end{eqnarray}
therefore for $n\geq 1$,
\begin{eqnarray}\label{IIn}
\lambda_n=-2\frac{\alpha_0}{\alpha_2}\nu_0\nu_n+\frac{2\alpha_2}{(R_2^2-R_1^2)^2}\delta_n,
\end{eqnarray}
where
\begin{eqnarray}\label{def-deltan}
\delta_n &=&\sum_{\tiny\begin{array}{c}n_1+n_2=n-3\\n_1,n_2\geq
    -1\end{array}}\left(-d\lambda_{n_1}'\lambda_{n_2}-2y_1\lambda_{n_1}''\lambda_{n_2}+y_1\lambda_{n_1}'\lambda_{n_2}'\right)\nonumber\\
&&+\sum_{\tiny\begin{array}{c}n_1+n_2=n-2\\n_1,n_2\geq
    -1\end{array}}\left(2R_1^2\lambda_{n_1}''\lambda_{n_2}-R_1^2\lambda_{n_1}'\lambda_{n_2}'+y_1\lambda_{n_1}\lambda_{n_2}\right)+(R_2^2-R_1^2)\sum_{\tiny\begin{array}{c}n_1+n_2=n-1\\n_1,n_2\geq
    0\end{array}}\lambda_{n_1}\lambda_{n_2}\nonumber\\
&&-2\alpha_2\sum_{\tiny\begin{array}{c}n_1+n_2+n_3=n-2\\n-1\geq n_1,n_2,n_3\geq
    -1\end{array}}\lambda_{n_1}\lambda_{n_2}\lambda_{n_3}
-2\alpha_0\sum_{\tiny\begin{array}{c}n_1+n_2+n_3+n_4=n-2\\n_1,n_2\geq
    -1,\ n-1\geq n_3,n_4\geq 0\end{array}}\lambda_{n_1}\lambda_{n_2}\nu_{n_3}\nu_{n_4}.
\end{eqnarray}
At this stage, we have constructed $\lambda_{-1}$, $\nu_0$ and
$\lambda_0$, which are given respectively by (\ref{I-1}), (\ref{nu00})
and (\ref{II0}). For $n\geq 1$, the $\lambda_n$'s and the $\nu_n$'s are
constructed by induction as follows. Let $n\geq 1$, and assume that
the $\lambda_k$'s and the $\nu_k$'s are
known for every $k\leq n-1$. Then, plugging (\ref{IIn}) and (\ref{II0}) into
(\ref{preIn}), $\nu_n$ has to solve
\begin{eqnarray}\label{In}
T\nu_n=F_n,
\end{eqnarray}
where
\begin{eqnarray}\label{def-T}
T=-4R_1^2\partial^2_{y_1}+W(y_1), \quad W(y_1)=6\alpha_1\Gamma_{12}\nu_0^2-\Gamma_2y_1
\end{eqnarray}
and
\begin{eqnarray}\label{def-Fn}
F_n=-\frac{4\alpha_0\alpha_2\nu_0}{(R_2^2-R_1^2)^2}\delta_n
-2d\nu_{n-1}'-4y_1\nu_{n-1}''-2\alpha_1\!\!\!\!\!\!\!\!\!\!\!\sum_{\tiny\begin{array}{c}n_1+n_2+n_3=n\\0\leq
  n_1,n_2,n_3\leq
  n-1\end{array}}\!\!\!\!\!\!\nu_{n_1}\nu_{n_2}\nu_{n_3}-2\alpha_0\!\!\!\!\!\!\!\!\sum_{\tiny\begin{array}{c}n_1+n_2=n,\\1\leq
  n_1,n_2\leq
n-1\end{array}}\!\!\!\!\!\!\lambda_{n_1}\nu_{n_2}.
\end{eqnarray}
Note that only $\lambda_k$'s and $\nu_k$'s for
$k\leq n-1$ appear in (\ref{def-Fn}) and (\ref{def-deltan}). Once
(\ref{In}) has been solved, $\lambda_n$ is given by (\ref{IIn}). In
order to invert $T$ in (\ref{In}), one needs to
understand the behaviour of $F_n(y_1)$ as $y_1\to\pm\infty$. Thus,
$\delta_n$, $F_n$, $\nu_n$ and $\lambda_n$ will be constructed recursively in such a way
that for every $n\geq 1$,
\begin{eqnarray}\label{asymp-deln}
\delta_n(y_1)\underset{y_1\to -\infty}{\approx}y_1^{n-2}\sum_{0\leq
  m\leq (n-2)/3}\widetilde{D}_{n,m}y_1^{-3m},\quad\delta_n(y_1)\underset{y_1\to
  +\infty}{\approx}y_1^{n-2}\sum_{m=0}^{+\infty}D_{n,m}y_1^{-3m},
\end{eqnarray}
\begin{eqnarray}\label{asymp-Fn}
F_n(y_1)\underset{y_1\to -\infty}{\approx}0,\quad F_n(y_1)\underset{y_1\to
  +\infty}{\approx}y_1^{n-3/2}\sum_{m=0}^{+\infty}F_{n,m}y_1^{-3m},
\end{eqnarray}
\begin{eqnarray}\label{asymp-nun}
\nu_n(y_1)\underset{y_1\to -\infty}{\approx}0,\quad\nu_n(y_1)\underset{y_1\to
  +\infty}{\approx}y_1^{n-5/2}\sum_{m=0}^{+\infty}N_{n,m}y_1^{-3m},
\end{eqnarray}
\begin{eqnarray}\label{asymp-lamn}
\lambda_n(y_1)\underset{y_1\to -\infty}{\approx}y_1^{n-2}\sum_{0\leq
  m\leq (n-2)/3}\widetilde{L}_{n,m}y_1^{-3m},\quad
\lambda_n(y_1)\underset{y_1\to
  +\infty}{\approx}y_1^{n-2}\sum_{m=0}^{+\infty}L_{n,m}y_1^{-3m},
\end{eqnarray}
where the $D_{n,m}$'s, $F_{n,m}$'s, $N_{n,m}$'s, $L_{n,m}$'s, $\widetilde{D}_{n,m}$'s and $\widetilde{L}_{n,m}$'s are some real coefficients.
Note that thanks to
(\ref{nu00}), (\ref{II0}), (\ref{asympnu0-}) and (\ref{asympnu0+}),
$\nu_0$ and $\lambda_0$ admit similar expansions. However, the power of
the leading term in the expansions they satisfy as $y_1\to +\infty$ (and
for $\lambda_0$, also as $y_1\to -\infty$) is
higher of three units to the one which would be given by
(\ref{asymp-nun}) and (\ref{asymp-lamn}) for $n=0$. More precisely, we have
\begin{eqnarray}\label{asymp-nu0}
\nu_0(y_1)\underset{y_1\to -\infty}{\approx}0,\quad\nu_0(y_1)\underset{y_1\to
  +\infty}{\approx}y_1^{1/2}\sum_{m=0}^{+\infty}N_{0,m}y_1^{-3m}
\end{eqnarray}
and
\begin{eqnarray}\label{asymp-lam0}
\lambda_0(y_1)\underset{y_1\to -\infty}{\approx}\frac{y_1}{2\alpha_2},\quad\lambda_0(y_1)\underset{y_1\to +\infty}{\approx}y_1\sum_{m=0}^{+\infty}L_{0,m}y_1^{-3m},
\end{eqnarray}
where 
\be\label{626'}
N_{0m}=\left(\frac{\Gamma_2}{2\alpha_1\Gamma_{12}}\right)^{1/2}\left(\frac{R_1^2}{\Gamma_2}\right)^ma_m
\ee
and
\be\label{626''}
L_{00}=\frac{\Gamma_1}{2\alpha_2\Gamma_{12}}\quad \text{and
  for}\ m\geq 1,\quad
L_{0m}=-\frac{\alpha_0}{\alpha_2}\sum_{m_1+m_2=m}N_{0m_1}N_{0m_2}.
\ee
Next, let us explain why $\delta_1$, $F_1$, $\nu_1$ and $\lambda_1$
admit asymptotic expansions like the ones given in (\ref{asymp-deln}),
(\ref{asymp-Fn}), (\ref{asymp-nun}),
(\ref{asymp-lamn}) and let us calculate explicitely the first terms in
these expansions. Thanks to (\ref{def-deltan}) for $n=1$ as well as
(\ref{I-1}), (\ref{II0}), (\ref{asymp-lam0}), (\ref{626'}), (\ref{626''}),we have
\begin{eqnarray}\label{def-del1}
\delta_1 &=&2y_1\lambda_{-1}\lambda_0+2R_1^2\lambda_0''\lambda_{-1}+(R_2^2-R_1^2)\lambda_0^2-6\alpha_2\lambda_0^2\lambda_{-1}-4\alpha_0\lambda_{-1}\lambda_0\nu_0^2\nonumber\\
&=&2R_1^2\lambda_0''\lambda_{-1}\underset{y_1\to +\infty}{=}\frac{3\alpha_0R_1^4(R_2^2-R_1^2)}{\alpha_1\alpha_2^2\Gamma_{12}}y_1^{-4}+O(y_1^{-7}).
\end{eqnarray}
Thus, the asymptotics as $y_1\to +\infty$ in (\ref{asymp-deln}) holds with
\be\label{d10d11}
D_{10}=0\quad {\rm  and}\quad D_{11}=\frac{3\alpha_0R_1^4(R_2^2-R_1^2)}{\alpha_1\alpha_2^2\Gamma_{12}}.
\ee
From (\ref{asymp-lam0}) and (\ref{def-del1}), we also infer that
$\delta_1\underset{y_1\to -\infty}{\approx}0$, which is the asymptotics as $y_1\to -\infty$ in (\ref{asymp-deln}).

Then, (\ref{def-Fn}) yields
\begin{eqnarray}\label{F1}
\lefteqn{F_1=-\frac{4\alpha_0\alpha_2\nu_0}{(R_2^2-R_1^2)^2}\delta_1
-2d\nu_{0}'-4y_1\nu_{0}''}\\
&\underset{y_1\to +\infty}{=}&\left(\frac{\Gamma_2}{2\alpha_1\Gamma_{12}}\right)^{1/2}(1-d)y_1^{-1/2}+\frac{R_1^2}{(2\alpha_1\Gamma_2\Gamma_{12})^{1/2}}\left(\frac{5}{2}(7-d)-\frac{12\alpha_0^2\Gamma_2R_1^2}{\alpha_1\alpha_2\Gamma_{12}(R_2^2-R_1^2)}\right)y_1^{-7/2}+O(y_1^{-13/2}).\nonumber
\end{eqnarray}
Thus, (\ref{asymp-nu0}) and (\ref{asymp-deln}) for $n=1$ imply that
(\ref{asymp-Fn}) for $n=1$ holds with
$$F_{10}=\left(\frac{\Gamma_2}{2\alpha_1\Gamma_{12}}\right)^{1/2}(1-d)\quad
{\rm  and}\quad F_{11}=\frac{R_1^2}{(2\alpha_1\Gamma_2\Gamma_{12})^{1/2}}\left(\frac{5}{2}(7-d)-\frac{12\alpha_0^2\Gamma_2R_1^2}{\alpha_1\alpha_2\Gamma_{12}(R_2^2-R_1^2)}\right).$$
In order to calculate $\nu_1$ from (\ref{In}), let us first notice
that the function $W$ defined in (\ref{def-T}) coincides, up to a
rescaling, to the function $W_0(y)=3\gamma_0(y)^2-y$ which was studied
in \cite{GP}. On the other side, $W$ can be expressed in terms of
$\lambda_0$ thanks to (\ref{II0}). Namely,
\begin{eqnarray}\label{asymp-W0}
W(y_1)=\Gamma_2^{2/3}R_1^{2/3}W_0\left(\frac{\Gamma_2^{1/3}y_1}{R_1^{2/3}}\right)=\left(\frac{3\alpha_1\Gamma_1}{\alpha_0}+2\Gamma_2\right)y_1-\frac{6\alpha_1\alpha_2\Gamma_{12}}{\alpha_0}\lambda_0(y_1).
\end{eqnarray}
In particular, there exists $C>0$ such that $W(y_1)\geq C$ for every
$y_1\in \R$, and $W$ admits the asymptotic expansions
\begin{eqnarray}\label{asymp-W}
W(y_1)\underset{y_1\to -\infty}{\approx}-\Gamma_2 y_1,\quad
W(y_1)\underset{y_1\to +\infty}{\approx}y_1\left(2\Gamma_2-\frac{6\alpha_1\alpha_2\Gamma_{12}}{\alpha_0}\sum_{m=1}^{+\infty}L_{0,m}y_1^{-3m}\right)
\end{eqnarray}
In the case $d=1$, since $F_{10}=0$, $F_1\in L^2(\R)$, and $\nu_1$ is
obtained by inversion of $T$, which is a Schr\"odinger operator on
$L^2(\R)$. Moreover, thanks to (\ref{asymp-W}) and
the positiveness of $W$, Lemma 2.1 in \cite{GP} implies that
the solution $\nu_1$ to (\ref{In}) admits asymptotic expansions like
the ones given in (\ref{asymp-nun}), with
$$N_{1,0}=0,\quad
N_{1,1}=-\frac{6R_1^2}{(2\alpha_1\Gamma_{12}\Gamma_2)^{1/2}}\left(\frac{\alpha_0^2R_1^2}{\alpha_1\alpha_2(R_2^2-R_1^2)\Gamma_{12}}-\frac{5}{4\Gamma_2}\right).$$
In the cases $d=2,3$, $F_{10}\neq 0$, and therefore $F_1\not\in
L^2(\R)$. We construct the solution $\nu_1$ to (\ref{In}) by using the
same trick as in \cite{GP}. Namely, we look for $\nu_1$ under the
form
$$\nu_1=\frac{F_{1,0}y_1^{-1/2}}{W(y_1)}\Phi(y_1)+\widetilde{\nu}_1,$$
where $\Phi\in\mathcal{C}^\infty(\R)$ is such that $\Phi(y_1)\equiv
0$ for $y_1\leq 1/2$ and $\Phi(y_1)\equiv
1$ for $y_1\geq 1$, in such a way that
\begin{eqnarray}\label{eq-nu1t}
(-4R_1^2\partial_{y_1}^2+W(y_1))\widetilde{\nu}_1&=&F_1-F_{1,0}y_1^{-1/2}\Phi(y_1)+4R_1^2\frac{d^2}{dy_1^2}\left[\frac{F_{1,0}y_1^{-1/2}}{W(y_1)}\Phi(y_1)\right].
\end{eqnarray}
The right hand side of (\ref{eq-nu1t}) behaves now like $O(y_1^{-7/2})$
as $y_1\to +\infty$, and its behaviour at $-\infty$ is the same as the
one of $F_1$, therefore the right hand side in (\ref{eq-nu1t}) belongs to
$L^2(\R)$, and (\ref{eq-nu1t}) has a unique solution $\widetilde{\nu}_1$
in $L^2(\R)$. Moreover, again thanks to Lemma 2.1 in \cite{GP}, we
deduce the existence of asymptotic expansions for $\widetilde{\nu}_1$ as
$y_1\to\pm\infty$, with
$\widetilde{\nu}_1(y_1)\underset{y_1\to+\infty}{=}O(y_1^{-9/2})$. These expansions
 for $\widetilde{\nu}_1$ imply that $\nu_1$ has expansions like in
(\ref{asymp-nun}), with
$$N_{1,0}=\frac{F_{1,0}}{2\Gamma_2}=\frac{1-d}{2(2\alpha_1\Gamma_2\Gamma_{12})^{1/2}}.$$
Then, from (\ref{IIn}) and (\ref{def-del1}),
\be\label{lambda1}
\lambda_1 &=&-2\frac{\alpha_0}{\alpha_2}\nu_0\nu_1+\frac{4\alpha_2
  R_1^2}{(R_2^2-R_1^2)^2}\lambda_0''\lambda_{-1}
\ee
and thanks to ({\ref{asymp-nu0}), ({\ref{def-del1}) and ({\ref{asymp-nun}) for $n=1$, $\lambda_1$ as asymptotic expansions
like in (\ref{asymp-lamn}) (in particular,
$\lambda_1(y_1)\underset{y_1\to -\infty}{\approx}0$), with
$$L_{1,0}=\frac{\alpha_0(d-1)}{2\alpha_1\alpha_2\Gamma_{12}},\quad{\rm
  and\ if\ }d=1,\quad L_{1,1}=\frac{3\alpha_0R_1^2}{2\alpha_1\alpha_2\Gamma_{12}}\left(\frac{4R_1^2}{(R_2^2-R_1^2)\Gamma_{12}}-\frac{5}{\Gamma_2}\right).$$
For $n=2$, (\ref{def-deltan}) and (\ref{def-Fn}) give similarly (after
simplifications involving also (\ref{I-1}), (\ref{II0}) and (\ref{lambda1}))
\be
\delta_2 & = &
-d\lambda_0'\lambda_{-1}-2y_1\lambda_0''\lambda_{-1}-2R_1^2\lambda_0''\lambda_0-R_1^2\lambda_0'^2+2R_1^2\lambda_1''\lambda_{-1}-2\alpha_0\lambda_{-1}^2\nu_1^2
\ee
and
\be
F_2&=&
-\frac{4\alpha_0\alpha_2\nu_0}{(R_2^2-R_1^2)^2}\delta_2-2d\nu_1'-4y_1\nu_1''-6\alpha_1\nu_0\nu_1^2-2\alpha_0\lambda_1\nu_1
\ee
which implies, thanks to the expansions calculated previously for
$\lambda_0$, $\nu_0$, $\lambda_1$ and $\nu_1$ that $\delta_2$ and
$F_2$ satisfy respectively (\ref{asymp-deln}) and (\ref{asymp-Fn})
for $n=2$, with
\be\label{D20}
D_{2,0}=-\frac{\Gamma_1\left(d\Gamma_{12}(R_2^2-R_1^2)+\Gamma_1R_1^2\right)}{4\alpha_2^2\Gamma_{12}^2},\quad
\widetilde{D}_{2,0}=-\frac{d(R_2^2-R_1^2)+R_1^2}{4\alpha_2^2}
\ee
and 
$$F_{2,0}=\frac{\alpha_0}{\alpha_2}\left(\frac{\Gamma_2}{2\alpha_1\Gamma_{12}}\right)^{1/2}\frac{\Gamma_1}{\Gamma_{12}^2}\frac{(d\Gamma_{12}(R_2^2-R_1^2)+\Gamma_1R_1^2)}{(R_2^2-R_1^2)^2}.$$
In order to solve (\ref{In}) for $n=2$, we look for $\nu_2$ under the form
\be\label{dec-nu2}
\nu_2 &=& \frac{F_{2,0}}{W(y_1)}y_1^{1/2}\Phi(y_1)+\widetilde{\nu}_2
\ee
Then $\nu_2$ solves (\ref{In}) for $n=2$ if and only if
$\widetilde{\nu}_2$ solves
\be
(-4R_1^2\partial_{y_1}^2+W(y_1))\widetilde{\nu}_2&=&F_2-F_{2,0}y_1^{1/2}\Phi(y_1)+4R_1^2\frac{d^2}{dy_1^2}\left[\frac{F_{2,0}y_1^{1/2}}{W(y_1)}\Phi(y_1)\right]\underset{y_1\to+\infty}{=}O(y_1^{-5/2}).\nonumber
\ee
In particular, the right hand side in this equation belongs to
$L^2(\R)$. Thus by inversion of $T$ like for $n=1$ and $d=2,3$, and coming
back to (\ref{dec-nu2}), $\nu_2$ satisfies (\ref{asymp-nun}) for
$n=2$, with 
$$N_{2,0}=\frac{F_{2,0}}{2\Gamma_2}=\frac{\alpha_0}{\alpha_2}\frac{1}{\left(2\alpha_1\Gamma_2\Gamma_{12}\right)^{1/2}}\frac{\Gamma_1}{\Gamma_{12}^2}\frac{(d\Gamma_{12}(R_2^2-R_1^2)+\Gamma_1R_1^2)}{2(R_2^2-R_1^2)^2}.$$
As a result, from (\ref{IIn}) for $n=2$,
 $\lambda_2$ satisfies (\ref{asymp-lamn}) for $n=2$, with
\be
L_{2,0}=
-\frac{\Gamma_1\left(d\Gamma_{12}(R_2^2-R_1^2)+\Gamma_1R_1^2\right)}{2(R_2^2-R_1^2)^2\Gamma_{12}^3\alpha_2}\quad{\rm
  and}\quad \widetilde{L}_{2,0}=-\frac{d(R_2^2-R_1^2)+R_1^2}{2\alpha_2(R_2^2-R_1^2)}.
\ee
Next, let us fix $n\geq 3$. We assume that we have constructed the
$\nu_k$'s for $k\in\{1\cdots n-1\}$, and that asymptotic expansions
(\ref{asymp-deln}), (\ref{asymp-Fn}), (\ref{asymp-nun}),
(\ref{asymp-lamn}) with $n$ replaced by each of these $k$'s are satisfied. 
Then it is clear from (\ref{def-deltan}) and (\ref{def-Fn}) that $F_n\approx 0$ as $y_1\to -\infty$ as
indicated in (\ref{asymp-Fn}). In order to
study the asymptotic expansion of $\delta_n$ as $y_1\to +\infty$, let
us first focus on the first sum in the right hand side of
(\ref{def-deltan}). If $n_1,n_2\geq 1$ and $n_1+n_2=n-3$, then it follows from
(\ref{asymp-lamn}) that
$h_{n_1n_2}(y_1):=-d\lambda_{n_1}'\lambda_{n_2}-2y_1\lambda_{n_1}''\lambda_{n_2}+y_1\lambda_{n_1}'\lambda_{n_2}'$
admits an asymptotic expansions
which can be written as
\begin{eqnarray}\label{DA-type}
y_1^{n-l}\sum_{m=0}^{+\infty}c_my_1^{-3m},
\end{eqnarray}
for some coefficients $(c_m)_{m\in\N}$, with
$l=8$. Thus, we deduce that
$$\sum_{\tiny\begin{array}{c}n_1+n_2=n-3\\n_1,n_2\geq
    1\end{array}}h_{n_1n_2}(y_1)\underset{y_1\to
  +\infty}{\approx}y_1^{n-8}\sum_{m=0}^{+\infty}\check{D}_{n,m+2}y_1^{-3m}=y_1^{n-2}\sum_{m=2}^{+\infty}\check{D}_{n,m}y_1^{-3m},$$
for some coefficients $(\check{D}_{n,m})_{m\geq 2}$. Similarly,
$h_{n_1n_2}(y_1)$ has an asymptotic expansion which can be written as
(\ref{DA-type}) with $l=5$ if $n_1\geq 1$ and $n_2\leq 0$ or $n_1\leq
0$ and $n_2\geq 1$, and with $l=2$ if $n_1,n_2\leq 0$. As a result,
the first sum in the right hand side of (\ref{def-deltan}) admits an asymptotic
expansion as $y_1\to +\infty$ like (\ref{DA-type}) with $l=2$, and in order to calculate
the leading term $cy_1^{n-2}$ in this expansion, one has to consider
only the terms of the sum corresponding to indices
$n_1,n_2\in\{-1,0\}$. The same kind of arguments applied to the other
terms in the right hand side of (\ref{def-deltan}) yields the
asymptotic expansion of $\delta_n$ as $y_1\to +\infty$ given in
(\ref{asymp-deln}). Moreover, in order to express $D_{n0}$, the only
terms of the right hand side of (\ref{def-deltan}) which have to be considered
are written in the calculation below, where we use (\ref{asymp-lam0}), (\ref{I-1}),
(\ref{asymp-lamn}), (\ref{asymp-nu0}) and (\ref{asymp-nun})
\begin{eqnarray}\label{DLNp}
\delta_n&\underset{y_1\to
  +\infty}{=}&\left(-d\lambda_{0}'\lambda_{n-3}-2y_1\lambda_{0}''\lambda_{n-3}+y_1\lambda_{0}'\lambda_{n-3}'\right)\mathbf{1}_{\{n=3\}}\nonumber\\
&&+2y_1\lambda_{-1}\lambda_{n-1}+2y_1\lambda_{0}\lambda_{n-2} +2(R_2^2-R_1^2)\lambda_0\lambda_{n-1}
-12\alpha_2\lambda_{-1}\lambda_0\lambda_{n-1}-6\alpha_2\lambda_{0}^2\lambda_{n-2}\nonumber\\
&&
-8\alpha_0\lambda_{-1}\lambda_0\nu_0\nu_{n-1}-4\alpha_0\lambda_{-1}\lambda_{n-1}\nu_0^2-4\alpha_0\lambda_0^2\nu_0\nu_{n-2}-4\alpha_0\lambda_{0}\lambda_{n-2}\nu_0^2+O(y_1^{n-5})\nonumber\\
&\underset{y_1\to
  +\infty}{=}&
\left(-d\frac{\Gamma_1}{2\alpha_2\Gamma_{12}}\frac{\Gamma_1}{2\alpha_2\Gamma_{12}}y_1+y_1\frac{\Gamma_1}{2\alpha_2\Gamma_{12}}\frac{\Gamma_1}{2\alpha_2\Gamma_{12}}\right)\mathbf{1}_{\{n=3\}}\nonumber\\
&&+2y_1\frac{R_2^2-R_1^2}{2\alpha_2}y_1^{n-3}L_{n-1,0}+2y_1\left(\frac{\Gamma_1}{2\alpha_2\Gamma_{12}}y_1\right)y_1^{n-4}L_{n-2,0}\nonumber\\
&&+2(R_2^2-R_1^2)\left(\frac{\Gamma_1}{2\alpha_2\Gamma_{12}}y_1\right)y_1^{n-3}L_{n-1,0}\nonumber\\
&&
-12\alpha_2\frac{R_2^2-R_1^2}{2\alpha_2}\left(\frac{\Gamma_1}{2\alpha_2\Gamma_{12}}y_1\right)y_1^{n-3}L_{n-1,0}-6\alpha_2\left(\frac{\Gamma_1}{2\alpha_2\Gamma_{12}}y_1\right)^2y_1^{n-4}L_{n-2,0}\nonumber\\
&&
-8\alpha_0\frac{R_2^2-R_1^2}{2\alpha_2}\left(\frac{\Gamma_1}{2\alpha_2\Gamma_{12}}y_1\right)\left(\left(\frac{\Gamma_2}{2\alpha_1\Gamma_{12}}\right)^{1/2}y_1^{1/2}\right)y_1^{n-7/2}N_{n-1,0}\nonumber\\
&&-4\alpha_0\frac{R_2^2-R_1^2}{2\alpha_2}y_1^{n-3}L_{n-1,0}\left(\left(\frac{\Gamma_2}{2\alpha_1\Gamma_{12}}\right)^{1/2}y_1^{1/2}\right)^2\nonumber\\
&&
-4\alpha_0\left(\frac{\Gamma_1}{2\alpha_2\Gamma_{12}}y_1\right)^2\left(\left(\frac{\Gamma_2}{2\alpha_1\Gamma_{12}}\right)^{1/2}y_1^{1/2}\right)y_1^{n-9/2}N_{n-2,0}\nonumber\\
&&-4\alpha_0\left(\frac{\Gamma_1}{2\alpha_2\Gamma_{12}}y_1\right)y_1^{n-4}L_{n-2,0}\left(\left(\frac{\Gamma_2}{2\alpha_1\Gamma_{12}}\right)^{1/2}y_1^{1/2}\right)^2+O(y_1^{n-5})\nonumber\\
&\underset{y_1\to
  +\infty}{=}& 
y_1\frac{(1-d)\Gamma_1^2}{4\alpha_2^2\Gamma_{12}^2}\mathbf{1}_{\{n=3\}}
+y_1^{n-2}\left(\frac{R_2^2-R_1^2}{\alpha_2}L_{n-1,0}+\frac{\Gamma_1}{\alpha_2\Gamma_{12}}L_{n-2,0}\right)+(R_2^2-R_1^2)y_1^{n-2}\frac{\Gamma_1}{\alpha_2\Gamma_{12}}L_{n-1,0}\nonumber\\
&&
-\frac{3(R_2^2-R_1^2)\Gamma_1}{\alpha_2\Gamma_{12}}y_1^{n-2}L_{n-1,0}-\frac{3\Gamma_1^2}{2\alpha_2\Gamma_{12}^2}y_1^{n-2}L_{n-2,0}
-2\alpha_0\frac{(R_2^2-R_1^2)\Gamma_1\Gamma_2^{1/2}N_{n-1,0}}{\alpha_2^2(2\alpha_1)^{1/2}\Gamma_{12}^{3/2}}y_1^{n-2}\nonumber\\
&&-\alpha_0\frac{(R_2^2-R_1^2)\Gamma_2L_{n-1,0}}{\alpha_1\alpha_2\Gamma_{12}}y_1^{n-2}
-\frac{\alpha_0\Gamma_1^2\Gamma_2^{1/2}N_{n-2,0}}{\alpha_2^2(2\alpha_1)^{1/2}\Gamma_{12}^{5/2}}y_1^{n-2}
-\frac{\alpha_0\Gamma_1\Gamma_2L_{n-2,0}}{\alpha_1\alpha_2\Gamma_{12}^2}y_1^{n-2}
+O(y_1^{n-5})\nonumber\\
&\underset{y_1\to
  +\infty}{=}& 
y_1^{n-2}D_{n,0}+O(y_1^{n-5}),
\end{eqnarray}
with
\be\label{Dn0}
D_{n,0}&=&\frac{(1-d)\Gamma_1^2}{4\alpha_2^2\Gamma_{12}^2}\mathbf{1}_{\{n=3\}}-\frac{(R_2^2-R_1^2)\Gamma_1}{\alpha_2\Gamma_{12}}L_{n-1,0}-\frac{\Gamma_1^2}{2\alpha_2\Gamma_{12}^2}L_{n-2,0}\nonumber\\
&&-2\frac{\alpha_0(R_2^2-R_1^2)\Gamma_1\Gamma_2^{1/2}}{\alpha_2^2(2\alpha_1)^{1/2}\Gamma_{12}^{3/2}}N_{n-1,0}-\frac{\alpha_0\Gamma_1^2\Gamma_2^{1/2}}{\alpha_2^2(2\alpha_1)^{1/2}\Gamma_{12}^{5/2}}N_{n-2,0}.
\ee
The existence of an asymptotic expansion of $\delta_n(y_1)$ as $y_1\to
-\infty$ like the one given in (\ref{asymp-deln}) follows from
(\ref{def-deltan}) similarly as for the expansion at
$y_1=+\infty$. Moreover, like in (\ref{DLNp}), we obtain
\begin{eqnarray}\label{DLN}
\delta_n&\underset{y_1\to
  -\infty}{=}&\left(-d\lambda_{0}'\lambda_{n-3}+y_1\lambda_{0}'\lambda_{n-3}'\right)\mathbf{1}_{\{n=3\}}\nonumber\\
&&+2y_1\lambda_{-1}\lambda_{n-1}+2y_1\lambda_{0}\lambda_{n-2} +2(R_2^2-R_1^2)\lambda_0\lambda_{n-1}
-12\alpha_2\lambda_{-1}\lambda_0\lambda_{n-1}-6\alpha_2\lambda_{0}^2\lambda_{n-2}+O(y_1^{n-5})\nonumber\\
&=&\widetilde{D}_{n,0}y_1^{n-2}+O(y_1^{n-5}),
\end{eqnarray}
with
\be\label{Dn0t}
\widetilde{D}_{n,0}=\frac{1-d}{4\alpha_2^2}\mathbf{1}_{\{n=3\}}-\left(\widetilde{L}_{n-1,0}\frac{R_2^2-R_1^2}{\alpha_2}+\frac{\widetilde{L}_{n-2,0}}{2\alpha_2}\right).
\ee
Then, (\ref{asymp-Fn}) follows from (\ref{def-Fn}), (\ref{asymp-nu0}),
(\ref{asymp-deln}) and the recursion assumption. Moreover, 
\be\label{FD}
F_{n,0}=-\frac{4\alpha_0\alpha_2}{(R_2^2-R_1^2)^2}\left(\frac{\Gamma_2}{2\alpha_1\Gamma_{12}}\right)^{1/2}D_{n,0}.
\ee
Then, using the same trick as for $n=2$, we look for a solution
$\nu_n$ of (\ref{In}) under the form
$$\nu_n=\frac{F_{n,0}}{W(y_1)}y_1^{n-3/2}\Phi(y_1)+\widetilde{\nu_n}.$$
$\nu_n$ solves (\ref{In}) if and only if $\widetilde{\nu_n}$
solves
$$T\widetilde{\nu_n}=\widetilde{F_n},\quad{\rm where}\quad\widetilde{F_n}=F_n-F_{n,0}y_1^{n-3/2}\Phi(y_1)+4R_1^2\frac{d^2}{dy_1^2}\left[\frac{F_{n,0}y_1^{n-3/2}}{W(y_1)}\Phi(y_1)\right].$$
The function $\widetilde{F_n}$ defined just above admits expansions as
$y_1\to \pm\infty$ which are similar to those satisfied by $F_n$ and
given in
(\ref{asymp-Fn}), except that in the expansion of $\widetilde{F_n}$ as
$y_1\to +\infty$, the power of $y_1$ in the leading term is smaller from
three units than the one of $F_n$. By iterating this process
a finite number of times, we are brought back to solve an equation
like (\ref{In}), but with a right hand side which is in
$L^2(\R)$. Thanks to Lemma 2.1 in \cite{GP} and (\ref{asymp-W}), it turns out that $\nu_n$
satisfies (\ref{asymp-nun}), where the coefficient in the leading term
as $y_1\to +\infty$ is
\be\label{ND}
N_{n,0}=\frac{F_{n,0}}{2\Gamma_2}=-\frac{2\alpha_0\alpha_2}{(R_2^2-R_1^2)^2\left(2\alpha_1\Gamma_{12}\Gamma_2\right)^{1/2}}D_{n,0}.
\ee
Finally, from (\ref{IIn}), (\ref{asymp-nu0}), (\ref{asymp-nun}), (\ref{asymp-deln}) and (\ref{ND}), we deduce that $\lambda_n$ satisfies
(\ref{asymp-lamn}), with
\be\label{LD}
L_{n,0}&=&-2\frac{\alpha_0}{\alpha_2}\left(\frac{\Gamma_2}{2\alpha_1\Gamma_{12}}\right)^{1/2}N_{0,n}+\frac{2\alpha_2}{(R_2^2-R_1^2)^2}D_{n0}=\frac{2\alpha_2}{(R_2^2-R_1^2)^2\Gamma_{12}}D_{n0},
\ee
and
\be\label{Ln0t}
\widetilde{L}_{n,0}&=&\frac{2\alpha_2\widetilde{D}_{n,0}}{(R_2^2-R_1^2)^2},
\ee
which completes the recursion and proves that (\ref{asymp-deln}),
(\ref{asymp-Fn}), (\ref{asymp-nun}) and (\ref{asymp-lamn}) hold for
every $n\geq 1$. In addition, one can compute explicitely the
coefficients of the leading terms in the expansions of $\delta_n$,
$F_n$, $\nu_n$ and $\lambda_n$ as $y_1\to \pm\infty$. Indeed, as
$y_1\to +\infty$, according
to (\ref{Dn0}), (\ref{ND}) and (\ref{LD}), we have, for every $n\geq
3$,
\be\label{recD}
D_{n,0}&=&\frac{(1-d)\Gamma_1^2}{4\alpha_2^2\Gamma_{12}^2}\mathbf{1}_{\{n=3\}}-\frac{(R_2^2-R_1^2)\Gamma_1}{\alpha_2\Gamma_{12}}\frac{2\alpha_2}{(R_2^2-R_1^2)^2\Gamma_{12}}D_{n-1,0}-\frac{\Gamma_1^2}{2\alpha_2\Gamma_{12}^2}\frac{2\alpha_2}{(R_2^2-R_1^2)^2\Gamma_{12}}D_{n-2,0}\nonumber\\
&&
+2\frac{\alpha_0(R_2^2-R_1^2)\Gamma_1\Gamma_2^{1/2}}{\alpha_2^2(2\alpha_1)^{1/2}\Gamma_{12}^{3/2}}\frac{2\alpha_0\alpha_2}{(R_2^2-R_1^2)^2\left(2\alpha_1\Gamma_{12}\Gamma_2\right)^{1/2}}D_{n-1,0}\nonumber\\
&&+\frac{\alpha_0\Gamma_1^2\Gamma_2^{1/2}}{\alpha_2^2(2\alpha_1)^{1/2}\Gamma_{12}^{5/2}}\frac{2\alpha_0\alpha_2}{(R_2^2-R_1^2)^2\left(2\alpha_1\Gamma_{12}\Gamma_2\right)^{1/2}}D_{n-2,0}\nonumber\\
&=&\frac{(1-d)\Gamma_1^2}{4\alpha_2^2\Gamma_{12}^2}\mathbf{1}_{\{n=3\}}-\frac{2\Gamma_1}{(R_2^2-R_1^2)\Gamma_{12}^2}D_{n-1,0}-\frac{\Gamma_1^2}{(R_2^2-R_1^2)^2\Gamma_{12}^3}D_{n-2,0}\nonumber\\
&&
+\frac{2\alpha_0^2\Gamma_1}{\alpha_1\alpha_2\Gamma_{12}^{2}(R_2^2-R_1^2)}D_{n-1,0}+\frac{\alpha_0^2\Gamma_1^2}{\alpha_1\alpha_2\Gamma_{12}^{3}(R_2^2-R_1^2)^2}D_{n-2,0}\nonumber\\
&=&\frac{(1-d)\Gamma_1^2}{4\alpha_2^2\Gamma_{12}^2}\mathbf{1}_{\{n=3\}}-\frac{2\Gamma_1}{(R_2^2-R_1^2)\Gamma_{12}}D_{n-1,0}-\frac{\Gamma_1^2}{(R_2^2-R_1^2)^2\Gamma_{12}^2}D_{n-2,0}.\nonumber
\ee
Thus, defining for every $n\in \N$
\be\label{def-dn}
d_n & =&
\left(\frac{(R_2^2-R_1^2)\Gamma_{12}}{\Gamma_1}\right)^nD_{n0},
\ee
we infer thanks to (\ref{d10d11}) and (\ref{D20})
$$d_1=0,\quad
d_2=-(R_2^2-R_1^2)^2\frac{d\Gamma_{12}(R_2^2-R_1^2)+\Gamma_1R_1^2}{4\alpha_2^2\Gamma_1},\quad
d_3=(R_2^2-R_1^2)^2\frac{(1+d)\Gamma_{12}(R_2^2-R_1^2)+2\Gamma_1R_1^2}{4\alpha_2^2\Gamma_1},$$
and
$$\forall n\geq 4,\ d_n=-2d_{n-1}-d_{n-2}.$$
It follows that for $n\geq 2$,
$$d_n=-\frac{(-1)^{n}(R_2^2-R_1^2)^2\left(\Gamma_{12}(R_2^2-R_1^2)(d+n-2)+(n-1)\Gamma_1R_1^2\right)}{4\alpha_2^2\Gamma_1},$$
and therefore
$$D_{n0}=-\left(\frac{-\Gamma_1}{\Gamma_{12}(R_2^2-R_1^2)}\right)^{n}\frac{(R_2^2-R_1^2)^{2}\left(\Gamma_{12}(R_2^2-R_1^2)(n+d-2)+(n-1)\Gamma_1R_1^2\right)}{4\alpha_2^2\Gamma_1}.$$
Coming back to (\ref{ND}) and (\ref{LD}), we get, for $n\geq 3$,
\be\label{Nn0}
N_{n,0}=\left(\frac{-\Gamma_1}{\Gamma_{12}(R_2^2-R_1^2)}\right)^n\frac{\alpha_0(\Gamma_{12}(R_2^2-R_1^2)(n+d-2)+\Gamma_1R_1^2(n-1))}{2\alpha_2\Gamma_1\left(2\alpha_1\Gamma_{12}\Gamma_2\right)^{1/2}}
\ee
and
\be\label{Ln0}
L_{n,0}&=&-\left(\frac{-\Gamma_1}{\Gamma_{12}(R_2^2-R_1^2)}\right)^n\frac{\Gamma_{12}(R_2^2-R_1^2)(n+d-2)+\Gamma_1R_1^2(n-1)}{2\alpha_2\Gamma_1\Gamma_{12}}.
\ee

Similarly, as $y_1\to -\infty$, for $n\geq 3$, from (\ref{Dn0t}) and
(\ref{Ln0t}) we get
\begin{eqnarray*}
\widetilde{D}_{n,0}=\frac{1-d}{4\alpha_2^2}\mathbf{1}_{\{n=3\}}-\frac{2}{R_2^2-R_1^2}\widetilde{D}_{n-1,0}-\frac{1}{(R_2^2-R_1^2)^2}\widetilde{D}_{n-2,0}.
\end{eqnarray*}
Since from (\ref{asymp-deln}) for $n=1$ and (\ref{D20}), we have
$$\widetilde{D}_{1,0}=0,\quad \widetilde{D}_{2,0}=-\frac{d(R_2^2-R_1^2)+R_1^2}{4\alpha_2^2},$$
we deduce
$$\widetilde{D}_{3,0}=\frac{d(R_2^2-R_1^2)+R_1^2+R_2^2}{4\alpha_2^2(R_2^2-R_1^2)}$$
and for $n\geq 4$,
$$\widetilde{D}_{n,0}=\frac{(-1)^{n+1}}{4\alpha_2^2(R_2^2-R_1^2)^{n-2}}((d-2+n)(R_2^2-R_1^2)+(n-1)R_1^2),$$
and therefore thanks to (\ref{Ln0t}), we obtain
$$\widetilde{L}_{n,0}=\frac{(-1)^{n+1}}{2\alpha_2(R_2^2-R_1^2)^{n}}((d-2+n)(R_2^2-R_1^2)+(n-1)R_1^2).$$
The main results obtained in this section are summarized in the following proposition.
\begin{prop}\label{summary}
\begin{eqnarray*}
\lambda_{-1}(y_1)&=&\frac{R_2^2-R_1^2}{2\alpha_2}\\
\nu_0(y_1)&=&\frac{R_1^{1/3}\Gamma_2^{1/3}}{(2\alpha_1)^{1/2}\Gamma_{12}^{1/2}}\gamma_0\left(\frac{\Gamma_2^{1/3}y_1}{R_1^{2/3}}\right)\underset{y_1\to
  +\infty}{=}\left(\frac{\Gamma_2}{2\alpha_1\Gamma_{12}}\right)^{1/2}y_1^{1/2}+O(y_1^{-5/2})\\
\nu_0(y_1)&\underset{y_1\to-\infty}{\approx}&0\\
\lambda_{0}(y_1)&=&\frac{y_1}{2\alpha_2}-\frac{\alpha_0}{\alpha_2}\nu_{0}(y_1)^2\underset{y_1\to
  +\infty}{=}\frac{\Gamma_1}{2\alpha_2\Gamma_{12}}y_1+O(y_1^{-2})\\
\lambda_{0}(y_1)&\underset{y_1\to-\infty}{\approx}&\frac{y_1}{2\alpha_2}\\
\nu_1(y_1)&=&T^{-1}\left(-\frac{8\alpha_0\alpha_2R_1^2\lambda_0''\lambda_{-1}\nu_0}{(R_2^2-R_1^2)^2}-2d\nu_{0}'-4y_1\nu_{0}''\right)\nonumber\\
&\underset{y_1\to
  +\infty}{=}&\left\{\begin{array}{ll}\frac{1-d}{2(2\alpha_1\Gamma_2\Gamma_{12})^{1/2}}y_1^{-3/2}+O(y_1^{-9/2})
& {\rm
    if\ }d=2,3\\-\frac{6R_1^2}{(2\alpha_1\Gamma_{12}\Gamma_2)^{1/2}}\left(\frac{\alpha_0^2R_1^2}{\alpha_1\alpha_2(R_2^2-R_1^2)\Gamma_{12}}-\frac{5}{4\Gamma_2}\right)y_1^{-9/2}+O(y_1^{-15/2})&
   {\rm if\ }d=1\end{array}\right.\\
\nu_{1}(y_1)&\underset{y_1\to-\infty}{\approx}&0\\
\lambda_1(y_1) &=&-2\frac{\alpha_0}{\alpha_2}\nu_0\nu_1+\frac{4\alpha_2
  R_1^2}{(R_2^2-R_1^2)^2}\lambda_0''\lambda_{-1}\nonumber\\
&\underset{y_1\to
  +\infty}{=}&\left\{\begin{array}{ll}\frac{\alpha_0(d-1)}{2\alpha_1\alpha_2\Gamma_{12}}y_1^{-1}+O(y_1^{-4})& {\rm
    if\ }d=2,3\\
\frac{3\alpha_0R_1^2}{2\alpha_1\alpha_2\Gamma_{12}}\left(\frac{4R_1^2}{(R_2^2-R_1^2)\Gamma_{12}}-\frac{5}{\Gamma_2}\right)y_1^{-4}+O(y_1^{-7})& {\rm
    if\ }d=1.\end{array}\right.\\
\lambda_{1}(y_1)&\underset{y_1\to-\infty}{\approx}&0
\end{eqnarray*}
and for $n\geq 2$,
\begin{eqnarray*}
\nu_n(y_1)
&\underset{y_1\to
  +\infty}{=}&\left(\frac{-\Gamma_1}{\Gamma_{12}(R_2^2-R_1^2)}\right)^n\frac{\alpha_0(\Gamma_{12}(R_2^2-R_1^2)(n+d-2)+\Gamma_1R_1^2(n-1))}{2\alpha_2\Gamma_1\left(2\alpha_1\Gamma_{12}\Gamma_2\right)^{1/2}}y_1^{n-5/2}+O(y_1^{n-11/2})\nonumber\\
\nu_n(y_1)
&\underset{y_1\to
  -\infty}{\approx}&0\\
\lambda_n(y_1)&\underset{y_1\to
  +\infty}{=}&-\left(\frac{-\Gamma_1}{\Gamma_{12}(R_2^2-R_1^2)}\right)^n\frac{\Gamma_{12}(R_2^2-R_1^2)(n+d-2)+\Gamma_1R_1^2(n-1)}{2\alpha_2\Gamma_1\Gamma_{12}}y_1^{n-2}+O(y_1^{n-5})\\
\lambda_n(y_1)&\underset{y_1\to
  -\infty}{=}&\frac{(-1)^{n+1}}{2\alpha_2(R_2^2-R_1^2)^n}((d-2+n)(R_2^2-R_1^2)+(n-1)R_1^2)y_1^{n-2}+O(y_1^{n-5}).
\end{eqnarray*}
\end{prop}
\section{Truncation of the asymptotic expansions}\label{truncation}
In section \ref{formal}, we have explained how to calculate asymptotic expansions into powers of $\eps$ of $\omega$, $\tau$, $\nu$, 
$\lambda$ and $\mu$ in such a way that (\ref{ansatznorest}),
(\ref{apD0}), (\ref{apD1}), (\ref{apD2}) and (\ref{multiscale}) provide formally solutions to (\ref{sys}) at any order. However, we have not said anything about the convergence of these formal series. In this section, we prove that the truncations of the formal series
at a finite order provide approximate solutions to (\ref{sys}) at a
arbitrarily high order in terms of powers of $\eps$. More precisely,
$M$, $N$ and $L$ are three fixed positive integers, and we set in all
the section 
\be\label{troncat}
\omega(z)=\sum_{m=0}^M\eps^{2m}\omega_m(z),&& \tau(z)=\sum_{m=0}^M\eps^{2m}\tau_m(z),\nonumber\\
\nu(y_1)=\sum_{n=0}^N\eps^{2n/3}\nu_n(y_1),&&
\lambda(y_1)=\sum_{n=-1}^N\eps^{2n/3}\lambda_n(y_1),\nonumber\\
&&\mu(y_2)=\sum_{n=0}^L\eps^{2n/3}\mu_n(y_2),
\ee
where the $\omega_m$'s, $\tau_m$'s, $\nu_n$'s, $\lambda_n$'s and
$\mu_n$'s are the ones calculated in Section \ref{formal}. The way
integers $M$, $N$ and $L$ are chosen is explained in Sections
\ref{sec-omega-tau} and \ref{secmuasymp} below.

\subsection{Consistency of the ansatz}
Ansatz (\ref{ansatz}) requires the calculation of $\lambda(y_1)^{1/2}$
for $x\in {\rm Supp} \chi_\eps\subset D_1$.  So it makes sense to
combine (\ref{multiscale}) and (\ref{troncat}) only if the function
$\lambda$ given by (\ref{troncat}) satisfies $\lambda(y_1)\geq 0$ for $x\in   {\rm Supp} \chi_\eps$. We next show that the last inequality indeed holds for $x\in D_1$.

\begin{lem}\label{minor-lam}
Let $N>0$ and $\lambda$ given by (\ref{troncat}). There exists
$C>0$ (which might depend on $N$) such that for $\eps\in(0,1]$
    sufficiently small,
for every $x\in D_1$, $$\lambda(y_1)\geq C\eps^{-2/3}.$$
\end{lem}

\begin{Proof}
Let $x\in  D_1$. Then $y_2\geq
(R_2^2-R_1^2)/\eps^{2/3}-2\eps^{\beta-2/3}$, $-2\eps^{\beta-2/3}\leq
y_1\leq 2\eps^{\beta-2/3}$, and since $\gamma_0$ is increasing
  and $\gamma_0(y)\underset{y\to +\infty}{\sim} \sqrt{y}$, we get on the one side
$$\eps^{-2/3}\lambda_{-1}+\lambda_0(y_1)=\frac{y_2}{2\alpha_2}-\frac{\alpha_0}{\alpha_2}\nu_0(y_1)^2\geq\frac{R_2^2-R_1^2}{2\alpha_2\eps^{2/3}}-\frac{\eps^{\beta-2/3}}{\alpha_2}-\frac{\alpha_0}{\alpha_2}\nu_0(2\eps^{\beta-2/3})^2=\frac{R_2^2-R_1^2}{2\alpha_2\eps^{2/3}}+O(\eps^{\beta-2/3}),$$
whereas for $n\geq 1$, thanks to (\ref{asymp-lamn})
\be
|\eps^{2n/3}\lambda_n(y_1)| &\leq &
c_n\eps^{2n/3}\left(\mathbf{1}_{\{|y_1|\leq
  1\}}+|y_1|^{n-2}\mathbf{1}_{\{|y_1|\geq 1\}}\right)\leq
\tilde{c}_n\eps^{2n/3}\max(1,\eps^{(\beta-2/3)(n-2)})\nonumber\\
&\leq & \tilde{c}_n \max(\eps^{2n/3},\eps^{\beta(n-2)+4/3})=O(\eps^{2/3}),
\ee
for some $c_n>0$ and $\tilde{c}_n=2^{n-2}c_n$. As a result, for $\eps$
sufficiently small, we have $\lambda(y_1)\geq \frac{R_2^2-R_1^2}{4\alpha_2\eps^{2/3}}$ for every $x\in D_1$.
.
\end{Proof}

\subsection{Truncation of $(\omega,\tau)$ in $D_0$}
In this section, we prove that (\ref{troncat}) provides an approximate
solution to (\ref{sys}) in $D_0$ at an arbitrarily high order. For
convenience, we use the same notation $\omega$ for the functions
$z\mapsto \omega(z)$ and $x\mapsto \omega(z)=\omega(R_1^2-|x|^2)$.
\begin{lem}\label{est-eq-om-tau}
Let $M\geq 1$ be an integer, $\beta\in(0,2/3)$ and $\omega,\tau$ given by (\ref{troncat}).
Then
$$\left\|\eps^2\Delta\omega+\frac{\alpha_0}{\alpha_2}(R_2^2-R_1^2)\omega+z\omega-2\alpha_1\omega^3-2\alpha_0\tau^2\omega\right\|_{L^\infty(D_0)}=O(\eps^{(2-3\beta)M+2-3\beta/2})$$
and
$$\left\|\eps^2\Delta\tau+(R_2^2-R_1^2+z)\tau-2\alpha_2\tau^3-2\alpha_0\omega^2\tau\right\|_{L^\infty(D_0)}=O(\eps^{(2-3\beta)M+2-2\beta}).$$
\end{lem}

\begin{Proof}
Thanks to (\ref{eq-omega0}), (\ref{eq-tau0}) and (\ref{47}), we have
\be\label{f01-D0}
\lefteqn{\eps^2\Delta\omega+\frac{\alpha_0}{\alpha_2}(R_2^2-R_1^2)\omega+z\omega-2\alpha_1\omega^3-2\alpha_0\tau^2\omega}\nonumber\\
&=&\sum_{m=1}^{ M+1}\eps^{2m}\Delta\omega_{m-1}+\frac{\alpha_0}{\alpha_2}(R_2^2-R_1^2)\sum_{m=0}^{ M}\eps^{2m}\omega_m+z\sum_{m=0}^{ M}\eps^{2m}\omega_m.\nonumber\\
&&-2\alpha_1\sum_{m=0}^{3M}\eps^{2m}\sum_{\tiny\begin{array}{c}m_1+m_2+m_3=m\\0\leq
    m_1,m_2,m_3\leq M\end{array}}\omega_{m_1}\omega_{m_2}\omega_{m_3}-2\alpha_0\sum_{m=0}^{ 3M}\eps^{2m}\sum_{\tiny\begin{array}{c}m_1+m_2+m_3=m\\0\leq
    m_1,m_2,m_3\leq
    M\end{array}}\tau_{m_1}\tau_{m_2}\omega_{m_3}\nonumber\\
&=&\eps^{2(M+1)}\Delta\omega_{M}-2\alpha_1\sum_{m=M+1}^{ 3M}\eps^{2m}\sum_{\tiny\begin{array}{c}m_1+m_2+m_3=m\\0\leq
    m_1,m_2,m_3\leq
    M\end{array}}\omega_{m_1}\omega_{m_2}\omega_{m_3}\nonumber\\
&&-2\alpha_0\sum_{m=M+1}^{ 3M}\eps^{2m}\sum_{\tiny\begin{array}{c}m_1+m_2+m_3=m\\0\leq
    m_1,m_2,m_3\leq
    M\end{array}}\tau_{m_1}\tau_{m_2}\omega_{m_3}.
\ee
From Lemma \ref{est-omega-tau}, (\ref{eq-omega0}),
(\ref{eq-tau0}) and Remark \ref{rem-22}, we infer that for
every $x\in D_0$, 
\be\label{est-f01-D0}
\lefteqn{|\eps^2\Delta\omega+\frac{\alpha_0}{\alpha_2}(R_2^2-R_1^2)\omega+z\omega-2\alpha_1\omega^3-2\alpha_0\tau^2\omega|}\\
&\lesssim&\eps^{2(M+1)}z^{-3/2-3 M}+\sum_{m=
    M+1}^{ 3M}\eps^{2m}z^{3/2-3m}+\sum_{m= M+1}^{
  3M}\eps^{2m}z^{5/2-3m}\lesssim\eps^{(2-3\beta) M+2-3\beta/2}.\nonumber
\ee
Similarly,
\be\label{f02-D0}
\lefteqn{\eps^2\Delta\tau+(R_2^2-R_1^2+z)\tau-2\alpha_2\tau^3-2\alpha_0\omega^2\tau}\nonumber\\
&=&\sum_{m=1}^{M+1}\eps^{2m}\Delta\tau_{m-1}+(R_2^2-R_1^2+z)\sum_{m=0}^M\eps^{2m}\tau_m\nonumber\\
&&-2\alpha_2\sum_{m=0}^{3M}\eps^{2m}\sum_{\tiny\begin{array}{c}m_1+m_2+m_3=m\\0\leq
    m_1,m_2,m_3\leq
    M\end{array}}\tau_{m_1}\tau_{m_2}\tau_{m_3}-2\alpha_0\sum_{m=0}^{3M}\eps^{2m}\sum_{\tiny\begin{array}{c}m_1+m_2+m_3=m\\0\leq
    m_1,m_2,m_3\leq
    M\end{array}}\omega_{m_1}\omega_{m_2}\tau_{m_3}\nonumber\\
&=&\eps^{2(M+1)}\Delta\tau_M-2\alpha_2\sum_{m=M+1}^{3M}\eps^{2m}\sum_{\tiny\begin{array}{c}m_1+m_2+m_3=m\\0\leq
    m_1,m_2,m_3\leq
    M\end{array}}\tau_{m_1}\tau_{m_2}\tau_{m_3}\nonumber\\
&&-2\alpha_0\sum_{m=M+1}^{3M}\eps^{2m}\sum_{\tiny\begin{array}{c}m_1+m_2+m_3=m\\0\leq
    m_1,m_2,m_3\leq
    M\end{array}}\omega_{m_1}\omega_{m_2}\tau_{m_3},
\ee
thus for $x\in D_0$,
\be\label{est-f02-D0}
\lefteqn{|\eps^2\Delta\tau+(R_2^2-R_1^2+z)\tau-2\alpha_2\tau^3-2\alpha_0\omega^2\tau|}\nonumber\\&\lesssim&\eps^{(2-3\beta)M+2-\beta}+\sum_{m=M+1}^{3M}\eps^{2m}z^{1-3m}
\lesssim\eps^{(2-3\beta)M+2-2\beta}.
\ee
\end{Proof}

\subsection{Truncation of $(\eps^{1/3}\nu,\eps^{1/3}\lambda^{1/2})$ in $D_1$}
\begin{lem}\label{lem-tronc-nu-l}
Let $N\geq 4$ be an integer, and $\nu,\lambda$ given by (\ref{troncat}).
Then
\be
\left\|\eps^2\Delta\left(\eps^{1/3}\nu\right)+\left(\frac{\alpha_0}{\alpha_2}(R_2^2-R_1^2)+z\right)\eps^{1/3}\nu-2\alpha_1\left(\eps^{1/3}\nu\right)^3-2\alpha_0\left(\eps^{1/3}\lambda^{1/2}\right)^2\eps^{1/3}\nu\right\|_{L^\infty(D_1)}=O(\eps^{\beta N+4-7\beta/2})\nonumber
\ee
and
\be
\left\|\eps^2\Delta\left(\eps^{1/3}\lambda^{1/2}\right)+(R_2^2-R_1^2+z)\eps^{1/3}\lambda^{1/2}-2\alpha_2\left(\eps^{1/3}\lambda^{1/2}\right)^3-2\alpha_0\left(\eps^{1/3}\nu\right)^2\eps^{1/3}\lambda^{1/2}\right\|_{L^\infty(D_1)}=O(\eps^{\beta N+2-\beta}).\nonumber
\ee
\end{lem}

\begin{Proof}
Using (\ref{troncat}), (\ref{I0}) and (\ref{preIn}) for $n\in
\{1,\cdots, N\}$, we get
\be\label{g1-D1}
\lefteqn{\eps^{-1}\left(\eps^2\Delta\left(\eps^{1/3}\nu\right)+\left(\frac{\alpha_0}{\alpha_2}(R_2^2-R_1^2)+z\right)\eps^{1/3}\nu-2\alpha_1\left(\eps^{1/3}\nu\right)^3-2\alpha_0\left(\eps^{1/3}\lambda^{1/2}\right)^2\eps^{1/3}\nu\right)}\nonumber\\
&=&\eps^{4/3}\Delta\nu+\frac{\alpha_0}{\alpha_2}\frac{R_2^2-R_1^2}{\eps^{2/3}}\nu+y_1\nu-2\alpha_1\nu^3-2\alpha_0\lambda\nu\nonumber\\
&=&-2d\eps^{2/3}\nu'+4R_1^2\nu''-4\eps^{2/3}y_1\nu''
+\frac{\alpha_0}{\alpha_2}\frac{R_2^2-R_1^2}{\eps^{2/3}}\nu+y_1\nu
-2\alpha_1\nu^3-2\alpha_0\lambda\nu\nonumber\\
&=&-2d\sum_{n=1}^{N+1}\eps^{2n/3}\nu_{n-1}'+4R_1^2\sum_{n=0}^N\eps^{2n/3}\nu_n''-4y_1\sum_{n=1}^{N+1}\eps^{2n/3}\nu_{n-1}''+\frac{\alpha_0}{\alpha_2}(R_2^2-R_1^2)\sum_{n=-1}^{N-1}\eps^{2n/3}\nu_{n+1}\nonumber\\
&&+y_1\sum_{n=0}^N\eps^{2n/3}\nu_n-2\alpha_1\sum_{n=0}^{3N}\eps^{2n/3}\sum_{\tiny\begin{array}{c}n_1+n_2+n_3=n\\0\leq
    n_1,n_2,n_3\leq
N\end{array}}\nu_{n_1}\nu_{n_2}\nu_{n_3}-2\alpha_0\sum_{n=-1}^{2N}\eps^{2n/3}\sum_{\tiny\begin{array}{c}n_1+n_2=n\\-1\leq
    n_1\leq
N\\0\leq n_2\leq N\end{array}}\lambda_{n_1}\nu_{n_2}\nonumber\\
&=&-2d\sum_{n=1}^{N+1}\eps^{2n/3}\nu_{n-1}'+4R_1^2\sum_{n=0}^N\eps^{2n/3}\nu_n''-4y_1\sum_{n=1}^{N+1}\eps^{2n/3}\nu_{n-1}''\nonumber\\
&&+y_1\sum_{n=0}^N\eps^{2n/3}\nu_n-2\alpha_1\sum_{n=0}^{3N}\eps^{2n/3}\sum_{\tiny\begin{array}{c}n_1+n_2+n_3=n\\0\leq
    n_1,n_2,n_3\leq
N\end{array}}\nu_{n_1}\nu_{n_2}\nu_{n_3}-2\alpha_0\sum_{n=0}^{2N}\eps^{2n/3}\sum_{\tiny\begin{array}{c}n_1+n_2=n\\0\leq
    n_1,n_2\leq
N\end{array}}\lambda_{n_1}\nu_{n_2}\nonumber\\
&=&-2d\eps^{2(N+1)/3}\nu_{N}'-4y_1\eps^{2(N+1)/3}\nu_{N}''\nonumber\\
&&-2\alpha_1\sum_{n=N+1}^{3N}\eps^{2n/3}\sum_{\tiny\begin{array}{c}n_1+n_2+n_3=n\\0\leq
    n_1,n_2,n_3\leq
N\end{array}}\nu_{n_1}\nu_{n_2}\nu_{n_3}-2\alpha_0\sum_{n=N+1}^{2N}\eps^{2n/3}\sum_{\tiny\begin{array}{c}n_1+n_2=n\\0\leq
    n_1,n_2\leq
N\end{array}}\lambda_{n_1}\nu_{n_2}\\
&=&\eps^{2(N+1)/3}\left[-2d\nu_{N}'-4y_1\nu_{N}''-2\alpha_1\sum_{n=0}^{2N-1}\eps^{2n/3}\sum_{\tiny\begin{array}{c}n_1+n_2+n_3=n+N+1\\0\leq
    n_1,n_2,n_3\leq
N\end{array}}\nu_{n_1}\nu_{n_2}\nu_{n_3}\right.\nonumber\\
&&\left.-2\alpha_0\sum_{n=0}^{N-1}\eps^{2n/3}\sum_{\tiny\begin{array}{c}n_1+n_2=n+N+1\\0\leq
    n_1,n_2\leq
N\end{array}}\lambda_{n_1}\nu_{n_2}\right].
\nonumber
\ee
Thus, if we note that for $x\in
D_1=\left\{x\in\R^d|-2\eps^{\beta-2/3}\leq y_1\leq
2\eps^{\beta-2/3}\right\}$, $|\eps^{2/3}y_1|\lesssim\eps^\beta\to 0$
as $\eps\to 0$, we have thanks to (\ref{asymp-nun}), (\ref{asymp-lamn}) and (\ref{asymp-nu0})
\be\label{est-g1-D1}
\lefteqn{\eps\left|\eps^{4/3}\Delta\nu+\frac{\alpha_0}{\alpha_2}\frac{R_2^2-R_1^2}{\eps^{2/3}}\nu+y_1\nu-2\alpha_1\nu^3-2\alpha_0\lambda\nu\right|}\nonumber\\
&\lesssim&\eps^{2N/3+5/3}\max(1,y_1)^{N-7/2}+\eps^{2N/3+5/3}\sum_{n=0}^{2N-1}\eps^{2n/3}\max(1,y_1)^{n+N+1-5/2-5/2+1/2}\nonumber\\
&&+\eps^{2N/3+5/3}\sum_{n=0}^{N-1}\eps^{2n/3}\max(1,y_1)^{n+N+1-5/2-2}\nonumber\\
&\lesssim&\eps^{2N/3+5/3}\max(1,y_1)^{N-7/2}\lesssim \eps^{2N/3+5/3}\max(1,\eps^{(\beta-2/3)(N-7/2)})=\eps^{\beta N+4-7\beta/2},
\ee
where the last equality holds because $N\geq 4$. The first estimate of
the lemma is proved.
Similarly, from (\ref{troncat}), (\ref{I-1}), (\ref{II0}) and
(\ref{star}), we deduce
\be\label{h1-D1}
\lefteqn{\eps^2\Delta\left(\eps^{1/3}\lambda^{1/2}\right)+(R_2^2-R_1^2+z)\eps^{1/3}\lambda^{1/2}-2\alpha_2\left(\eps^{1/3}\lambda^{1/2}\right)^3-2\alpha_0\left(\eps^{1/3}\nu\right)^2\eps^{1/3}\lambda^{1/2}}\nonumber\\
&=&\eps\lambda^{-3/2}\left[-d\eps^{2/3}\lambda\lambda'-(R_1^2-\eps^{2/3}y_1){\lambda'}^2+2(R_1^2-\eps^{2/3}y_1)\lambda\lambda''+y_2\lambda^{2}-2\alpha_2\lambda^{3}-2\alpha_0\nu^2\lambda^{2}\right]\nonumber\\
&=&\eps^{-1/3}\lambda^{-3/2}\left[-d\sum_{n=1}^{2N+3}\eps^{2n/3}\sum_{\tiny\begin{array}{c}n_1+n_2=n-3\\-1\leq
    n_1,n_2\leq
    N\end{array}}\lambda_{n_1}'\lambda_{n_2}
+2R_1^2\sum_{n=0}^{2N+2}\eps^{2n/3}\sum_{\tiny\begin{array}{c}n_1+n_2=n-2\\-1\leq
    n_1,n_2\leq
    N\end{array}}\lambda_{n_1}''\lambda_{n_2}\right.\nonumber\\
&&-2y_1\sum_{n=1}^{2N+3}\eps^{2n/3}\sum_{\tiny\begin{array}{c}n_1+n_2=n-3\\-1\leq
    n_1,n_2\leq
    N\end{array}}\lambda_{n_1}''\lambda_{n_2}
-R_1^2\sum_{n=0}^{2N+2}\eps^{2n/3}\sum_{\tiny\begin{array}{c}n_1+n_2=n-2\\-1\leq
    n_1,n_2\leq
    N\end{array}}\lambda_{n_1}'\lambda_{n_2}'\nonumber\\
&&+y_1\sum_{n=1}^{2N+3}\eps^{2n/3}\sum_{\tiny\begin{array}{c}n_1+n_2=n-3\\-1\leq
    n_1,n_2\leq
    N\end{array}}\lambda_{n_1}'\lambda_{n_2}'
+(R_2^2-R_1^2)\sum_{n=-1}^{2N+1}\eps^{2n/3}\sum_{\tiny\begin{array}{c}n_1+n_2=n-1\\-1\leq
    n_1,n_2\leq
    N\end{array}}\lambda_{n_1}\lambda_{n_2}\nonumber\\
&&+y_1\sum_{n=0}^{2N+2}\eps^{2n/3}\sum_{\tiny\begin{array}{c}n_1+n_2=n-2\\-1\leq
    n_1,n_2\leq
    N\end{array}}\lambda_{n_1}\lambda_{n_2}
-2\alpha_2\sum_{n=-1}^{3N+2}\eps^{2n/3}\sum_{\tiny\begin{array}{c}n_1+n_2+n_3=n-2\\-1\leq
  n_1,n_2,n_3\leq N\end{array}}\lambda_{n_1}\lambda_{n_2}\lambda_{n_3}\nonumber\\
&&\left.-2\alpha_0\sum_{n=0}^{4N+2}\eps^{2n/3}\sum_{\tiny\begin{array}{c}n_1+n_2+n_3+n_4=n-2\\-1\leq
  n_1,n_2\leq N\\
0\leq n_3,n_4\leq
N\end{array}}\lambda_{n_1}\lambda_{n_2}\nu_{n_3}\nu_{n_4}\right]\nonumber\\
&=&\eps^{-1/3}\lambda^{-3/2}\left[-d\sum_{n=N+1}^{2N+3}\eps^{2n/3}\sum_{\tiny\begin{array}{c}n_1+n_2=n-3\\-1\leq
    n_1,n_2\leq
    N\end{array}}\lambda_{n_1}'\lambda_{n_2}
+2R_1^2\sum_{n=N+1}^{2N+2}\eps^{2n/3}\sum_{\tiny\begin{array}{c}n_1+n_2=n-2\\-1\leq
    n_1,n_2\leq
    N\end{array}}\lambda_{n_1}''\lambda_{n_2}\right.\nonumber\\
&&-2y_1\sum_{n=N+1}^{2N+3}\eps^{2n/3}\sum_{\tiny\begin{array}{c}n_1+n_2=n-3\\-1\leq
    n_1,n_2\leq
    N\end{array}}\lambda_{n_1}''\lambda_{n_2}
-R_1^2\sum_{n=N+1}^{2N+2}\eps^{2n/3}\sum_{\tiny\begin{array}{c}n_1+n_2=n-2\\-1\leq
    n_1,n_2\leq
    N\end{array}}\lambda_{n_1}'\lambda_{n_2}'\nonumber\\
&&+y_1\sum_{n=N+1}^{2N+3}\eps^{2n/3}\sum_{\tiny\begin{array}{c}n_1+n_2=n-3\\-1\leq
    n_1,n_2\leq
    N\end{array}}\lambda_{n_1}'\lambda_{n_2}'
+(R_2^2-R_1^2)\sum_{n=N+1}^{2N+1}\eps^{2n/3}\sum_{\tiny\begin{array}{c}n_1+n_2=n-1\\-1\leq
    n_1,n_2\leq
    N\end{array}}\lambda_{n_1}\lambda_{n_2}\nonumber\\
&&+y_1\sum_{n=N+1}^{2N+2}\eps^{2n/3}\sum_{\tiny\begin{array}{c}n_1+n_2=n-2\\-1\leq
    n_1,n_2\leq
    N\end{array}}\lambda_{n_1}\lambda_{n_2}
-2\alpha_2\sum_{n=N+1}^{3N+2}\eps^{2n/3}\sum_{\tiny\begin{array}{c}n_1+n_2+n_3=n-2\\-1\leq
  n_1,n_2,n_3\leq N\end{array}}\lambda_{n_1}\lambda_{n_2}\lambda_{n_3}\nonumber\\
&&\left.-2\alpha_0\sum_{n=N+1}^{4N+2}\eps^{2n/3}\sum_{\tiny\begin{array}{c}n_1+n_2+n_3+n_4=n-2\\-1\leq
  n_1,n_2\leq N\\
0\leq n_3,n_4\leq N\end{array}}\lambda_{n_1}\lambda_{n_2}\nu_{n_3}\nu_{n_4}\right]
\ee
In order to estimate this quantity, we consider separately each sum appearing in the bracket in
the right hand side of (\ref{h1-D1}). Let us focus for instance on
the first one. If $n\geq N+1$,
$n_1+n_2=n-3$ and $n_1,n_2\geq 1$, then we infer from
(\ref{asymp-lamn}) that for $x\in D_1$ (which implies
$|y_1|\lesssim\eps^{\beta-2/3}$), we have
\begin{eqnarray*}
\eps^{2n/3}|\lambda_{n_1}'\lambda_{n_2}| &\lesssim&\eps^{2n/3}\max(1,|y_1|)^{n-8}\lesssim\max(\eps^{2n/3},\eps^{2n/3+(\beta-2/3)(n-8)})=\max(\eps^{2n/3},\eps^{\beta
  n+8(2/3-\beta)})\\
&\lesssim&\max(\eps^{2(N+1)/3},\eps^{\beta
  (N+1)+8(2/3-\beta)}).
\end{eqnarray*}
If one of the two indices $n_1,n_2$
belongs to $\{-1,0\}$, whereas the other one is larger than or equal
to 1, we infer similarly thanks to (\ref{I-1}), (\ref{II0}) and
(\ref{asymp-lamn}) that  
$$\eps^{2n/3}|\lambda_{n_1}'\lambda_{n_2}|\lesssim\max(\eps^{2(N+1)/3},\eps^{\beta
  (N+1)+5(2/3-\beta)}).$$
Finally, If $N\geq 3$, the conditions $n_1+n_2=n-3\geq N-2$ excludes
the case where both $n_1$ and $n_2$ belong to $\{-1,0\}$. Using
similar arguments as well as Lemma \ref{minor-lam}, we deduce that for
$x\in D_1$ and $N$ large enough,
\be\label{est-h1-D1}
\lefteqn{\left|\eps^2\Delta\left(\eps^{1/3}\lambda^{1/2}\right)+(R_2^2-R_1^2+z)\eps^{1/3}\lambda^{1/2}-2\alpha_2\left(\eps^{1/3}\lambda^{1/2}\right)^3-2\alpha_0\left(\eps^{1/3}\nu\right)^2\eps^{1/3}\lambda^{1/2}\right|}\nonumber\\
&\lesssim &\eps^{-1/3}\eps\max(\eps^{2(N+1)/3},\eps^{\beta
  (N+1)+2(2/3-\beta)})\lesssim\eps^{\beta N+2-\beta}.\ \ \ \ \ \ \ \ \ \ \ \ \ \ \ \ \ \ \  \ \ \ \ \ \ \ \  \ \ \ \ \ \ \  \ \ \ \ \ \ \ \  \ \ 
\ee
\end{Proof}

\subsection{Truncation of $(0,\eps^{1/3}\mu)$ in $D_2$}
\begin{lem}\label{lem-est-h2-D2}
Let $L\geq 1$ be an integer and $\mu$ be given by (\ref{troncat}).
There exists $C>0$ such that for $x\in\R^d$ and $\eps\in ]0,1]$,
\be\label{pre-est-h2-D2}
\left|\eps^2\Delta\left(\eps^{1/3}\mu\right)+(R_2^2-R_1^2+z)\eps^{1/3}\mu-2\alpha_2\left(\eps^{1/3}\mu\right)^3\right|\leq \frac{C\eps^{2L/3+5/3}}{1+|y_2|^{2L+1/2}},
\ee
where $y_2=(R_2^2-|x|^2)/\eps^{2/3}$.
\end{lem}

\begin{cor}\label{cor1}
Under the same assumptions, there is $h\in L^2\cap L^\infty(\R^d)$ such that for
every $x\in\R^d$ and $\eps\in ]0,1]$,
\be\label{est-h2-D2}
\left|\eps^2\Delta\left(\eps^{1/3}\mu\right)+(R_2^2-R_1^2+z)\eps^{1/3}\mu-2\alpha_2\left(\eps^{1/3}\mu\right)^3\right|\leq \eps^{2L/3+5/3}h(x).
\ee
\end{cor}

\begin{cor}\label{cor2}
Under the same assumptions, there is $C>0$ such that for $x\in D_1\cap
D_2$ and $\eps\in ]0,1]$,
\be\label{post-est-h2-D2}
\left|\eps^2\Delta\left(\eps^{1/3}\mu\right)+(R_2^2-R_1^2+z)\eps^{1/3}\mu-2\alpha_2\left(\eps^{1/3}\mu\right)^3\right|\leq C\eps^{2L+2}.
\ee
\end{cor}

\paragraph{Proof of Lemma \ref{lem-est-h2-D2}.}
Taking into account the equations satisfied by the $\mu_n$'s, namely
\be\label{eq-mu0}
4R_2^2\mu_0''+y_2\mu_0-2\alpha_2\mu_0^3=0
\ee
for $n=0$ and
\be\label{eq-mun}
4R_2^2\mu_n''=2\alpha_2\sum_{n_1+n_2+n_3=n}\mu_{n_1}\mu_{n_2}\mu_{n_3}+2d\mu_{n-1}'+4y_2\mu_{n-1}''-y_2\mu_n
\ee
for $n\geq 1$, we infer
\be\label{def-g2}
\lefteqn{\eps^2\Delta\left(\eps^{1/3}\mu\right)+(R_2^2-R_1^2+z)\eps^{1/3}\mu-2\alpha_2\left(\eps^{1/3}\mu\right)^3}\nonumber\\
&=& \eps(\eps^{4/3}\Delta \mu+y_2\mu-2\alpha_2\mu^3)\nonumber\\
&=&\eps
\left(-2d\eps^{2(L+1)/3}\mu_L'-4y_2\eps^{2(L+1)/3}\mu_L''-2\alpha_2\sum_{n=L+1}^{3L}\eps^{2n/3}\sum_{\tiny\begin{array}{c}n_1+n_2+n_3=n\\0\leq
    n_1,n_2,n_3\leq
    L \end{array}}\mu_{n_1}\mu_{n_2}\mu_{n_3}\right)\nonumber\\
&=&\eps^{2L/3+5/3}\left(-2d\mu_L'-4y_2\mu_L''-2\alpha_2\sum_{n=0}^{2L-1}\eps^{2n/3}\sum_{\tiny\begin{array}{c}n_1+n_2+n_3=n+L+1\\0\leq
    n_1,n_2,n_3\leq
    L \end{array}}\mu_{n_1}\mu_{n_2}\mu_{n_3}\right).
\ee
Let us define for $y\in \R$
$$h_0(y)=(1+|y|^{2L+1/2})\max\left(|\mu_L'(y)|,|y\mu_L''(y)|,\underset{\tiny\begin{array}{c}0\leq
    n\leq2L-1\\n_1+n_2+n_3=n+L+1\\0\leq n_1,n_2,n_3\leq L\end{array}}{\max}|\mu_{n_1}(y)\mu_{n_2}(y)\mu_{n_3}(y)|\right).$$
Thanks to (\ref{mun}) and Propositions
\ref{proposition-Painleve} and \ref{proposition-mun}, $h_0$ is
uniformly bounded on $\R$. The lemma follows.\hfill$\rule{.3\baselineskip}{.35\baselineskip}$

\paragraph{Proof of Corollary \ref{cor1}.}
For $x\in \R^d$ and $\eps\leq 1$, one has 
$$\frac{1}{1+|y_2|^{2L+1/2}}=\frac{1}{1+\left(\frac{R_2^2-|x|^2)}{\eps^{2/3}}\right)^{2L+1/2}}\leq
h(x)=\left\{\begin{array}{lll} 1&{\rm if} &|x|^2\leq
2R_2^2\\ \frac{1}{1+\left(\frac{|x|^2}{2}\right)^{2L+1/2}}&{\rm if} &|x|^2\geq
2R_2^2\end{array}\right.$$
The corollary follows, since $L\geq 1$ and $d\leq 3$ imply $h\in
L^2(\R^d)$.\hfill$\rule{.3\baselineskip}{.35\baselineskip}$

\paragraph{Proof of Corollary \ref{cor2}.}
The corollary follows from Lemma \ref{pre-est-h2-D2} and from the
inequality
$$\frac{1}{1+|y_2|^{2L+1/2}}\lesssim \eps^{4L/3+1/3},$$
that holds for $x\in D_1\cap D_2$.\hfill$\rule{.3\baselineskip}{.35\baselineskip}$

\subsection{Comparison of $(\omega,\tau)$ and
  $\eps^{1/3}(\nu,\lambda^{1/2})$ in $D_0\cap D_1$}\label{sec-omega-tau}
\begin{lem}\label{lem-omega-tau}
Let $N\in\N^*$, $M\geq\frac{\beta}{2-3\beta}N$, and $\omega$, $\tau$ given by (\ref{troncat}). 
Then for every $l\geq 0$,
\be\label{est-om-12}
\left\|\omega^{(l)}-\sum_{\substack{(m,n)\in\N^2\\(2-3\beta)m+\beta
    n\leq\beta N}}
  \eps^{2m}w_{m,n}\frac{d^l}{dz^l}\left(z^{1/2-3m+n}\right)\right\|_{L^\infty(D_0\cap
  D_1)}\underset{\eps\to 0}{=}o\left(\eps^{\beta(N+1/2-l)}\right)
\ee
and
\be\label{est-tau-12}
\left\|\tau^{(l)}-\lambda_{-1}^{1/2}\mathbf{1}_{l=0}-\sum_{\substack{(m,n)\in\N^2\\ (2-3\beta)m+\beta
    n\leq \beta N}}\eps^{2m}t_{m,n}\ \frac{d^l}{dz^l}\left(z^{1+n-3m}\right)\right\|_{L^\infty(D_0\cap
  D_1)}\underset{\eps\to 0}{=}o(\eps^{\beta(N+1-l)}).
\ee
where the $w_{m,n}$'s and the $t_{m,n}$'s are defined in
Lemma \ref{est-omega-tau} and (\ref{asymp-tau0}).
\end{lem}

\begin{Proof}
From Lemma \ref{est-omega-tau}, for every $l\geq 0$,
\be
\omega^{(l)}(z)=\sum_{m=0}^M\eps^{2m}\omega_m^{(l)}(z)&\underset{z\to
  0}{\approx}&\sum_{m=0}^M\eps^{2m}\sum_{n=0}^\infty
w_{m,n}\frac{d^l}{dz^l}\left(z^{1/2-3m+n}\right)\nonumber\\
&\underset{z\to
  0}{\approx}&\sum_{k=-3M}^\infty\ \sum_{\substack{(m,n)\in\{0,\cdots,M\}\times\N\\ n-3m=k}}\eps^{2m}w_{m,n}\ \frac{d^l}{dz^l}\left(z^{1/2+k}\right).
\ee
Thus, since $x\in D_0\cap D_1$ implies $\eps^{\beta}\leq z\leq
2\eps^{\beta}\to 0$ as $\eps\to 0$,
\be
\omega^{(l)}(z)&\underset{z\to
  0}{=}&\sum_{k=-3M}^N\big(\sum_{\substack{(m,n)\in\{0,\cdots,M\}\times\N\\ n-3m=k}}\eps^{2m}w_{m,n}\big)\frac{d^l}{dz^l}\left(z^{1/2+k}\right)+o(z^{1/2+N-l})\nonumber\\
&\underset{z\to
  0}{=}&\sum_{\substack{(m,n)\in\{0,\cdots,M\}\times\N\\ n-3m\leq
    N}}\eps^{2m}w_{m,n}\frac{d^l}{dz^l}\left(z^{1/2+n-3m}\right)+o(z^{1/2+N-l})\nonumber\\
&\underset{z\to
  0}{=}&\sum_{\substack{(m,n)\in\{0,\cdots,M\}\times\N\\ n-3m\leq
    N\\2m+\beta(1/2+n-3m-l)\leq\beta(1/2+N-l)}}\eps^{2m}w_{m,n}\frac{d^l}{dz^l}\left(z^{1/2+n-3m}\right)+o_{L^\infty(D_0\cap
  D_1)}(\eps^{\beta(1/2+N-l)})\nonumber\\
&\underset{z\to
  0}{=}&\sum_{\substack{(m,n)\in\N^2\\ (2-3\beta)m+\beta
    n\leq\beta N}}\eps^{2m}w_{m,n}\frac{d^l}{dz^l}\left(z^{1/2+n-3m}\right)+o_{L^\infty(D_0\cap
  D_1)}(\eps^{\beta(1/2+N-l)}).
\ee
Note that the assumption on $M$ in the statement of the Lemma ensures
that the set $\{(m,n)\in \N^2, (2-3\beta)m+\beta n\leq \beta N\}$ is a
  triangle included in the rectangle $\{0,\cdots, M\}\times
    \{0,\cdots,N\}$.  
Similarly, we infer from Lemma \ref{est-omega-tau} and
(\ref{asymp-tau0}) that
\be
\tau^{(l)}(z)=\sum_{m=0}^M\eps^{2m}\tau_m^{(l)}(z) &\underset{z\to
  0}{\approx}&\lambda_{-1}^{1/2}\mathbf{1}_{l=0}+\sum_{m=0}^M\eps^{2m}\sum_{n=0}^{\infty}t_{m,n}\frac{d^l}{dz^l}\left(z^{1-3m+n}\right)\nonumber\\
&\underset{z\to
  0}{\approx}&\lambda_{-1}^{1/2}\mathbf{1}_{l=0}+\sum_{k=-3M}^\infty\ \sum_{\substack{(m,n)\in\{0,\cdots,M\}\times\N\\ n-3m=k}}\eps^{2m}t_{m,n}\ \frac{d^l}{dz^l}\left(z^{1+k}\right).\nonumber\\
&\underset{z\to
  0}{=}&\lambda_{-1}^{1/2}\mathbf{1}_{l=0}+\sum_{k=-3M}^N\ \sum_{\substack{(m,n)\in\{0,\cdots,M\}\times\N\\ n-3m=k}}\eps^{2m}t_{m,n}\ \frac{d^l}{dz^l}\left(z^{1+k}\right)+o(z^{N+1-l}).\nonumber
\ee
Thus,
\be
\tau^{(l)}(z)&\underset{\eps\to
  0}{=}&\lambda_{-1}^{1/2}\mathbf{1}_{l=0}+\sum_{\substack{(m,n)\in\{0,\cdots,M\}\times\N\\ 2m+(1+n-3m-l)\beta\leq (N+1-l)\beta}}\eps^{2m}t_{m,n}\ \frac{d^l}{dz^l}\left(z^{1+n-3m}\right)+o_{L^\infty(D_0\cap
  D_1)}(\eps^{\beta(N+1-l)})\nonumber\\
&\underset{\eps\to
  0}{=}&\lambda_{-1}^{1/2}\mathbf{1}_{l=0}+\sum_{\substack{(m,n)\in\N^2\\ (2-3\beta)m+\beta
    n\leq \beta N}}\eps^{2m}t_{m,n}\ \frac{d^l}{dz^l}\left(z^{1+n-3m}\right)+o_{L^\infty(D_0\cap
  D_1)}(\eps^{\beta(N+1-l)}).
\ee
\end{Proof}

\begin{lem}\label{lem-lam-nu}
Let $N\geq 1$. We assume that $\beta\in (0,2/3)\backslash \Q$. There exist two families of numbers
$(n_{m,n})_{m\geq 0,n\geq 0}$ and $(l_{m,n})_{m\geq 0,n\geq 0}$ which do not depend on $N$ such
that if $\nu$ and
$\lambda$ are given by (\ref{troncat}), then for $l=0,1,2$,
\be
\left\|\frac{d^l}{dz^l}\left(\eps^{1/3}\nu(y_1)\right)-\sum_{\substack{(m,n)\in\N^2\\(2-3\beta)m+\beta
    n\leq\beta N}}
  \eps^{2m}n_{m,n}\frac{d^l}{dz^l}\left(z^{1/2-3m+n}\right)\right\|_{L^\infty(D_0\cap
  D_1)}\underset{\eps\to 0}{=}o\left(\eps^{\beta(N+1/2-l)}\right)
\ee
and
\be\label{d2}
\left\|\frac{d^l}{dz^l}\left(\eps^{1/3}\lambda(y_1)^{1/2}\right)-\lambda_{-1}^{1/2}\mathbf{1}_{l=0}-\!\!\!\!\!\!\!\!\!\!\!\sum_{\substack{(m,n)\in\N^2\\ (2-3\beta)m+\beta
    n\leq \beta N}}\!\!\!\!\eps^{2m}l_{m,n}\ \frac{d^l}{dz^l}\left(z^{1+n-3m}\right)\right\|_{L^\infty(D_0\cap
  D_1)}\!\!\!\!\!\!\!\!\!\!\!&\underset{\eps\to
  0}{=}&\!\!\!o(\eps^{\beta(N+1-l)}).\ \ \ \ 
\ee
\end{lem}

\begin{Proof}
For $x\in D_0\cap D_1$, we have $2\eps^{\beta-2/3}\geq y_1\geq\eps^{\beta-2/3}\to +\infty$
as $\eps\to 0$. Thus, we infer from (\ref{asymp-nun}) and
(\ref{asymp-nu0}) that for every $l\geq 0$,
\be
\lefteqn{\frac{d^l}{dy_1^l}\left(\eps^{1/3}\nu(y_1)\right)}\nonumber\\
&\underset{y_1\to +\infty}{\approx}&\eps^{1/3}\sum_{m=0}^\infty
N_{0,m}\frac{d^l}{dy_1^l}\left(y_1^{1/2-3m}\right)+\eps^{1/3}\sum_{n=1}^N\eps^{2n/3}\sum_{m=0}^\infty
N_{n,m}\frac{d^l}{dy_1^l}\left(y_1^{n-5/2-3m}\right)\nonumber\\
&\underset{y_1\to +\infty}{\approx}&\sum_{m=0}^\infty
N_{0,m}\eps^{1/3}\frac{d^l}{dy_1^l}\left(y_1^{1/2-3m}\right)+\sum_{k=-N}^\infty\ \sum_{\substack{(n,m)\in\{1,\cdots,N\}\times\N\\ 3m-n=k}}N_{n,m}\eps^{1/3+2n/3}\frac{d^l}{dy_1^l}\left(y_1^{-5/2-k}\right)\nonumber\\
&\underset{y_1\to +\infty}{=}&\sum_{0\leq m\leq\frac{\beta N}{2-3\beta}}
N_{0,m}\eps^{1/3}\frac{d^l}{dy_1^l}\left(y_1^{1/2-3m}\right)+\eps^{1/3}o(y_1^{1/2-\frac{3\beta N}{2-3\beta}-l})\nonumber\\
&&+\sum_{-N\leq k\leq\frac{3\beta
    N-2}{2-3\beta}-3}\ \sum_{\substack{(n,m)\in\{1,\cdots,N\}\times\N\\ 3m-n=k}}N_{n,m}\eps^{1/3+2n/3}\frac{d^l}{dy_1^l}\left(y_1^{-5/2-k}\right)+\eps
o(y_1^{-5/2-\frac{3\beta
  N-2}{2-3\beta}+3-l}).\nonumber
\ee
Thus, for $x\in D_0\cap D_1$, we have
\be
\lefteqn{\frac{d^l}{dy_1^l}\left(\eps^{1/3}\nu(y_1)\right)
}\nonumber\\
&=&\eps^{2l/3}\sum_{0\leq
  m\leq\frac{\beta}{2-3\beta}N}\eps^{2m}N_{0,m}\frac{d^l}{dz^l}\left(z^{1/2-3m}\right)+\eps^{2l/3}\sum_{\substack{(n,m)\in\{1,\cdots,N\}\times\N\\ 3m-n\leq \frac{3\beta
    N-2}{2-3\beta}-3}}N_{n,m}\eps^{2+2m}\frac{d^l}{dz^l}\left(z^{-5/2+n-3m}\right)\nonumber\\
&&+o_{L^\infty(D_0\cap
  D_1)}(\eps^{\beta(N+1/2)+(2/3-\beta)l})\nonumber\\
&=&\eps^{2l/3}\sum_{0\leq
  m\leq\frac{\beta}{2-3\beta}N}\eps^{2m}N_{0,m}\frac{d^l}{dz^l}\left(z^{1/2-3m}\right)+\eps^{2l/3}\sum_{\substack{(n,m)\in\{1,\cdots,N\}\times\N^*\\ 3m-n\leq \frac{3\beta
    N-2}{2-3\beta}}}N_{n,m-1}\eps^{2m}\frac{d^l}{dz^l}\left(z^{1/2+n-3m}\right)\nonumber\\
&&+o_{L^\infty(D_0\cap
  D_1)}(\eps^{\beta(N+1/2)+(2/3-\beta)l})\nonumber\\
&=&\eps^{2l/3}\sum_{0\leq
  m\leq\frac{\beta}{2-3\beta}N}\eps^{2m}N_{0,m}\frac{d^l}{dz^l}\left(z^{1/2-3m}\right)+\eps^{2l/3}\sum_{\substack{(n,m)\in{\N^*}^2\\
(2-3\beta)m+\beta n\leq \beta
    N}}N_{n,m-1}\eps^{2m}\frac{d^l}{dz^l}\left(z^{1/2+n-3m}\right)\nonumber\\
&&+o_{L^\infty(D_0\cap
  D_1)}(\eps^{\beta(N+1/2)+(2/3-\beta)l}),
\ee
where in the last equality, we have neglected all the terms in the sum
over $(n,m)$ which can be incorporated in the rest
term, and we have used that the condition 
\be\label{cond1}
(n,m)\in{\N^*}^2,\quad (2-3\beta)m+\beta n\leq \beta
    N
\ee
clearly implies $n\leq N$ (even $n<N$, in fact), as well as 
\be\label{cond2}
3m-n\leq \frac{3\beta
    N-2}{2-3\beta}.
\ee
Indeed, (\ref{cond1}) can be rewritten as
\be\label{rewrite-cond}
\frac{1}{3}(3m-n)+n\left(\frac{1}{3}+\frac{\beta}{2-3\beta}\right)\leq\frac{\beta
  N}{2-3\beta},
\ee
which yields (\ref{cond2}) if we take into account that $n\geq 1$.
The result follows from the change of variable $z=\eps^{2/3}y_1$, with
$$n_{m,n}=\left\{\begin{array}{lll}N_{0,m} & {\rm if} & n=0 \\
0 & {\rm if} & n\geq 1 {\rm\ and\ }
m=0\\
N_{n,m-1} &
{\rm if} & n\geq 1 {\rm\ and\ } m\geq 1.\\
 \end{array}\right.$$
Similarly as for $\eps^{1/3}\nu$, we have
\be
\lefteqn{\frac{d^l}{dy_1^l}\left(\eps^{2/3}\lambda(y_1)\right)}\nonumber\\
&\underset{y_1\to +\infty}{\approx}&\lambda_{-1}\mathbf{1}_{\{l=0\}}+\eps^{2/3}\sum_{m=0}^\infty
L_{0,m}\frac{d^l}{dy_1^l}\left(y_1^{1-3m}\right)+\eps^{2/3}\sum_{n=1}^N\eps^{2n/3}\sum_{m=0}^\infty
L_{n,m}\frac{d^l}{dy_1^l}\left(y_1^{n-2-3m}\right)\nonumber\\
&\underset{y_1\to +\infty}{\approx}&\lambda_{-1}\mathbf{1}_{\{l=0\}}+\eps^{2/3}\sum_{m=0}^\infty
L_{0,m}\frac{d^l}{dy_1^l}\left(y_1^{1-3m}\right)+\eps^{2/3}\sum_{n=1}^N\eps^{2n/3}\sum_{m=1}^\infty
L_{n,m-1}\frac{d^l}{dy_1^l}\left(y_1^{1+n-3m}\right)\nonumber\\
&\underset{y_1\to +\infty}{\approx}&\lambda_{-1}\mathbf{1}_{\{l=0\}}+\sum_{n=0}^N\eps^{2(n+1)/3}\sum_{m=0}^\infty
\check{L}_{n,m}\frac{d^l}{dy_1^l}\left(y_1^{1+n-3m}\right),
\ee
with
$$\check{L}_{n,m}=\left\{\begin{array}{lll}L_{0,m} & {\rm if} & n=0 \\
0 & {\rm if} & n\geq 1 {\rm\ and\ }
m=0\\
L_{n,m-1} &
{\rm if} & n\geq 1 {\rm\ and\ } m\geq 1.\\
 \end{array}\right.$$
Thus,
\be\label{ddl-lambda}
\lefteqn{\frac{d^l}{dy_1^l}\left(\eps^{2/3}\lambda(y_1)\right)}\nonumber\\
&\underset{y_1\to +\infty}{\approx}&\lambda_{-1}\mathbf{1}_{\{l=0\}}+\sum_{ k=-N}^\infty\ \sum_{\substack{(n,m)\in\{0,\cdots,N\}\times\N\\ 3m-n=k}}\check{L}_{n,m}\eps^{2(n+1)/3}\frac{d^l}{dy_1^l}\left(y_1^{1-k}\right)\nonumber\\
&\underset{y_1\to +\infty}{=}&\lambda_{-1}\mathbf{1}_{\{l=0\}}+\sum_{-N\leq k\leq\frac{3\beta
    N}{2-3\beta}}\ \sum_{\substack{(n,m)\in\{0,\cdots,N\}\times\N\\ 3m-n=k}}\check{L}_{n,m}\eps^{2(n+1)/3}\frac{d^l}{dy_1^l}\left(y_1^{1-k}\right)+\eps^{2/3}o(y_1^{1-\frac{3\beta
    N}{2-3\beta}-l})\nonumber\\
&\underset{y_1\to
  +\infty}{=}&\lambda_{-1}\mathbf{1}_{\{l=0\}}+\sum_{\substack{(n,m)\in\{0,\cdots,N\}\times\N\\ 3m-n\leq\frac{3\beta
    N}{2-3\beta}
}}\check{L}_{n,m}\eps^{2(n+1)/3}\frac{d^l}{dy_1^l}\left(y_1^{1-3m+n}\right)+\eps^{2/3}o(y_1^{1-\frac{3\beta N}{2-3\beta}-l}).
\ee
Therefore for $x\in D_0\cap D_1$,
\be\label{724}
\lefteqn{\frac{d^l}{dy_1^l}\left(\eps^{2/3}\lambda(y_1)\right)}\nonumber\\
&\underset{y_1\to
  +\infty}{=}&\lambda_{-1}\mathbf{1}_{\{l=0\}}+\sum_{\substack{(n,m)\in\{0,\cdots,N\}\times\N\\ 3m-n\leq\frac{3\beta
    N}{2-3\beta}
}}\check{L}_{n,m}\eps^{2m}\eps^{2l/3}\frac{d^l}{dz^l}\left(z^{1-3m+n}\right)+\eps^{2l/3}o_{L^\infty(D_0\cap
  D_1)}(\eps^{\beta(N+1-l)}).\ \ \ \ \nonumber\\
&\underset{y_1\to
  +\infty}{=}&\lambda_{-1}\mathbf{1}_{\{l=0\}}+\eps^{2l/3}\sum_{\substack{(n,m)\in\N^2\\ (2-3\beta)m+\beta
    n\leq\beta
    N}}\check{L}_{n,m}\eps^{2m}\frac{d^l}{dz^l} \left(z^{1-3m+n}\right)+\eps^{2l/3}o_{L^\infty(D_0\cap D_1)}(\eps^{\beta(N+1-l)}),
\ee
thanks to the same remark as in (\ref{rewrite-cond}). In particular,
for $l=0$, we get
\be\label{dev-pow-dem}
\lefteqn{\eps^{1/3}\lambda(y_1)^{1/2}}\\
&=&\left(\lambda_{-1}+\sum_{\substack{(n,m)\in\N^2\\ (2-3\beta)m+\beta
    n\leq\beta
    N}}\check{L}_{n,m}\eps^{2m}z^{1-3m+n}+o_{L^\infty(D_0\cap
  D_1)}(\eps^{\beta(N+1)})\right)^{1/2}\nonumber\\
&=&\lambda_{-1}^{1/2}+\sum_{k=1}^{N+1}c_k\left(\sum_{\substack{(n,m)\in\N^2\\ (2-3\beta)m+\beta
    n\leq\beta
    N}}\check{L}_{n,m}\eps^{2m}z^{1-3m+n}\right)^k+o_{L^\infty(D_0\cap
  D_1)}(\eps^{\beta(N+1)})\nonumber\\
&=&\lambda_{-1}^{1/2}+\sum_{k=1}^{N+1}c_kz^k\!\!\!\!\!\!\!\!\!\!\!\!\!\!\!\!\!\!\sum_{\substack{\left((n_1,m_1),\cdots,(n_k,m_k)\right)\in(\N^2)^k\\ \forall
    j\in\{1,\cdots,k\},\ (2-3\beta)m_j+\beta
    n_j\leq\beta
    N}}\prod_{j=1}^k\check{L}_{n_j,m_j}(\eps^2z^{-3})^{m_1+\cdots+m_k}z^{n_1+\cdots+n_k}+o_{L^\infty(D_0\cap
  D_1)}(\eps^{\beta(N+1)})\nonumber\\
&=&\lambda_{-1}^{1/2}+\sum_{k=1}^{N+1}c_kz^k\!\!\!\!\!\!\!\!\!\sum_{\substack{(n,m)\in\N^2\\(2-3\beta)m+\beta
    n\leq\beta(N+1-k)}}\!\!\!\!\!\!\!\!\!\eps^{2m}z^{n-3m}\!\!\!\!\!\!\!\!\!\!\!\!\!\!\!\!\!\!\sum_{\substack{\left((n_1,m_1),\cdots,(n_k,m_k)\right)\in(\N^2)^k\\ n_1+\cdots+n_k=n\\m_1+\cdots+m_k=m}}\prod_{j=1}^k\check{L}_{n_j,m_j}+o_{L^\infty(D_0\cap
  D_1)}(\eps^{\beta(N+1)})\nonumber\\
&=&\lambda_{-1}^{1/2}+\sum_{\substack{(n,m)\in\N^2\\(2-3\beta)m+\beta
    n\leq\beta N}}\eps^{2m}z^{n-3m+1}\underbrace{\sum_{k=1}^{n+1}c_k\sum_{\substack{\left((n_1,m_1),\cdots,(n_k,m_k)\right)\in(\N^2)^k\\ n_1+\cdots+n_k=n-k+1\\m_1+\cdots+m_k=m}}\prod_{j=1}^k\check{L}_{n_j,m_j}}_{=:l_{m,n}}+o_{L^\infty(D_0\cap
  D_1)}(\eps^{\beta(N+1)}),\nonumber
\ee
where the $c_k$'s are some real coefficients. So, we have proved
(\ref{d2}) for $l=0$. In order to prove 
(\ref{d2}) for $l=1$, we first write
$$\frac{d}{d
  y_1}\left(\eps^{1/3}\lambda(y_1)^{1/2}\right)=\frac{1}{2}\frac{d}{d
  y_1}\left(\eps^{2/3}\lambda(y_1)\right)\left(\eps^{2/3}\lambda(y_1)\right)^{-1/2}.$$
Then, note that $\left(\eps^{2/3}\lambda(y_1)\right)^{-1/2}$ has the
same kind of asymptotic expansion as the one that appears in the right
hand side of (\ref{dev-pow-dem}). Indeed, the same calculation can be done
with the power $1/2$ replaced by $-1/2$, which only changes the values of
the $c_k$'s. Thus, for some coefficients
$(\alpha_{m,n})_{m,n\in\N^2}$, we have
\be\label{dev-pow--dem}
\left(\eps^{2/3}\lambda(y_1)\right)^{-1/2}&=&\lambda_{-1}^{-1/2}+\sum_{\substack{(n,m)\in\N^2\\(2-3\beta)m+\beta
    n\leq\beta N}}\eps^{2m}z^{n-3m+1}\alpha_{m,n}+o_{L^\infty(D_0\cap
  D_1)}(\eps^{\beta(N+1)})
\ee
From the product of this expansion with
(\ref{724}) for $l=1$, we infer that 
\be\label{727}
\frac{d}{dy_1}\left(\eps^{1/3}\lambda(y_1)^{1/2}\right)&=&\frac{1}{2}\left(\eps^{2/3}\sum_{\substack{(n,m)\in\N^2\\ (2-3\beta)m+\beta
    n\leq\beta
    N}}\check{L}_{n,m}\eps^{2m}\frac{d}{dz}
\left(z^{1-3m+n}\right)+\eps^{2/3}o_{L^\infty(D_0\cap
  D_1)}(\eps^{\beta N})\right)\nonumber\\
&&\times\left(\lambda_{-1}^{-1/2}+\sum_{\substack{(n,m)\in\N^2\\(2-3\beta)m+\beta
    n\leq\beta N}}\eps^{2m}z^{n-3m+1}\alpha_{m,n}+o_{L^\infty(D_0\cap
  D_1)}(\eps^{\beta(N+1)})\right)\nonumber\\
&=&\frac{\lambda_{-1}^{-1/2}\eps^{2/3}}{2}\sum_{\substack{(n,m)\in\N^2\\ (2-3\beta)m+\beta
    n\leq\beta
    N}}\check{L}_{n,m}\eps^{2m}
(1-3m+n)z^{n-3m}\nonumber\\
&&+\frac{\eps^{2/3}}{2}\sum_{\substack{(n,m)\in\N^*\times\N\\(2-3\beta)m+\beta
    n\leq\beta
    N}}\eps^{2m}z^{n-3m}\sum_{\substack{n_1,n_2,m_1,m_2\in\N\\ n_1+n_2=n-1\\m_1+m_2=m}}\check{L}_{n_1,m_1}(1-3m_1+n_1)\alpha_{m_2,n_2}\nonumber\\
&&+o_{L^\infty(D_0\cap
  D_1)}(\eps^{\beta N+2/3})\nonumber\\
&&=\eps^{2/3}\sum_{\substack{(m,n)\in\N^2\\ (2-3\beta)m+\beta
    n\leq \beta N}}l_{m,n}'\eps^{2m}z^{n-3m}+o_{L^\infty(D_0\cap
  D_1)}(\eps^{\beta N+2/3}),
\ee
for some coefficients $l_{m,n}'\in\R$. In order to prove (\ref{d2}), it is now
sufficient to establish that for every $m,n\geq 0$, the $l_{m,n}'$'s
and the $l_{m,n}$'s, defined respectively in (\ref{727}) and (\ref{dev-pow-dem}), are related by 
$l_{m,n}'=(1+n-3m)l_{m,n}$. For this purpose, we note, for
$z\in[\eps^\beta,2\eps^\beta]$
$\theta(z)=\eps^{1/3}\lambda(y_1)^{1/2}$, such that according to
(\ref{dev-pow-dem}) and (\ref{727}),
\be\label{devthet}
\theta(z)=\lambda_{-1}^{1/2}+\sum_{\substack{(m,n)\in\N^2\\ (2-3\beta)m+\beta
    n\leq \beta N}}l_{m,n}\eps^{2m}z^{n-3m+1}+o_{L^\infty(D_0\cap
  D_1)}(\eps^{\beta(N+1)})
\ee
and
\be\label{devthet'}
\theta'(z)=\sum_{\substack{(m,n)\in\N^2\\ (2-3\beta)m+\beta
    n\leq \beta N}}l_{m,n}'\eps^{2m}z^{n-3m}+o_{L^\infty(D_0\cap
  D_1)}(\eps^{\beta N}).
\ee
Then, we have on the one side from (\ref{devthet})
\be\label{comparthet}
\theta(2\eps^\beta)-\theta(\eps^\beta)&=&\sum_{\substack{(m,n)\in\N^2\\ (2-3\beta)m+\beta
    n\leq \beta N}}l_{m,n}(2^{n-3m+1}-1)\eps^{(2-3\beta)m+\beta(n+1)}+o_{L^\infty(D_0\cap
  D_1)}(\eps^{\beta(N+1)}),
\ee
whereas on the other side, thanks to (\ref{devthet'}),
\be\label{comparthet'}
\theta(2\eps^\beta)-\theta(\eps^\beta)&=&
\sum_{\substack{(m,n)\in\N^2\\ (2-3\beta)m+\beta
    n\leq \beta N}}l_{m,n}'\eps^{2m}\int_{\eps^\beta}^{2\eps^\beta}z^{n-3m}dz+o_{L^\infty(D_0\cap
  D_1)}(\eps^{\beta (N+1)})\nonumber\\
&=&
\sum_{\substack{(m,n)\in\N^2\\ (2-3\beta)m+\beta
    n\leq \beta N\\n-3m\neq
    -1}}\frac{l_{m,n}'}{n-3m+1}(2^{n-3m+1}-1)\eps^{(2-3\beta)m+\beta(n+1)}\nonumber\\
&&+\ln(2)\sum_{\substack{(m,n)\in\N^2\\ (2-3\beta)m+\beta
    n\leq \beta N\\n-3m= -1}}l_{m,n}'\eps^{2m}+o_{L^\infty(D_0\cap
  D_1)}(\eps^{\beta (N+1)}).
\ee
Since $\beta$ is not rational, the family of functions of the variable
$\eps$, $(\eps^{2m+\beta(n+1-3m)})_{(m,n)\in\N^2}$ is linearly
independent, and we deduce by comparison of (\ref{comparthet}) and
(\ref{comparthet'}) that $l_{m,n}'=(n-3m+1)l_{m,n}$, in both cases
$n-3m+1\neq 0$ and $n-3m+1= 0$. (\ref{d2}) for $l=1$ follows. The
proof for $l=2$ is similar.
\end{Proof}

\begin{lem}\label{lem-compar}
Let $N\geq 2$ be an integer, $\eps_0>0$ and
$\beta\in (0,2/3)\string\ \mathbb{Q}$. Let
$(\theta_1,\theta_2)_{0<\eps\leq\eps_0}$ be a sequence of pairs of
regular functions defined for $z\in [\eps^\beta,2\eps^\beta]$, such
that 
\be\label{eq-theta1}
\left\|\eps^2\Delta\theta_1+\frac{\alpha_0}{\alpha_2}(R_2^2-R_1^2)\theta_1+z\theta_1-2\alpha_1\theta_1^3-2\alpha_0\theta_2^2\theta_1\right\|_{L^\infty(D_0\cap D_1)}=o\left(\eps^{\beta(N+1/2)}\right)
\ee
and
\be\label{eq-theta2}
\left\|\eps^2\Delta\theta_2+(R_2^2-R_1^2)\theta_2+z\theta_2-2\alpha_2\theta_2^3-2\alpha_0\theta_1^2\theta_2\right\|_{L^\infty(D_0\cap D_1)}=o\left(\eps^{\beta(N+1)}\right)
\ee
are satisfied, where $\Delta\theta_j$ refers to
$\sum_{k=1}^d\frac{\partial^2}{\partial
  x_k^2}\left(\theta_j(R_1^2-|x|^2)\right)=-2d\theta_j'(z)+4(R_1^2-z)\theta_j''(z)$
(with $z=R_1^2-|x|^2$). We assume that there exists two families of real numbers $p_{m,n}$, $q_{m,n}$, defined for every
  $(m,n)\in\N^2$ such that $(2-3\beta)m+\beta n\leq\beta N$, such that
\be\label{theta1}
\forall l\in\{0,1,2\},&&\left\|\theta_1^{(l)}-\sum_{\substack{(m,n)\in\N^2\\(2-3\beta)m+\beta
    n\leq\beta N}}
  \eps^{2m}p_{m,n}\frac{d^l}{dz^l}\left(z^{1/2-3m+n}\right)\right\|_{L^\infty(D_0\cap
  D_1)}\!\!\!\!\underset{\eps\to
    0}{=}o\left(\eps^{\beta(N+1/2-l)}\right)\ \ \ \ \ \ 
\ee
and
\be\label{theta2}
\forall l\in\{0,1,2\},&&\!\!\!\!\!\left\|\theta_2^{(l)}-\lambda_{-1}^{1/2}\mathbf{1}_{\{l=0\}}-\!\!\!\!\!\!\sum_{\substack{(m,n)\in\N^2\\ (2-3\beta)m+\beta
    n\leq \beta N}}\!\!\!\!\eps^{2m}q_{m,n}\ \frac{d^l}{dz^l}\left(z^{1+n-3m}\right)\right\|_{L^\infty(D_0\cap
  D_1)}\!\!\!\!\!\!\!\!\!\!\!\underset{\eps\to
  0}{=}o(\eps^{\beta(N+1-l)}).\ \ \ \ \ 
\ee
 Then,
equations (\ref{theta1}), (\ref{theta2}), (\ref{eq-theta1}) and
(\ref{eq-theta2}) entirely determine the values of the $p_{m,n}$'s and the $q_{m,n}$'s for $(2-3\beta)m+\beta
n\leq\beta (N-1)$. Moreover, these coefficients do not depend on $N$
or $\beta$.
\end{lem}

\begin{Proof}
For convenience, for every $(m,n)\in\N^2$, we denote
$p_{m,n}'=(1/2-3m+n)p_{m,n}$,
$p_{m,n}''=(-1/2-3m+n)(1/2-3m+n)p_{m,n}$, 
$q_{m,n}'=(1+n-3m)q_{m,n}$ and $q_{m,n}''=(n-3m)(1+n-3m)q_{m,n}$. For functions
$(\theta_1,\theta_2)$ that satisfy (\ref{theta1}) and (\ref{theta2}), let
us calculate the function that appears in the left hand side of
(\ref{eq-theta1}), evaluated at $z=\eps^\beta$. In the calculation
below, implicitely, $\theta_j=\theta_j(\eps^\beta)$.
\be\label{dev-eq-theta1}
\lefteqn{\eps^2\Delta\theta_1+\frac{\alpha_0}{\alpha_2}(R_2^2-R_1^2)\theta_1+z\theta_1-2\alpha_1\theta_1^3-2\alpha_0\theta_2^2\theta_1}\nonumber\\
&=&-2\sum_{\substack{m,n\geq 0\\(2-3\beta)m+\beta n\leq \beta N}}(d
p_{m,n}'+2p_{m,n}'')\eps^{(2-3\beta)(m+1)+\beta(
  n+2)+\beta/2}\nonumber\\
&&+4R_1^2\sum_{\substack{m,n\geq 0\\(2-3\beta)m+\beta
    n\leq \beta N}}p_{m,n}''\eps^{(2-3\beta)(m+1)+\beta(
  n+1)+\beta/2}\nonumber\\
&&+\frac{\alpha_0}{\alpha_2}(R_2^2-R_1^2)\sum_{\substack{m,n\geq 0\\(2-3\beta)m+\beta
    n\leq \beta N}}p_{m,n}\eps^{(2-3\beta)m+\beta
  n+\beta/2}+\sum_{\substack{m,n\geq 0\\(2-3\beta)m+\beta
    n\leq \beta N}}p_{m,n}\eps^{(2-3\beta)m+\beta
  (n+1)+\beta/2}\nonumber\\
&&-2\alpha_1\sum_{\substack{m,n\geq 0\\(2-3\beta)m+\beta
    n\leq \beta N}}\big(\sum_{\substack{m_1,m_2,m_3,n_1,n_2,n_3\geq 0\\m_1+m_2+m_3=m\\n_1+n_2+n_3=n}}p_{m_1,n_1}p_{m_2,n_2}p_{m_3,n_3}\big)\eps^{(2-3\beta)m+\beta
  n+3\beta/2}\nonumber\\
&&-2\alpha_0\lambda_{-1}\sum_{\substack{m,n\geq 0\\(2-3\beta)m+\beta
    n\leq \beta N}}p_{m,n}\eps^{(2-3\beta)m+\beta
  n+\beta/2}\nonumber\\
&&-4\alpha_0\lambda_{-1}^{1/2}\sum_{\substack{m,n\geq 0\\(2-3\beta)m+\beta
    n\leq \beta N}}\big(\sum_{\substack{m_1,m_2,n_1,n_2\geq
    0\\n_1+n_2=n\\m_1+m_2=m}}q_{m_1,n_1}p_{m_2,n_2}\big)\eps^{(2-3\beta)m+\beta
  n+3\beta/2}\nonumber\\
&&-2\alpha_0\sum_{\substack{m,n\geq 0\\(2-3\beta)m+\beta
    n\leq \beta N}}\big(\sum_{\substack{ m_1,m_2,m_3,n_1,n_2,n_3\geq
    0,\\n_1+n_2+n_3=n\\m_1+m_2+m_3=m}}q_{m_1,n_1}q_{m_2,n_2}p_{m_3,n_3}\big)\eps^{(2-3\beta)m+\beta
  n+5\beta/2}+o(\eps^{\beta N+\beta/2})\nonumber\\
&=&-2\sum_{\substack{m\geq 1,n\geq 2\\(2-3\beta)m+\beta n\leq \beta N}}(d
p_{m-1,n-2}'+2p_{m-1,n-2}'')\eps^{(2-3\beta)m+\beta n+\beta/2}\nonumber\\
&&+4R_1^2\sum_{\substack{m,n\geq 1\\(2-3\beta)m+\beta
    n\leq \beta N}}p_{m-1,n-1}''\eps^{(2-3\beta)m+\beta n+\beta/2}+\sum_{\substack{m\geq 0,n\geq 1\\(2-3\beta)m+\beta
    n\leq \beta N}}p_{m,n-1}\eps^{(2-3\beta)m+\beta n+\beta/2}\nonumber\\
&&-2\alpha_1\sum_{\substack{m\geq 0,n\geq 1\\(2-3\beta)m+\beta
    n\leq \beta N}}\big(\sum_{\substack{m_1,m_2,m_3,n_1,n_2,n_3\geq 0\\m_1+m_2+m_3=m\\n_1+n_2+n_3=n-1}}p_{m_1,n_1}p_{m_2,n_2}p_{m_3,n_3}\big)\eps^{(2-3\beta)m+\beta
  n+\beta/2}\nonumber\\
&&-4\alpha_0\lambda_{-1}^{1/2}\sum_{\substack{m\geq 0,n\geq 1\\(2-3\beta)m+\beta
    n\leq \beta N}}\big(\sum_{\substack{n_1,n_2,m_1,m_2\geq
    0,\\n_1+n_2=n-1\\m_1+m_2=m}}q_{m_1,n_1}p_{m_2,n_2}\big)\eps^{(2-3\beta)m+\beta
  n+\beta/2}\\
&&-2\alpha_0\sum_{\substack{m\geq 0, n\geq 2\\(2-3\beta)m+\beta
    n\leq \beta N}}\big(\sum_{\substack{n_1,n_2,n_3,m_1,m_2,m_3\geq
    0\\n_1+n_2+n_3=n-2\\m_1+m_2+m_3=m}}q_{m_1,n_1}q_{m_2,n_2}p_{m_3,n_3}\big)\eps^{(2-3\beta)m+\beta
  n+\beta/2}+o(\eps^{\beta N+\beta/2}),\nonumber
\ee
where we have used (\ref{I-1}). Since $\beta$ is not rational, the functions
$\left((0,\eps_0)\ni\eps\mapsto \eps^{(2-3\beta)m+\beta
  n}\right)_{m,n\in\N^2}$ are two by two distinct, and therefore
linearly independent. According to (\ref{eq-theta1}), we deduce from
(\ref{dev-eq-theta1}):\\
\begin{itemize}
\item for $m=0$ and $1\leq n\leq N$,
\begin{eqnarray*}
\lefteqn{p_{0,n-1}-2\alpha_1\sum_{\substack{n_1,n_2,n_3\geq
    0\\n_1+n_2+n_3=n-1}}p_{0,n_1}p_{0,n_2}p_{0,n_3}}\nonumber\\
&&-4\alpha_0\lambda_{-1}^{1/2}\sum_{\substack{
    n_1,n_2\geq
    0\\n_1+n_2=n-1}}q_{0,n_1}p_{0,n_2}-2\alpha_0\sum_{\substack{
    n_1,n_2,n_3\geq
    0\\n_1+n_2+n_3=n-2}}q_{0,n_1}q_{0,n_2}p_{0,n_3}\ =\ 0,
\end{eqnarray*}
which can be rewritten as
\be\label{I01}
p_{0,0}-2\alpha_1p_{0,0}^3-4\alpha_0\lambda_{-1}^{1/2}q_{0,0}p_{0,0}&=&0
\ee
for $n=1$, and
\be\label{I0n}
\lefteqn{p_{0,n-1}\left(1-6\alpha_1p_{0,0}^2-4\alpha_0\lambda_{-1}^{1/2}q_{0,0}\right)-4\alpha_0\lambda_{-1}^{1/2}p_{0,0}q_{0,n-1}}\\
&=&2\alpha_1\sum_{\substack{0\leq
    n_1,n_2,n_3<n-1\\n_1+n_2+n_3=n-1}}p_{0,n_1}p_{0,n_2}p_{0,n_3}+4\alpha_0\lambda_{-1}^{1/2}\!\!\!\sum_{\substack{
    0\leq n_1,n_2<n-1\\n_1+n_2=n-1}}\!\!\!q_{0,n_1}p_{0,n_2}+2\alpha_0\!\!\!\sum_{\substack{
    n_1,n_2,n_3\geq
    0\\n_1+n_2+n_3=n-2}}\!\!\!q_{0,n_1}q_{0,n_2}p_{0,n_3}\nonumber
\ee
for $n\geq 2$.
\item for $ 1\leq m\leq \beta(N-1)/(2-3\beta)$ and $n=1$,
\begin{eqnarray*}
4R_1^2p_{m-1,0}''+p_{m,0}-2\alpha_1\sum_{\substack{m_1,m_2,m_3\geq
    0\\m_1+m_2+m_3=m}}p_{m_1,0}p_{m_2,0}p_{m_3,0}
-4\alpha_0\lambda_{-1}^{1/2}\sum_{\substack{m_1,m_2\geq 0\\m_1+m_2=m}}q_{m_1,0}p_{m_2,0}&=&0,\nonumber
\end{eqnarray*}
which can be rewritten as
\be\label{Im1}
\lefteqn{p_{m,0}\left(1-6\alpha_1p_{0,0}^2-4\alpha_0\lambda_{-1}^{1/2}q_{0,0}\right)-4\alpha_0\lambda_{-1}^{1/2}p_{0,0}q_{m,0}}\nonumber\\
&=&-4R_1^2p_{m-1,0}''+2\alpha_1\sum_{\substack{0\leq m_1,m_2,m_3<m\\m_1+m_2+m_3=m}}p_{m_1,0}p_{m_2,0}p_{m_3,0}
+4\alpha_0\lambda_{-1}^{1/2}\sum_{\substack{0\leq m_1,m_2<m\\m_1+m_2=m}}q_{m_1,0}p_{m_2,0}.
\ee
\item for  $m\geq 1$ and $n\geq 2$ such that $(2-3\beta)m+\beta n\leq \beta N$,
\begin{eqnarray*}
-2(d
p_{m-1,n-2}'+2p_{m-1,n-2}'')+4R_1^2p_{m-1,n-1}''+p_{m,n-1}\nonumber\\
-2\alpha_1\sum_{\substack{m_1,m_2,m_3,n_1,n_2,n_3\geq 0\\m_1+m_2+m_3=m\\n_1+n_2+n_3=n-1}}p_{m_1,n_1}p_{m_2,n_2}p_{m_3,n_3}\\
-4\alpha_0\lambda_{-1}^{1/2}\sum_{\substack{m_1,m_2,n_1,n_2\geq 0\\n_1+n_2=n-1\\m_1+m_2=m}}q_{m_1,n_1}p_{m_2,n_2}
-2\alpha_0\sum_{\substack{m_1,m_2,m_3,n_1,n_2,n_3\geq 0\\n_1+n_2+n_3=n-2\\m_1+m_2+m_3=m}}q_{m_1,n_1}q_{m_2,n_2}p_{m_3,n_3}&=&0,\nonumber
\end{eqnarray*}
which can be rewritten as
\be\label{Imn}
\lefteqn{\left(1-6\alpha_1p_{0,0}^2-4\alpha_0\lambda_{-1}^{1/2}q_{0,0}\right)p_{m,n-1}-4\alpha_0\lambda_{-1}^{1/2}p_{0,0}q_{m,n-1}}\\
&=&2(d
p_{m-1,n-2}'+2p_{m-1,n-2}'')-4R_1^2p_{m-1,n-1}''
+2\alpha_1\sum_{\substack{m_1,m_2,m_3,n_1,n_2,n_3\geq
    0\\m_1+m_2+m_3=m\\n_1+n_2+n_3=n-1\\\forall j\in\{1,2,3\},(m_j,n_j)\neq(m,n-1)}}p_{m_1,n_1}p_{m_2,n_2}p_{m_3,n_3}\nonumber\\
&&+4\alpha_0\lambda_{-1}^{1/2}\sum_{\substack{m_1,m_2,n_1,n_2\geq 0,\\n_1+n_2=n-1\\m_1+m_2=m\\\forall j\in\{1,2\},(m_j,n_j)\neq(m,n-1)}}q_{m_1,n_1}p_{m_2,n_2}+2\alpha_0\sum_{\substack{m_1,m_2,m_3,n_1,n_2,n_3,m_3\geq 0,\\n_1+n_2+n_3=n-2\\m_1+m_2+m_3=m}}q_{m_1,n_1}q_{m_2,n_2}p_{m_3,n_3}.\nonumber
\ee
\end{itemize}
Next, we perform the same kind of calculations with the function that
appears in the left hand side of (\ref{eq-theta2}).
\be
\lefteqn{\eps^2\Delta\theta_2+(R_2^2-R_1^2)\theta_2+z\theta_2-2\alpha_2\theta_2^3-2\alpha_0\theta_1^2\theta_2}\nonumber\\
&=&-2\sum_{\substack{m,n\geq 0\\(2-3\beta)m+\beta n\leq \beta N}}(d
q_{m,n}'+2q_{m,n}'')\eps^{(2-3\beta)(m+1)+\beta(
  n+2)+\beta}\nonumber\\
&&+4R_1^2\sum_{\substack{m,n\geq 0\\(2-3\beta)m+\beta
    n\leq \beta N}}q_{m,n}''\eps^{(2-3\beta)(m+1)+\beta(
  n+1)+\beta}\nonumber\\
&&+(R_2^2-R_1^2)\lambda_{-1}^{1/2}+(R_2^2-R_1^2)\sum_{\substack{m,n\geq 0\\(2-3\beta)m+\beta
    n\leq \beta N}}q_{m,n}\eps^{(2-3\beta)m+\beta
  n+\beta}\nonumber\\
&&+\eps^\beta\lambda_{-1}^{1/2}+\sum_{\substack{m\geq 0,n\geq 0\\(2-3\beta)m+\beta
    n\leq \beta N}}q_{m,n}\eps^{(2-3\beta)m+\beta
  (n+1)+\beta}\nonumber\\
&&-2\alpha_2\lambda_{-1}^{3/2}-6\alpha_2\lambda_{-1}\sum_{\substack{m\geq 0,n\geq 0\\(2-3\beta)m+\beta
    n\leq \beta N}}q_{m,n}\eps^{(2-3\beta)m+\beta
  n+\beta}\nonumber\\
&&-6\alpha_2\lambda_{-1}^{1/2}\sum_{\substack{m,n\geq 0\\(2-3\beta)m+\beta
    n\leq \beta N}}\big(\sum_{\substack{m_1,m_2,n_1,n_2\geq
  0\\m_1+m_2=m\\n_1+n_2=n}}q_{m_1,n_1}q_{m_2,n_2}\big)\eps^{(2-3\beta)m+\beta (n+1)+\beta}\nonumber\\
&&-2\alpha_2\sum_{\substack{m\geq
    0,n\geq 0\\(2-3\beta)m+\beta
    n\leq \beta N}}\big(\sum_{\substack{m_1,m_2,m_3,n_1,n_2,n_3\geq
  0,\\n_1+n_2+n_3=n\\m_1+m_2+m_3=m}}q_{m_1,n_1}q_{m_2,n_2}q_{m_3,n_3}\big)\eps^{(2-3\beta)m+\beta (n+2)+\beta}\nonumber\\
&&-2\alpha_0\lambda_{-1}^{1/2}\sum_{\substack{m,n\geq 0\\(2-3\beta)m+\beta
    n\leq \beta N}}\big(\sum_{\substack{m_1,m_2,n_1,n_2\geq
  0\\m_1+m_2=m\\n_1+n_2=n}}p_{m_1,n_1}p_{m_2,n_2}\big)\eps^{(2-3\beta)m+\beta
  n+\beta}\nonumber\\
&&-2\alpha_0\sum_{\substack{m,n\geq 0\\(2-3\beta)m+\beta
    n\leq \beta N}}\big(\sum_{\substack{m_1,m_2,m_3,n_1,n_2,n_3\geq
  0\\m_1+m_2+m_3=m\\n_1+n_2+n_3=n}}p_{m_1,n_1}p_{m_2,n_2}q_{m_3,n_3}\big)\eps^{(2-3\beta)m+\beta
  (n+1)+\beta}+o(\eps^{\beta(N+1)}).\nonumber
\ee
Thus, changing the indices and throwing away all the terms that can be
incorporated in the rest,
\be\label{dev-eq-theta2}
\lefteqn{\eps^2\Delta\theta_2+(R_2^2-R_1^2)\theta_2+z\theta_2-2\alpha_2\theta_2^3-2\alpha_0\theta_1^2\theta_2}\nonumber\\
&=&-2\sum_{\substack{m\geq 1,n\geq 2\\(2-3\beta)m+\beta n\leq \beta N}}(d
q_{m-1,n-2}'+2q_{m-1,n-2}'')\eps^{(2-3\beta)m+\beta
  n+\beta}\nonumber\\
&&+4R_1^2\sum_{\substack{m,n\geq 1\\(2-3\beta)m+\beta
    n\leq \beta N}}q_{m-1,n-1}''\eps^{(2-3\beta)m+\beta
  n+\beta}\nonumber\\
&&+(R_2^2-R_1^2)\sum_{\substack{m,n\geq 0\\(2-3\beta)m+\beta
    n\leq \beta N}}q_{m,n}\eps^{(2-3\beta)m+\beta
  n+\beta}\nonumber\\
&&+\eps^\beta\lambda_{-1}^{1/2}+\sum_{\substack{m\geq 0,n\geq 1\\(2-3\beta)m+\beta
    n\leq \beta N}}q_{m,n-1}\eps^{(2-3\beta)m+\beta
  n+\beta}\nonumber\\
&&-6\alpha_2\lambda_{-1}\sum_{\substack{m\geq 0,n\geq 0\\(2-3\beta)m+\beta
    n\leq \beta N}}q_{m,n}\eps^{(2-3\beta)m+\beta
  n+\beta}\nonumber\\
&&-6\alpha_2\lambda_{-1}^{1/2}\sum_{\substack{m\geq 0,n\geq 1\\(2-3\beta)m+\beta
    n\leq \beta N}}\big(\sum_{\substack{m_1,m_2,n_1,n_2\geq
  0\\m_1+m_2=m\\n_1+n_2=n-1}}q_{m_1,n_1}q_{m_2,n_2}\big)\eps^{(2-3\beta)m+\beta n+\beta}\nonumber\\
&&-2\alpha_2\sum_{\substack{m\geq
    0,n\geq 2\\(2-3\beta)m+\beta
    n\leq \beta N}}\big(\sum_{\substack{m_1,m_2,m_3,n_1,n_2,n_3\geq
  0,\\n_1+n_2+n_3=n-2\\m_1+m_2+m_3=m}}q_{m_1,n_1}q_{m_2,n_2}q_{m_3,n_3}\big)\eps^{(2-3\beta)m+\beta n+\beta}\nonumber\\
&&-2\alpha_0\lambda_{-1}^{1/2}\sum_{\substack{m,n\geq 0\\(2-3\beta)m+\beta
    n\leq \beta N}}\big(\sum_{\substack{m_1,m_2,n_1,n_2\geq
  0\\m_1+m_2=m\\n_1+n_2=n}}p_{m_1,n_1}p_{m_2,n_2}\big)\eps^{(2-3\beta)m+\beta
  n+\beta}\\
&&-2\alpha_0\sum_{\substack{m\geq 0,n\geq 1\\(2-3\beta)m+\beta
    n\leq \beta N}}\big(\sum_{\substack{m_1,m_2,m_3,n_1,n_2,n_3\geq
  0\\m_1+m_2+m_3=m\\n_1+n_2+n_3=n-1}}p_{m_1,n_1}p_{m_2,n_2}q_{m_3,n_3}\big)\eps^{(2-3\beta)m+\beta
  n+\beta}+o(\eps^{\beta(N+1)}).\nonumber
\ee
According to (\ref{eq-theta2}), the right hand side of
(\ref{dev-eq-theta2}) is equal to 0, up to the rest term
$o(\eps^{\beta(N+1)})$. Thus, the linear independance of the family of
functions of $\eps$, $\left(\eps^{(2-3\beta)m+\beta n}\right)_{m,n\geq
  0}$ yields:
\begin{itemize}
\item for $m=0$, $n=0$, thanks to (\ref{I-1}), we get
\be\label{II00}
-2(R_2^2-R_1^2)q_{0,0}=\lambda_{-1}^{1/2}(2\alpha_0p_{0,0}^2-1),
\ee
\item for $m=0$ and $1\leq n\leq N$,
\begin{eqnarray*}
\lefteqn{(R_2^2-R_1^2)q_{0,n}+q_{0,n-1}-6\alpha_2\lambda_{-1}q_{0,n}-6\alpha_2\lambda_{-1}^{1/2}\sum_{\substack{n_1,n_2\geq
      0,\\n_1+n_2=n-1}}q_{0,n_1}q_{0,n_2}}\nonumber\\
&&-2\alpha_2\sum_{\substack{n_1,n_2,n_3\geq
      0,\\n_1+n_2+n_3=n-2}}q_{0,n_1}q_{0,n_2}q_{0,n_3}-2\alpha_0\lambda_{-1}^{1/2}\sum_{\substack{n_1,n_2\geq
  0\\n_1+n_2=n}}p_{0,n_1}p_{0,n_2}-2\alpha_0\sum_{\substack{n_1,n_2,n_3\geq
  0\\n_1+n_2+n_3=n-1}}p_{0,n_1}p_{0,n_2}q_{0,n_3}=0,
\end{eqnarray*}
which, using (\ref{I-1}), can be rewritten as
\be\label{II0n}
\lefteqn{-2(R_2^2-R_1^2)q_{0,n}-4\alpha_0p_{0,0}\lambda_{-1}^{1/2}p_{0,n}}\\
&=&-q_{0,n-1}+6\alpha_2\lambda_{-1}^{1/2}\sum_{\substack{n_1,n_2\geq
    0,\\n_1+n_2=n-1}}q_{0,n_1}q_{0,n_2}+2\alpha_2\sum_{\substack{n_1,n_2,n_3\geq
    0,\\n_1+n_2+n_3=n-2}}q_{0,n_1}q_{0,n_2}q_{0,n_3}\nonumber\\&&+2\alpha_0\lambda_{-1}^{1/2}\sum_{\substack{0\leq
    n_1,n_2<n\\n_1+n_2=n}}p_{0,n_1}p_{0,n_2}+2\alpha_0\sum_{\substack{n_1,n_2,n_3\geq
  0\\n_1+n_2+n_3=n-1}}p_{0,n_1}p_{0,n_2}q_{0,n_3},\nonumber
\ee
\item for $1\leq m\leq \beta N/(2-3\beta)$ and $n=0$,
\begin{eqnarray*}
(R_2^2-R_1^2)q_{m,0}
-6\alpha_2\lambda_{-1}q_{m,0}
-2\alpha_0\lambda_{-1}^{1/2}\sum_{\substack{m_1,m_2\geq 0,\\m_1+m_2=m}}p_{m_1,0}p_{m_2,0}=0,
\end{eqnarray*}
that is
\be\label{IIm0}
-2(R_2^2-R_1^2)q_{m,0}-4\alpha_0p_{0,0}\lambda_{-1}^{1/2}p_{m,0}
&=&2\alpha_0\lambda_{-1}^{1/2}\sum_{\substack{0\leq m_1,m_2<m,\\m_1+m_2=m}}p_{m_1,0}p_{m_2,0}.
\ee
\item for $m\geq 1$ and $n\geq 1$ such that $(2-3\beta)m+\beta n\leq \beta N$,
\begin{eqnarray*}
-2(dq_{m-1,n-2}'+2q_{m-1,n-2}'')\mathbf{1}_{\{n\geq 2\}}+4R_1^2q_{m-1,n-1}''+(R_2^2-R_1^2)q_{m,n}+q_{m,n-1}\nonumber\\
-6\alpha_2\lambda_{-1}q_{m,n}-6\alpha_2\lambda_{-1}^{1/2}\sum_{\substack{m_1,m_2,n_1,n_2\geq
  0\\m_1+m_2=m,\\n_1+n_2=n-1}}q_{m_1,n_1}q_{m_2,n_2}-2\alpha_2\sum_{\substack{m_1,m_2,m_3,n_1,n_2,n_3\geq
  0,\\n_1+n_2+n_3=n-2\\m_1+m_2+m_3=m}}q_{m_1,n_1}q_{m_2,n_2}q_{m_3,n_3}\nonumber\\
-2\alpha_0\lambda_{-1}^{1/2}\sum_{\substack{m_1,m_2,n_1,n_2\geq
  0\\m_1+m_2=m\\n_1+n_2=n}}p_{m_1,n_1}p_{m_2,n_2}
-2\alpha_0\sum_{\substack{m_1,m_2,m_3,n_1,n_2,n_3\geq
  0\\m_1+m_2+m_3=m\\n_1+n_2+n_3=n-1}}p_{m_1,n_1}p_{m_2,n_2}q_{m_3,n_3}=0\nonumber
\end{eqnarray*}
which can be rewritten as
\be\label{IImn}
\lefteqn{-2(R_2^2-R_1^2)q_{m,n}-4\alpha_0\lambda_{-1}^{1/2}p_{0,0}p_{m,n}}\nonumber\\
&=&2(dq_{m-1,n-2}'+2q_{m-1,n-2}'')\mathbf{1}_{\{n\geq 2\}}-4R_1^2q_{m-1,n-1}''-q_{m,n-1}
+6\alpha_2\lambda_{-1}^{1/2}\sum_{\substack{0\leq m_1,m_2,n_1,n_2\\n_1+n_2=n-1\\m_1+m_2=m}}q_{m_1,n_1}q_{m_2,n_2}\nonumber\\
&&+2\alpha_2\sum_{\substack{m_1,m_2,m_3,n_1,n_2,n_3\geq
  0\\n_1+n_2+n_3=n-2\\m_1+m_2+m_3=m}}q_{m_1,n_1}q_{m_2,n_2}q_{m_3,n_3}\\
&&+2\alpha_0\lambda_{-1}^{1/2}\sum_{\substack{m_1,m_2,n_1,n_2\geq
  0\\m_1+m_2=m\\n_1+n_2=n\\\forall j\in\{1,2\},(m_j,n_j)\neq (m,n)}}p_{m_1,n_1}p_{m_2,n_2}
+2\alpha_0\sum_{\substack{m_1,m_2,m_3,n_1,n_2,n_3\geq
  0\\m_1+m_2+m_3=m\\n_1+n_2+n_3=n-1}}p_{m_1,n_1}p_{m_2,n_2}q_{m_3,n_3}\nonumber
\ee
\end{itemize}
Next, we show that the system of equations satisfied by the
$p_{m,n}$'s and $q_{m,n}$'s has a unique solution such that
$p_{0,0}>0$. First, plugging
\be\label{preq01}
q_{0,0}&=&\frac{1-2\alpha_0p_{0,0}^2}{2(R_2^2-R_1^2)}\lambda_{-1}^{1/2},
\ee
(which comes from (\ref{II00})) into (\ref{I01}) and using
also (\ref{I-1}), we get 
\be\label{p00}
p_{0,0}&=&\left(\frac{\Gamma_2}{2\alpha_1\Gamma_{12}}\right)^{1/2},
\ee
and
\be\label{q01}
q_{0,0}&=& \frac{\Gamma_1}{2\Gamma_{12}(2\alpha_2(R_2^2-R_1^2))^{1/2}}.
\ee
Next, for $1\leq n\leq N-1$, the $q_{0,n}$'s and the $p_{0,n}$'s are constructed
recursively thanks to (\ref{II0n}) as well as (\ref{I0n}) with $n$ replaced
by $n+1$. We solve the system obtained by combination of these two equations by inverting the matrix
\begin{eqnarray*}
M=\left[\begin{array}{cc}-2(R_2^2-R_1^2)&-4\alpha_0p_{0,0}\lambda_{-1}^{1/2}\\-4\alpha_0p_{0,0}\lambda_{-1}^{1/2}&1-6\alpha_1p_{0,0}^2-4\alpha_0\lambda_{-1}^{1/2}q_{0,0}\end{array}\right]
&=&\left[\begin{array}{cc}-2(R_2^2-R_1^2)&-2\alpha_0\left(\frac{(R_2^2-R_1^2)\Gamma_2}{\alpha_1\alpha_2\Gamma_{12}}\right)^{1/2}\\-2\alpha_0\left(\frac{(R_2^2-R_1^2)\Gamma_2}{\alpha_1\alpha_2\Gamma_{12}}\right)^{1/2}&-2\frac{\Gamma_2}{\Gamma_{12}}\end{array}\right],
\end{eqnarray*}
where we have used (\ref{I-1}), (\ref{p00}) and (\ref{q01}). The
determinant of $M$ is 
$$\det M=4(R_2^2-R_1^2)\Gamma_2>0,$$ 
therefore $M$ is invertible, and there is a unique possible choice for
$(q_{0,n},p_{0,n})$ for $1\leq n\leq N-1$ such that the assumptions of
the lemma are satisfied. Then, for $1\leq m\leq
\beta(N-1)/(2-3\beta)$, the $q_{m,0}$'s and the $p_{m,0}$'s are constructed
recursively thanks to (\ref{Im1})  and (\ref{IIm0}) by inverting the
same matrix $M$. Finally, if $m\geq 1$, $n\geq 1$ and
$(2-3\beta)m+\beta n\leq \beta(N-1)$ and if the $q_{k,l}$'s and the
$p_{k,l}$'s are known for every $k\leq m$, $l\leq n$ and $(k,l)\neq
(m,n)$, $(q_{m,n},p_{m,n})$ is entirely determined because the
system made of (\ref{Imn}) for $n$ replaced by $n+1$ and (\ref{IImn})
has a unique solution thanks to the invertibility of $M$. This way, we prove recursively that the
assumptions of the lemma determine completely the values of the
coefficients $q_{m,n}$ and $p_{m,n}$, provided $(2-3\beta)m+\beta n\leq \beta(N-1)$.
\end{Proof}

\begin{lem}\label{lem-compar-om-nu-tau-l}
Let $N\geq 3$ be an integer, $M \geq \frac{\beta}{2-3\beta}N$, and $\omega$, $\tau$, $\nu$, $\lambda$ given by (\ref{troncat}). 
Then for $l=0,1,2$, we have
\be\label{comp-om-nu}
\left\|\frac{d^l}{dz^l}\left(\omega(z)-\eps^{1/3}\nu(y_1)\right)\right\|_{L^\infty(D_0\cap
  D_1)}=o(\eps^{\beta(N-1/2-l)})
\ee
and
\be\label{comp-tau-lam}
\left\|\frac{d^l}{dz^l}\left(\tau(z)-\eps^{1/3}\lambda(y_1)^{1/2}\right)\right\|_{L^\infty(D_0\cap
  D_1)}=o(\eps^{\beta (N-l)}).
\ee
\end{lem}

\begin{Proof}
The assumptions (\ref{theta1}) and (\ref{theta2}) made on
$(\theta_1(z),\theta_2(z))$ in Lemma \ref{lem-compar} are satisfied by
$(\omega(z),\tau(z))$ thanks to Lemma \ref{lem-omega-tau}, and also by
$(\eps^{1/3}\nu(y_1),\eps^{1/3}\lambda(y_1)^{1/2})$ thanks to Lemma
\ref{lem-lam-nu}. Assumptions (\ref{eq-theta1}) and (\ref{eq-theta2})
are satisfied by $(\omega(z),\tau(z))$ thanks to
(\ref{est-f01-D0}) and (\ref{est-f02-D0}), and they are also satisfied
by $(\eps^{1/3}\nu(y_1),\eps^{1/3}\lambda(y_1)^{1/2})$ thanks to
Lemma \ref{lem-tronc-nu-l}. Therefore, Lemma
\ref{lem-compar} ensures that for every $(m,n)\in\N^2$ such that
$(2-3\beta)m+\beta n\leq \beta(N-1)$, $w_{m,n}=n_{m,n}$ and
$t_{m,n}=l_{m,n}$. In particular, (\ref{comp-om-nu}) and
(\ref{comp-tau-lam}) are satisfied.
\end{Proof}

\subsection{Comparison of $\eps^{1/3}(\nu,\lambda^{1/2})$ and
  $(0,\eps^{1/3}\mu)$ in $D_1\cap D_2$}\label{secmuasymp}
We first give an expansion  of  $\eps^{1/3}(\nu,\lambda^{1/2})$ into powers of $\eps$ in $D_1\cap D_2$, as $\eps\to 0$.

\begin{lem}\label{-lem-lam-nu}
Let $N\geq 1$ be an integer. There exist a family of numbers $(\widetilde{l}_{m,n})_{m\geq 0,n\geq 0}$ which does not depend on $N$ such
that if $\nu$ and
$\lambda$ are given by (\ref{troncat}), then for every $\alpha>0$,
\be\label{-d1}
\left\|\frac{d^l}{dz^l}\left(\eps^{1/3}\nu(y_1)\right)\right\|_{L^\infty(D_1\cap
  D_2)}\underset{\eps\to 0}{=}o\left(\eps^{\alpha}\right)
\ee
and
\be\label{-d2}
\left\|\frac{d^l}{dz^l}\left(\eps^{1/3}\lambda(y_1)^{1/2}\right)-\lambda_{-1}^{1/2}\mathbf{1}_{l=0}-\!\!\!\!\!\!\!\!\!\!\!\sum_{\substack{(m,n)\in\N^2\\ (2-3\beta)m+\beta
    n\leq \beta N\\1+n-3m\geq 0}}\!\!\!\!\eps^{2m}\widetilde{l}_{m,n}\ \frac{d^l}{dz^l}\left(z^{1+n-3m}\right)\right\|_{L^\infty(D_1\cap
  D_2)}\!\!\!\!\!\!\!\!\!\!\!&\underset{\eps\to
  0}{=}&\!\!\!o(\eps^{\beta(N+1-l)}).\ \ \ \ 
\ee
\end{lem}

\begin{Proof}
For $x\in D_1\cap D_2$, we have $-2\eps^{\beta-2/3}\leq y_1\leq-\eps^{\beta-2/3}\to -\infty$
as $\eps\to 0$. Thus, (\ref{-d1}) follows from (\ref{asymp-nun}) and
(\ref{asymp-nu0}). As for (\ref{-d2}), we proceed like in the proof of
Lemma \ref{lem-lam-nu}. First, from (\ref{asymp-lamn}) and
(\ref{asymp-lam0}), we have
\be
\frac{d^l}{dy_1^l}\left(\eps^{2/3}\lambda(y_1)\right)
&\underset{y_1\to
  -\infty}{\approx}&\lambda_{-1}\mathbf{1}_{\{l=0\}}+\eps^{2/3}\frac{d^l}{dy_1^l}\left(y_1/(2\alpha_2)\right)+\eps^{2/3}\sum_{n=1}^N\eps^{2n/3}\sum_{0\leq
  m\leq(n-2)/3}
\widetilde{L}_{n,m}\frac{d^l}{dy_1^l}\left(y_1^{n-2-3m}\right)\nonumber\\
&\underset{y_1\to
  -\infty}{\approx}&\lambda_{-1}\mathbf{1}_{\{l=0\}}+\eps^{2/3}\frac{d^l}{dy_1^l}\left(y_1/(2\alpha_2)\right)+\eps^{2/3}\sum_{n=1}^N\eps^{2n/3}\sum_{1\leq
  m\leq(n+1)/3}
\widetilde{L}_{n,m-1}\frac{d^l}{dy_1^l}\left(y_1^{n+1-3m}\right)\nonumber\\
&\underset{y_1\to
  -\infty}{\approx}&\lambda_{-1}\mathbf{1}_{\{l=0\}}+\sum_{n=0}^N\eps^{2(n+1)/3}\sum_{0\leq
  m\leq(n+1)/3}
\check{\widetilde{L}}_{n,m}\frac{d^l}{dy_1^l}\left(y_1^{n+1-3m}\right)\nonumber\\
&\underset{y_1\to
  -\infty}{\approx}&\lambda_{-1}\mathbf{1}_{\{l=0\}}+\sum_{0\leq n\leq
  N, m\geq 0}\eps^{2(n+1)/3}\check{\widetilde{L}}_{n,m}\frac{d^l}{dy_1^l}\left(y_1^{n+1-3m}\right),
\ee
with
$$\check{\widetilde{L}}_{n,m}=\left\{\begin{array}{lll}1/(2\alpha_2) & {\rm if} & n=m=0 \\
0 & {\rm if} & n\geq 1 {\rm\ and\ }
m=0 {\rm \ or\ } n=0 {\rm\ and\ } m\geq 1\\
\widetilde{L}_{n,m-1} &
{\rm if} & n\geq 1 {\rm\ and\ } 1\leq m\leq (n+1)/3\\
0&{\rm if} & n\geq 1 {\rm\ and\ }  m> (n+1)/3.
 \end{array}\right.$$

\noindent Thus, for $x\in D_1\cap D_2$, throwing away the smallest terms,
\be\label{-724}
\frac{d^l}{dy_1^l}\left(\eps^{2/3}\lambda(y_1)\right)
&\underset{\eps\to  0}{=}&\lambda_{-1}\mathbf{1}_{\{l=0\}}+\eps^{2l/3}\!\!\!\!\!\!\!\!\sum_{\substack{(m,n)\in\N^2\\(2-3\beta)m+\beta
  n\leq \beta N}}\!\!\!\!\!\!\!\!\eps^{2m}
\check{\widetilde{L}}_{n,m}\frac{d^l}{dz^l}\left(z^{1+n-3m}\right)+o_{L^\infty(D_1\cap
  D_2)}(\eps^{\beta(N+1)+(2/3-\beta)l}).\nonumber
\ee
At this point, the calculation becomes similar to the one which was
performed for $y_1\to +\infty$ in the proof of Lemma
\ref{lem-lam-nu}. Indeed, we can deduce like in (\ref{dev-pow-dem})
that for $l=0$,
\be\label{-dev-pow-dem}
\eps^{1/3}\lambda(y_1)^{1/2}
&\underset{\eps\to 0}{=}&\lambda_{-1}^{1/2}+\sum_{\substack{(n,m)\in\N^2\\(2-3\beta)m+\beta
    n\leq\beta N}}\eps^{2m}z^{n-3m+1}\widetilde{l}_{m,n}+o_{L^\infty(D_1\cap
  D_2)}(\eps^{\beta(N+1)}),
\ee
where
\be\label{def-lmnt}
\widetilde{l}_{m,n}=\sum_{k=1}^{n+1}c_k\sum_{\substack{\left((n_1,m_1),\cdots,(n_k,m_k)\right)\in(\N^2)^k\\ n_1+\cdots+n_k=n-k+1\\m_1+\cdots+m_k=m}}\prod_{j=1}^k\check{\widetilde{L}}_{n_j,m_j}
\ee
for the same coefficients $c_k$ as in (\ref{724}). Note in particular that
$\widetilde{l}_{m,n}=0$ if $m>(n+1)/3$. Indeed, under this condition,
for every $k\in\{1,\cdots,N+1\}$, if $n_1,\cdots,n_k,m_1,\cdots,m_k$
are indices like in the second sum in (\ref{def-lmnt}), we have
$$m_1+\cdots+m_k=m>\frac{n+1}{3}=\frac{(n_1+1)+\cdots+(n_k+1)}{3},$$ 
therefore at least for one of the indices $j\in\{1,\cdots,k\}$, we
have $m_j>(n_j+1)/3$, which implies 
$$\prod_{j=1}^k\check{\widetilde{L}}_{n_j,m_j}=0,$$
for every $k\in \{1,\cdots,N+1\}$, and therefore
$\widetilde{l}_{m,n}=0$. This is the reason why we can add without
changing the result the condition $1+n-3m\geq 0$ in the sum that
appears in (\ref{-d2}) for $l=0$. The proof of (\ref{-d2}) for $l=1$
and $l=2$ is similar to the one which was done on $D_0\cap D_1$ in the
proof of Lemma \ref{lem-lam-nu}.
\end{Proof}

The next lemma provides an asymptotic expansion of $(0,\eps^{1/3}\mu)$ into powers of $\eps$ in $D_1\cap D_2$ as $\eps\to 0$.

\begin{lem}\label{lem-dev-mu+inf}
Let $L\geq 1$, and $\mu$ given by (\ref{troncat}). Then there exists a
family of numbers $(\alpha_{m,n})_{m,n\geq 0}$ such that for every
$l\in \{0,1,2\}$,
\be\label{dev-mu+inf}
\left\|\frac{d^l}{dz^l}\left(\eps^{1/3}\mu(y_2)\right)-\lambda_{-1}^{1/2}\mathbf{1}_{\{l=0\}}-\sum_{\substack{m,n\geq
    0\\\beta n+(2-3\beta)m\leq 2L-\beta\\1+n-3m\geq
    0}}\alpha_{m,n}\eps^{2m}\frac{d^l}{dz^l}\left(z^{1+n-3m}\right)\right\|_{L^\infty(D_1\cap
  D_2)}=o(\eps^{2L-\beta l}).
\ee
\end{lem}

\begin{Proof}
For $x\in D_1\cap D_2$, $(R_2^2-R_1^2)\eps^{-2/3}-\eps^{\beta-2/3}\geq
y_2\geq (R_2^2-R_1^2)\eps^{-2/3}-2\eps^{\beta-2/3}\to +\infty$ as
$\eps\to 0$. Thus, for $l=0,1,2$, thanks to (\ref{mun}),
(\ref{asympnu0+}) and Proposition \ref{proposition-mun}, using for
convenience the notations $g_{0,m}=a_m$, $g_{n,m}^{(0)}=g_{n,m}$,
$g_{n,m}^{(1)}=(1/2-2n-3m)g_{n,m}$ and
$g_{n,m}^{(2)}=(-1/2-2n-3m)g_{n,m}^{(1)}$, we infer
\be
\lefteqn{\frac{d^l}{dz^l}\left(\eps^{1/3}\mu(y_2)\right)}\nonumber\\
 &=&
\eps^{1/3}\sum_{n=0}^L\eps^{2n/3}\frac{d^l}{dz^l}\left(\mu_n(y_2)\right)\nonumber\\
&=&\frac{\eps^{1/3}}{(2\alpha_2)^{1/2}}\sum_{n=0}^L\eps^{2n/3}\sum_{m=0}^\infty g_{n,m}^{(l)}R_2^{2m}\eps^{-2l/3}y_2^{1/2-2n-3m-l}\nonumber\\
&=&\frac{\eps^{1/3}}{(2\alpha_2)^{1/2}}\sum_{n=0}^L\eps^{2n/3}\sum_{m=0}^{L-n}
g_{n,m}^{(l)}R_2^{2m}\eps^{-2l/3}y_2^{1/2-2n-3m-l}+\eps^{1/3}\sum_{n=0}^L\eps^{2(n-l)/3}o_{L^\infty(D_1\cap D_2)}(y_2^{1/2+n-3L-l})\nonumber\\
&=&\frac{\eps^{1/3}}{(2\alpha_2)^{1/2}}\sum_{n=0}^L\eps^{2n/3}\sum_{m=0}^{L-n}
g_{n,m}^{(l)}R_2^{2m}\eps^{-2l/3}\frac{(R_2^2-R_1^2)^{1/2-2n-3m-l}}{\eps^{1/3-4n/3-2m-2l/3}}\left(1+\frac{z}{R_2^2-R_1^2}\right)^{1/2-2n-3m-l}\!\!\!\!\!\!\!\!\!\!\!\!\!\!\!+o_{L^\infty(D_1\cap D_2)}(\eps^{2L})\nonumber\\
&=&\lambda_{-1}^{1/2}\sum_{n=0}^L\sum_{m=0}^{L-n}\eps^{2(n+m)}g_{n,m}^{(l)}R_2^{2m}(R_2^2-R_1^2)^{-2n-3m-l}\left(1+\frac{z}{R_2^2-R_1^2}\right)^{1/2-2n-3m-l}+o_{L^\infty(D_1\cap D_2)}(\eps^{2L})\nonumber\\
&=&\lambda_{-1}^{1/2}\sum_{j=0}^L\eps^{2j}\sum_{\substack{m,n\geq
    0,\\n+m=j}}g_{n,m}^{(l)}R_2^{2m}(R_2^2-R_1^2)^{-2n-3m-l}\sum_{\substack{k\geq
  0,\\\beta k+2j\leq 2L}}c_{k,l,m,n}z^k+o_{L^\infty(D_1\cap D_2)}(\eps^{2L})\\
&=&\lambda_{-1}^{1/2}g_{0,0}^{(l)}(R_2^2-R_1^2)^{-l}+\lambda_{-1}^{1/2}\!\!\!\!\!\sum_{\substack{j,k\geq
    0\\\beta k+2j\leq 2L,\\
(j,k)\neq (0,0)}}\!\!\!\!\!\eps^{2j}\sum_{\substack{m,n\geq
    0,\\n+m=j}}g_{n,m}^{(l)}R_2^{2m}(R_2^2-R_1^2)^{-2n-3m-l}c_{k,l,m,n}z^k+o_{L^\infty(D_1\cap D_2)}(\eps^{2L}),\nonumber
\ee
for some coefficients $(c_{k,l,m,n})_{k\geq 0}$ (with $c_{0,l,m,n}=1$,
$\forall l,m,n$). Then, we
change the variable $k$ in the sum into $p=3j+k-1$. Note that $p\in
\N$ since $(j,k)\in\N^2\string\ \{(0,0)\}$. Thus,
\be\label{pre-dev-mu}
\lefteqn{\frac{d^l}{dz^l}\left(\eps^{1/3}\mu(y_2)\right)
 =\lambda_{-1}^{1/2}g_{0,0}^{(l)}(R_2^2-R_1^2)^{-l}}\\
 &&+\lambda_{-1}^{1/2}\sum_{\substack{j,p\geq
    0\\\beta (p+1)+(2-3\beta)j\leq 2L\\p\geq 3j-1}}c_{1+p-3j}\sum_{\substack{m,n\geq
    0,\\n+m=j}}g_{n,m}^{(l)}R_2^{2m}(R_2^2-R_1^2)^{-2n-3m-l}\eps^{2j}z^{1+p-3j}+o_{L^\infty(D_1\cap D_2)}(\eps^{2L}).\nonumber
\ee
The result follows for $l=0$, since $g_{0,0}=1$, with
$$\alpha_{m,n}= \lambda_{-1}^{1/2}c_{1+n-3m}\sum_{\substack{k,i\geq
    0,\\k+i=m}}g_{i,k}^{(0)}R_2^{2k}(R_2^2-R_1^2)^{-2i-3k}.$$
For $l=1$, (\ref{pre-dev-mu}) gives the existence of some coefficients
$(\alpha_{m,n}')_{m,n}$ such that 
\be
\frac{d}{dz}\left(\eps^{1/3}\mu(y_2)\right)&=&\alpha_{0,0}'+\sum_{\substack{m,n\geq
    0\\\beta (n+1)+(2-3\beta)m\leq 2L\\1+n-3m\geq
    0}}\alpha_{m,n+1}'\eps^{2m}\frac{d}{dz}\left(z^{2+n-3m}\right)+o_{L^\infty(D_1\cap
  D_2)}(\eps^{2L}).\nonumber
\ee
Thus,
\be\label{dev-mu'}
\frac{d}{dz}\left(\eps^{1/3}\mu(y_2)\right)&=&\sum_{\substack{m,n\geq
    0\\\beta n+(2-3\beta)m\leq 2L\\n-3m\geq
    0}}\alpha_{m,n}'\eps^{2m}\frac{d}{dz}\left(z^{1+n-3m}\right)+o_{L^\infty(D_1\cap
  D_2)}(\eps^{2L})\nonumber\\
&=&\sum_{\substack{m,n\geq
    0\\\beta (n+1)+(2-3\beta)m\leq 2L\\n-3m\geq
    0}}\alpha_{m,n}'\eps^{2m}\frac{d}{dz}\left(z^{1+n-3m}\right)+o_{L^\infty(D_1\cap
  D_2)}(\eps^{2L-\beta})\nonumber\\
&=&\sum_{\substack{m,n\geq
    0\\\beta (n+1)+(2-3\beta)m\leq 2L\\1+n-3m\geq
    0}}\alpha_{m,n}'\eps^{2m}\frac{d}{dz}\left(z^{1+n-3m}\right)+o_{L^\infty(D_1\cap
  D_2)}(\eps^{2L-\beta}).
\ee
where in the first equality, we have changed the index of summation
$n$ by $n+1$, in the second equality, we have neglected some terms in
the sum, and in the last equality, the extra term we write in the sum
is in fact equal to 0. In order to prove that (\ref{dev-mu+inf}) also
holds for $l=1$, it remains to prove that for every pair of indices
$(m,n)$ appearing in the sum in (\ref{dev-mu+inf}) (except for
$1+n-3m=0$, for which the corresponding term in (\ref{dev-mu+inf}) for
$l=1$ is anyway equal to 0), we have
$\alpha_{m,n}'=\alpha_{m,n}$. This can be done by using the same trick
as in the proof of Lemma \ref{lem-lam-nu}. Namely, we have on the one
side thanks to (\ref{dev-mu+inf})
\be
\lefteqn{\eps^{1/3}\mu\left(\frac{R_2^2-R_1^2}{\eps^{2/3}}-\eps^{\beta-2/3}\right)-\eps^{1/3}\mu\left(\frac{R_2^2-R_1^2}{\eps^{2/3}}-2\eps^{\beta-2/3}\right)}\nonumber\\
&=&\sum_{\substack{m,n\geq
    0\\\beta n+(2-3\beta)m\leq 2L-\beta\\1+n-3m\geq
    0}}\alpha_{m,n}\eps^{2m}(-1)^{1+n-3m}\eps^{\beta(1+n-3m)}\left(1-2^{1+n-3m}\right)+o_{L^\infty(D_1\cap
  D_2)}(\eps^{2L}),
\ee
and on the other side, by integration of (\ref{dev-mu'}) between
$z=-2\eps^\beta$ and $z=-\eps^\beta$, we have the
same equality with $\alpha_{m,n}$ replaced by $\alpha_{m,n}'$. Since
$\beta$ has been chosen irrational, the linear independance of the
functions $\eps\mapsto\eps^{(2-3\beta)m+\beta (n+1)}$ implies that for
all the indices $(m,n)$ appearing in the sum (except for $1+n-3m=0$),
we have $\alpha_{m,n}=\alpha_{m,n}'$. The proof of (\ref{dev-mu+inf})
for $l=2$ is similar. 
\end{Proof}

\noindent The next lemma shows that the expansions of $\eps^{1/3}(\nu,\lambda^{1/2})$ and
  $(0,\eps^{1/3}\mu)$ calculated respectively in Lemmata \ref{-lem-lam-nu}  and \ref{lem-dev-mu+inf} are in fact the same.

\begin{lem}\label{lem-compar2}
Let $N\geq 1$ be an integer, $\eps_0>0$, and
$\beta\in (0,2/3)\string\ \mathbb{Q}$. Let $(\theta)_{0<\eps\leq\eps_0}$ be a sequence of regular functions defined for $z\in [-2\eps^\beta,-\eps^\beta]$ such that
\be\label{2eq-theta2}
\left\|\eps^2\Delta\theta+(R_2^2-R_1^2)\theta+z\theta-2\alpha_2\theta^3\right\|_{L^\infty(D_1\cap D_2)}=o\left(\eps^{\beta(N+1)}\right).
\ee
We assume that there exists a family of real numbers $q_{m,n}$, defined for every
  $(m,n)\in\N^2$ such that $(2-3\beta)m+\beta n\leq\beta N$, such that for 
$l\in\{0,1,2\}$, we have
\be\label{2theta2}
&&\!\!\!\!\!\left\|\theta^{(l)}-\lambda_{-1}^{1/2}\mathbf{1}_{\{l=0\}}-\!\!\!\!\!\!\sum_{\substack{(m,n)\in\N^2\\ (2-3\beta)m+\beta
    n\leq \beta N}}\!\!\!\!\eps^{2m}q_{m,n}\ \frac{d^l}{dz^l}\left(z^{1+n-3m}\right)\right\|_{L^\infty(D_1\cap
  D_2)}\!\!\!\!\!\!\!\!\!\!\!\underset{\eps\to
  0}{=}o(\eps^{\beta(N+1-l)}).\ \ \ \ \ 
\ee
Then,
equations (\ref{2theta2}) and
(\ref{2eq-theta2}) entirely determine the values of the $q_{m,n}$'s for $(2-3\beta)m+\beta
n\leq\beta (N-1)$. Moreover, these coefficients do not depend on $N$
or $\beta$.
\end{lem}

\begin{Proof}
For convenience, for every $(m,n)\in\N^2$, we denote
$q_{m,n}'=(1+n-3m)q_{m,n}$ and $q_{m,n}''=(n-3m)(1+n-3m)q_{m,n}$. For
a function
$\theta$ that satisfies (\ref{2theta2}), let
us calculate the function that appears in the left hand side of
(\ref{2eq-theta2}), evaluated at $z=-\eps^\beta$. We have
\be
\lefteqn{\eps^2\Delta\theta+(R_2^2-R_1^2)\theta+z\theta-2\alpha_2\theta^3}\nonumber\\
&=&-2\sum_{\substack{m,n\geq 0\\(2-3\beta)m+\beta n\leq \beta N}}(d
q_{m,n}'+2q_{m,n}'')(-1)^{n-3m}\eps^{(2-3\beta)(m+1)+\beta(
  n+2)+\beta}\nonumber\\
&&+4R_1^2\sum_{\substack{m,n\geq 0\\(2-3\beta)m+\beta
    n\leq \beta N}}q_{m,n}''(-1)^{n-3m-1}\eps^{(2-3\beta)(m+1)+\beta(
  n+1)+\beta}\nonumber\\
&&+(R_2^2-R_1^2)\lambda_{-1}^{1/2}+(R_2^2-R_1^2)\sum_{\substack{m,n\geq 0\\(2-3\beta)m+\beta
    n\leq \beta N}}q_{m,n}(-1)^{n-3m+1}\eps^{(2-3\beta)m+\beta
  n+\beta}\nonumber\\
&&-\eps^\beta\lambda_{-1}^{1/2}+\sum_{\substack{m\geq 0,n\geq 0\\(2-3\beta)m+\beta
    n\leq \beta N}}q_{m,n}(-1)^{n-3m}\eps^{(2-3\beta)m+\beta
  (n+1)+\beta}\nonumber\\
&&-2\alpha_2\lambda_{-1}^{3/2}-6\alpha_2\lambda_{-1}\sum_{\substack{m\geq 0,n\geq 0\\(2-3\beta)m+\beta
    n\leq \beta N}}q_{m,n}(-1)^{n-3m+1}\eps^{(2-3\beta)m+\beta
  n+\beta}\nonumber\\
&&-6\alpha_2\lambda_{-1}^{1/2}\sum_{\substack{m,n\geq 0\\(2-3\beta)m+\beta
    n\leq \beta N}}\big(\sum_{\substack{m_1,m_2,n_1,n_2\geq
  0\\m_1+m_2=m\\n_1+n_2=n}}q_{m_1,n_1}q_{m_2,n_2}\big)(-1)^{n-3m}\eps^{(2-3\beta)m+\beta (n+1)+\beta}\nonumber\\
&&-2\alpha_2\sum_{\substack{m\geq
    0,n\geq 0\\(2-3\beta)m+\beta
    n\leq \beta N}}\big(\sum_{\substack{m_1,m_2,m_3,n_1,n_2,n_3\geq
  0,\\n_1+n_2+n_3=n\\m_1+m_2+m_3=m}}q_{m_1,n_1}q_{m_2,n_2}q_{m_3,n_3}\big)(-1)^{n-3m+1}\eps^{(2-3\beta)m+\beta (n+2)+\beta}\nonumber
\ee
Thus, changing the indices and throwing away all the terms that can be
incorporated in the rest,
\be\label{2dev-eq-theta2}
\lefteqn{\eps^2\Delta\theta+(R_2^2-R_1^2)\theta+z\theta-2\alpha_2\theta^3-2\alpha_0\theta_1^2\theta}\nonumber\\
&=&-2\sum_{\substack{m\geq 1,n\geq 2\\(2-3\beta)m+\beta n\leq \beta N}}(d
q_{m-1,n-2}'+2q_{m-1,n-2}'')(-1)^{n-3m+1}\eps^{(2-3\beta)m+\beta
  n+\beta}\nonumber\\
&&+4R_1^2\sum_{\substack{m,n\geq 1\\(2-3\beta)m+\beta
    n\leq \beta N}}q_{m-1,n-1}''(-1)^{n-3m+1}\eps^{(2-3\beta)m+\beta
  n+\beta}\nonumber\\
&&+(R_2^2-R_1^2)\sum_{\substack{m,n\geq 0\\(2-3\beta)m+\beta
    n\leq \beta N}}q_{m,n}(-1)^{n-3m+1}\eps^{(2-3\beta)m+\beta
  n+\beta}\nonumber\\
&&-\eps^\beta\lambda_{-1}^{1/2}+\sum_{\substack{m\geq 0,n\geq 1\\(2-3\beta)m+\beta
    n\leq \beta N}}q_{m,n-1}(-1)^{n-3m+1}\eps^{(2-3\beta)m+\beta
  n+\beta}\nonumber\\
&&-6\alpha_2\lambda_{-1}\sum_{\substack{m\geq 0,n\geq 0\\(2-3\beta)m+\beta
    n\leq \beta N}}q_{m,n}(-1)^{n-3m+1}\eps^{(2-3\beta)m+\beta
  n+\beta}\nonumber\\
&&-6\alpha_2\lambda_{-1}^{1/2}\sum_{\substack{m\geq 0,n\geq 1\\(2-3\beta)m+\beta
    n\leq \beta N}}\big(\sum_{\substack{m_1,m_2,n_1,n_2\geq
  0\\m_1+m_2=m\\n_1+n_2=n-1}}q_{m_1,n_1}q_{m_2,n_2}\big)(-1)^{n-3m+1}\eps^{(2-3\beta)m+\beta n+\beta}\nonumber\\
&&-2\alpha_2\sum_{\substack{m\geq
    0,n\geq 2\\(2-3\beta)m+\beta
    n\leq \beta N}}\big(\sum_{\substack{m_1,m_2,m_3,n_1,n_2,n_3\geq
  0,\\n_1+n_2+n_3=n-2\\m_1+m_2+m_3=m}}q_{m_1,n_1}q_{m_2,n_2}q_{m_3,n_3}\big)(-1)^{n-3m+1}\eps^{(2-3\beta)m+\beta n+\beta}.\nonumber
\ee
According to (\ref{2eq-theta2}), the right hand side of
(\ref{2dev-eq-theta2}) is equal to 0, up to the rest term
$o(\eps^{\beta(N+1)})$. Thus, the linear independance of the family of
functions of $\eps$, $\left(\eps^{(2-3\beta)m+\beta n}\right)_{m,n\geq
  0}$ yields:
\begin{itemize}
\item for $m=0$, $n=0$, thanks to (\ref{I-1}), we get
\be\label{2II00}
2(R_2^2-R_1^2)q_{0,0}=\lambda_{-1}^{1/2},
\ee
\item for $m=0$ and $1\leq n\leq N$,
\begin{eqnarray*}
\!\!\!\!\!\!\!\!\!\!\!\!\!\!\!(R_2^2-R_1^2)q_{0,n}+q_{0,n-1}-6\alpha_2\lambda_{-1}q_{0,n}-6\alpha_2\lambda_{-1}^{1/2}\sum_{\substack{n_1,n_2\geq
      0,\\n_1+n_2=n-1}}q_{0,n_1}q_{0,n_2}-2\alpha_2\sum_{\substack{n_1,n_2,n_3\geq
      0,\\n_1+n_2+n_3=n-2}}q_{0,n_1}q_{0,n_2}q_{0,n_3}=0,
\end{eqnarray*}
which, using (\ref{I-1}), can be rewritten as
\be\label{2II0n}
\!\!\!\!\!\!-2(R_2^2-R_1^2)q_{0,n}
&=&-q_{0,n-1}+6\alpha_2\lambda_{-1}^{1/2}\sum_{\substack{n_1,n_2\geq
    0,\\n_1+n_2=n-1}}q_{0,n_1}q_{0,n_2}+2\alpha_2\sum_{\substack{n_1,n_2,n_3\geq
    0,\\n_1+n_2+n_3=n-2}}q_{0,n_1}q_{0,n_2}q_{0,n_3},\ \ \ \ \ \ \ \ 
\ee
\item for $1\leq m\leq \beta N/(2-3\beta)$ and $n=0$, we get
\be\label{2IIm0}
q_{m,0}&=&0.
\ee
\item for $m\geq 1$ and $n\geq 1$ such that $(2-3\beta)m+\beta n\leq \beta N$,
\begin{eqnarray*}
-2(dq_{m-1,n-2}'+2q_{m-1,n-2}'')\mathbf{1}_{\{n\geq 2\}}+4R_1^2q_{m-1,n-1}''+(R_2^2-R_1^2)q_{m,n}+q_{m,n-1}\nonumber\\
-6\alpha_2\lambda_{-1}q_{m,n}-6\alpha_2\lambda_{-1}^{1/2}\sum_{\substack{m_1,m_2,n_1,n_2\geq
  0\\m_1+m_2=m,\\n_1+n_2=n-1}}q_{m_1,n_1}q_{m_2,n_2}-2\alpha_2\sum_{\substack{m_1,m_2,m_3,n_1,n_2,n_3\geq
  0,\\n_1+n_2+n_3=n-2\\m_1+m_2+m_3=m}}q_{m_1,n_1}q_{m_2,n_2}q_{m_3,n_3}=0,\nonumber
\end{eqnarray*}
which can be rewritten as
\be\label{2IImn}
\lefteqn{-2(R_2^2-R_1^2)q_{m,n}\ =\ 2(dq_{m-1,n-2}'+2q_{m-1,n-2}'')\mathbf{1}_{\{n\geq 2\}}-4R_1^2q_{m-1,n-1}''-q_{m,n-1}}\nonumber\\
&&
+6\alpha_2\lambda_{-1}^{1/2}\sum_{\substack{0\leq m_1,m_2,n_1,n_2\\n_1+n_2=n-1\\m_1+m_2=m}}q_{m_1,n_1}q_{m_2,n_2}+2\alpha_2\sum_{\substack{m_1,m_2,m_3,n_1,n_2,n_3\geq
  0\\n_1+n_2+n_3=n-2\\m_1+m_2+m_3=m}}q_{m_1,n_1}q_{m_2,n_2}q_{m_3,n_3}
\ee
\end{itemize}
From (\ref{2II00}), (\ref{2II0n}), (\ref{2IIm0}) and (\ref{2IImn}), it
clearly follows that all the $q_{m,n}$'s for indices $(m,n)$ that
satisfy $(2-3\beta)m+\beta n\leq \beta N$ are completely determined.
\end{Proof}

\noindent Finally, we show that
$(\eps^{1/3}\nu(y_1),\eps^{1/3}\lambda(y_1)^{1/2})$ and
$(0,\eps^{1/3}\mu(y_2))$ are close one from another on $D_1\cap D_2$.

\begin{lem}\label{lem-compar-D1D2}
Let $N\geq 1$ be an integer, $L\geq \beta(N+1)/2$  and $\nu$,
$\lambda$, $\mu$ given by (\ref{troncat}). Then for $l\in\{0,1,2\}$,
\be\label{comp-0-nu}
\forall\alpha>0, \quad \left\|\frac{d^l}{dz^l}\left(\eps^{1/3}\nu(y_1)\right)\right\|_{L^\infty(D_1\cap
  D_2)}=o(\eps^\alpha)
\ee
and
\be\label{comp-lam-mu}
\left\|\frac{d^l}{dz^l}\left(\eps^{1/3}\lambda(y_1)^{1/2}-\eps^{1/3}\mu(y_2)\right)\right\|_{L^\infty(D_1\cap
  D_2)}=o(\eps^{\beta(N+1-l)}).
\ee
\end{lem}

\begin{Proof}
(\ref{comp-0-nu}) has already been proved in Lemma
  \ref{-lem-lam-nu}. $\theta=\eps^{1/3}\lambda(y_1)^{1/2}$ satisfies
  assumption (\ref{2theta2}) in Lemma \ref{lem-compar2} thanks to
  Lemma (\ref{-lem-lam-nu}) (with $q_{m,n}=0$ if
  $1+n-3m<0$). $\eps^{1/3}\lambda^{1/2}$ also satisfies the assumption
  (\ref{2eq-theta2}) thanks to Lemma \ref{lem-tronc-nu-l} and (\ref{comp-0-nu}).
The two assumptions (\ref{2theta2}) and (\ref{2eq-theta2}) of Lemma
\ref{lem-compar2} are also satisfied by $\theta=\eps^{1/3}\mu(y_2)$,
thanks respectively to Lemma \ref{lem-dev-mu+inf} and Corollary \ref{cor2}. Therefore, thanks to Lemma \ref{lem-compar2},
(\ref{-d2}) and (\ref{dev-mu+inf}), we deduce (\ref{comp-lam-mu}).
\end{Proof}

\section{Proof of Theorem \ref{main}}\label{rigorous}
\subsection{Derivation of the equations}\label{sec-deriv}
We look for solutions of (\ref{sys}) under the form given by the
ansatz (\ref{ansatz}), where $\beta\in (0,2/3)\backslash \Q$, $N$ is a
large integer, 
$M\geq \max(1,\beta N/(2-3\beta)$ and $L\geq \max(1,\beta(N+1)/2$. For the sake of simplicity, we rewrite this
ansatz as
\be\label{ansatz-rho}
\eta_1 &=&\eps^{1/3}\left(\rho_1+\eps^{2(N+1)/3}P\right),\\
\eta_2 &=&\eps^{1/3}\left(\rho_2+\eps^{2(N+1)/3}Q\right),\label{ansatz-rho2}
\ee
where
\be\label{def-rho}
\rho_1&=& \Phi_\eps\eps^{-1/3}\omega+\chi_\eps\nu,\\
\rho_2 &=&
\Phi_\eps\eps^{-1/3}\tau+\chi_\eps\lambda^{1/2}+\Psi_\eps\mu.\label{def-rho2}
\ee
Implicitely, $\rho_1$, $\rho_2$, $P$ and $Q$ are functions of
$x\in\R^d$, $\omega$ and $\tau$ are functions of $z=R_1^2-|x|^2$,
$\nu$ and $\lambda$ are functions of $y_1=z/\eps^{2/3}$ and $\mu$ is a
function of the variable $y_2=(R_2^2-|x|^2)/\eps^{2/3}$. $\nabla$ and $\Delta$ refer to derivatives with respect to $x\in \R^d$, whereas primes refer to derivatives with respect to variables $z$, $y_1$ or $y_2$, depending on the function which is concerned. For instance, we note $\nabla \omega$ for $\nabla\omega=-2x\omega'(R_1^2-|x|^2)=-2x\omega'(z)$.
Using this ansatz and these notations, the first equation in (\ref{sys}) becomes
\be
\lefteqn{\eps^{4/3}\Delta\rho_1+\eps^{2N/3+2}\Delta
P+\left(\frac{\alpha_0}{\alpha_2}\frac{R_2^2-R_1^2}{\eps^{2/3}}+y_1\right)(\rho_1+\eps^{2(N+1)/3}P)}\nonumber\\
&&-2\alpha_1(\rho_1+\eps^{2(N+1)/3}P)^3-2\alpha_0(\rho_2+\eps^{2(N+1)/3}Q)^2(\rho_1+\eps^{2(N+1)/3}P)=0.\nonumber
\ee
Reorganizing the different terms, we get
\be\label{eq-I-PQ}
\eps^{4/3}\Delta\rho_1
+\frac{\alpha_0}{\alpha_2}\frac{R_2^2-R_1^2}{\eps^{2/3}}\rho_1+y_1\rho_1
-2\alpha_1\rho_1^3-2\alpha_0\rho_2^2\rho_1\nonumber\\ 
+\eps^{2(N+1)/3}\left(
\eps^{4/3}\Delta
P+\frac{\alpha_0}{\alpha_2}\frac{R_2^2-R_1^2}{\eps^{2/3}}P
+y_1P-6\alpha_1\rho_1^2P-2\alpha_0\rho_2^2P-4\alpha_0\rho_1\rho_2 Q\right)\nonumber\\ 
+\eps^{4(N+1)/3}\left(-6\alpha_1\rho_1 P^2-4\alpha_0\rho_2
PQ-2\alpha_0\rho_1 Q^2\right)+\eps^{2(N+1)}\left(-2\alpha_1P^3-2\alpha_0PQ^2 \right)&=&0.
\ee
Similarly, the second equation in (\ref{sys}) writes
\be\label{eq-II-PQ}
\eps^{4/3}\Delta\rho_2
+y_2\rho_2
-2\alpha_2\rho_2^3-2\alpha_0\rho_1^2\rho_2\nonumber\\ 
+\eps^{2(N+1)/3}\left(
\eps^{4/3}\Delta
Q+y_2Q-6\alpha_2\rho_2^2Q-2\alpha_0\rho_1^2Q-4\alpha_0\rho_1\rho_2 P\right)\nonumber\\ 
+\eps^{4(N+1)/3}\left(-6\alpha_2\rho_2 Q^2-4\alpha_0\rho_1
PQ-2\alpha_0\rho_2 P^2\right)+\eps^{2(N+1)}\left(-2\alpha_2Q^3-2\alpha_0P^2Q \right)&=&0.
\ee
Equations (\ref{eq-I-PQ}) and (\ref{eq-II-PQ}) can be rewritten as the
system
\be\label{eq-reste}
A_\eps \left[\begin{array}{c}
P\\Q
\end{array}\right]=f_\eps^0(x)+f_\eps^2(x,P,Q)+f_\eps^3(x,P,Q),
\ee
where
$$A_\eps=\left[\begin{array}{cc}
-\eps^{4/3} \Delta+p_\eps(x)&r_\eps(x)\\
r_\eps(x)&-\eps^{4/3} \Delta+q_\eps(x)
\end{array}\right],$$
$$p_\eps(x)=-\frac{\alpha_0}{\alpha_2}\frac{R_2^2-R_1^2}{\eps^{2/3}}-y_1+6\alpha_1\rho_1^2+2\alpha_0\rho_2^2,$$
$$q_\eps(x)=-y_2+6\alpha_2\rho_2^2+2\alpha_0\rho_1^2,$$
$$r_\eps(x)=4\alpha_0\rho_1\rho_2,$$
\be
f_\eps^0(x)
&=&\eps^{-2(N+1)/3}\left[\begin{array}{c}
\eps^{4/3}\Delta\rho_1
+\frac{\alpha_0}{\alpha_2}\frac{R_2^2-R_1^2}{\eps^{2/3}}\rho_1+y_1\rho_1
-2\alpha_1\rho_1^3-2\alpha_0\rho_2^2\rho_1\\
\eps^{4/3}\Delta\rho_2
+y_2\rho_2
-2\alpha_2\rho_2^3-2\alpha_0\rho_1^2\rho_2
\end{array}\right],\nonumber
\ee
$$f_\eps^2(x,P,Q)=-2\eps^{2(N+1)/3}\left[\begin{array}{c}3\alpha_1\rho_1
P^2+2\alpha_0\rho_2 PQ+\alpha_0\rho_1
Q^2\\
3\alpha_2\rho_2 Q^2+2\alpha_0\rho_1
  PQ+\alpha_0\rho_2 P^2\end{array}\right],$$
  $$f_\eps^3(x,P,Q)=-2\eps^{4(N+1)/3}\left[\begin{array}{c}\alpha_1 P^3+\alpha_0 PQ^2\\\alpha_2
  Q^3+\alpha_0 P^2Q\end{array}\right].$$

\subsection{Estimate on the source term $f_\eps^0$}
Equation (\ref{eq-reste}) will be solved thanks to a fixed point
argument. For this purpose, we need to show that the source term
$f_\eps^0$ is small if functions $\omega$, $\tau$,
$\nu$, $\lambda$ and $\mu$ are given by (\ref{troncat}).
The first component of $f_\eps^0$ can be rewritten as
\be\label{f01}
\left[f_\eps^0\right]_1&=&\eps^{-2(N+1)/3}\left[\eps^{4/3}\Delta\rho_1
+\frac{\alpha_0}{\alpha_2}\frac{R_2^2-R_1^2}{\eps^{2/3}}\rho_1+y_1\rho_1
-2\alpha_1\rho_1^3-2\alpha_0\rho_2^2\rho_1\right]\nonumber\\
&=&\Phi_\eps\underbrace{\eps^{-2(N+1)/3}\eps^{-1}\left[\eps^2\Delta\omega+\frac{\alpha_0}{\alpha_2}(R_2^2-R_1^2)\omega+z\omega-2\alpha_1\omega^3-2\alpha_0\tau^2\omega\right]}_{g_0}\nonumber\\
&&+\chi_\eps\underbrace{\eps^{-2(N+1)/3}\left[\eps^{4/3}\Delta\nu+\frac{\alpha_0}{\alpha_2}\frac{R_2^2-R_1^2}{\eps^{2/3}}\nu+y_1\nu-2\alpha_1\nu^3-2\alpha_0\lambda\nu\right]}_{g_1}\nonumber\\
&&+\underbrace{2\eps^{-2(N-1)/3}\left[\nabla\Phi_\eps\nabla(\eps^{-1/3}\omega)+\nabla\chi_\eps\nabla\nu\right]}_{k_1}+\underbrace{\eps^{-2(N-1)/3}\left[\Delta\Phi_\eps(\eps^{-1/3}\omega)+\Delta\chi_\eps\nu\right]}_{k_2}\nonumber\\
&&+\underbrace{2\alpha_1\eps^{-2(N+1)/3}\left[\eps^{-1}\Phi_\eps\omega^3+\chi_\eps\nu^3-\rho_1^3\right]}_{l_1}+\underbrace{2\alpha_0\eps^{-2(N+1)/3}\left[\eps^{-1}\Phi_\eps\tau^2\omega+\chi_\eps\lambda\nu-\rho_1\rho_2^2\right]}_{l_2}.\ \ \ \ 
\ee
As for the second component of $f_\eps^0$, we have
\be\label{dec_f02}
\lefteqn{\left[f_\eps^0\right]_2  = \eps^{-2(N+1)/3}\left[\eps^{4/3}\Delta\rho_2
+y_2\rho_2
-2\alpha_2\rho_2^3-2\alpha_0\rho_1^2\rho_2\right]}\nonumber\\
&=&\Phi_\eps\underbrace{\eps^{-2(N+1)/3}\eps^{-1}\left[\eps^2\Delta\tau+(R_2^2-R_1^2+z)\tau-2\alpha_2\tau^3-2\alpha_0\omega^2\tau\right]}_{h_0}\\
&&+\chi_\eps\underbrace{\eps^{-2(N+1)/3}\left[\eps^{4/3}\Delta(\lambda^{1/2}) +y_2\lambda^{1/2}-2\alpha_2\lambda^{3/2}-2\alpha_0\nu^2\lambda^{1/2}\right]}_{h_1}
+\Psi_\eps\underbrace{\eps^{-2(N+1)/3}\left[\eps^{4/3}\Delta\mu+y_2\mu-2\alpha_2\mu^3\right]}_{h_2}\nonumber\\
&&+\underbrace{\eps^{-2(N-1)/3}\left[\Delta\Phi_\eps\eps^{-1/3}\tau+\Delta\chi_\eps\lambda^{1/2}+\Delta\Psi_\eps\mu\right]}_{k_3}+\underbrace{2\eps^{-2(N-1)/3}\left[\nabla\Phi_\eps\eps^{-1/3}\nabla\tau+\nabla\chi_\eps\nabla\lambda^{1/2}+\nabla\Psi_\eps\nabla\mu\right]}_{k_4}\nonumber\\
&&+\underbrace{2\alpha_2\eps^{-2(N+1)/3}\left[\eps^{-1}\Phi_\eps\tau^3+\chi_\eps\lambda^{3/2}+\Psi_\eps\mu^3-\rho_2^3\right]}_{l_3}+\underbrace{2\alpha_0\eps^{-2(N+1)/3}\left[\eps^{-1}\Phi_\eps\omega^2\tau+\chi_\eps\nu^2\lambda^{1/2}-\rho_1^2\rho_2\right]}_{l_4}.\nonumber
\ee
Thanks to Lemma \ref{est-eq-om-tau}, for $x\in {\rm Supp}\Phi_\eps\subset
D_0$, we have
\be\label{est-g0-h0}
|g_0|\lesssim\eps^{(2-3\beta)(M+1/2)-2N/3-2/3}\quad\text{and}\quad
|h_0|\lesssim\eps^{(2-3\beta)(M+1/2)-2N/3-2/3-\beta/2}.
\ee
From Lemma \ref{lem-tronc-nu-l}, for $x\in {\rm Supp}\chi_\eps\subset
D_1$, we obtain
\be\label{est-g1-h1}
|g_1|\lesssim\eps^{-(2-3\beta)(N/3-7/6)}\quad\text{and}\quad
|h_1|\lesssim\eps^{-(2-3\beta)N/3+1/3-\beta}.
\ee
Lemma \ref{lem-est-h2-D2} yields, for $x\in {\rm Supp}\Psi_\eps\subset
D_2$,
\be\label{est-h2}
|h_2|\lesssim\eps^{2(L-N)/3}h(x).
\ee
Next, let us estimate $k_1$. Note that $\nabla\Phi_\eps$ is supported
in $D_0\cap D_1$, whereas $\nabla\chi_\eps$ is supported
in $(D_0\cap D_1)\cup (D_1\cap D_2)$. Moreover, for $x\in D_0\cap D_1$, we
have
$$\nabla\Phi_\eps=-\nabla \chi_\eps=-2x\eps^{-\beta}\phi'\left(\frac{z-\eps^\beta}{2\eps^\beta-\eps^\beta}\right).$$
Thus,
$$|k_1|\lesssim
\eps^{-2N/3+1/3-\beta}\left\|\nabla(\omega-\eps^{1/3}\nu)\right\|_{L^\infty(D_0\cap
  D_1)}\mathbf{1}_{D_0\cap
  D_1}+\eps^{-2N/3+2/3-\beta}\left\|\nabla\nu\right\|_{L^\infty(D_1\cap D_2)}\mathbf{1}_{D_1\cap
  D_2}.$$
Then, thanks to Lemma \ref{lem-compar-om-nu-tau-l} and Lemma
\ref{lem-compar-D1D2}, 
\be\label{est-k1}
k_1=o_{L^\infty((D_0\cap D_1)\cup (D_1\cap D_2))}\left(\eps^{-(2-3\beta)N/3-5\beta/2+1/3}\right).
\ee
Similarly,
\be\label{est-k4}
\lefteqn{\|k_4\|_{L^\infty((D_0\cap D_1)\cup (D_1\cap D_2))}}\\
&\lesssim
&\eps^{-2N/3+2/3}\left(\eps^{-1/3-\beta}\left\|\frac{d}{dz}\left(\tau-\eps^{1/3}\lambda^{1/2}\right)\right\|_{L^\infty(D_0\cap
    D_1)}+\eps^{-1/3-\beta}\left\|\frac{d}{dz}\left(\eps^{1/3}\lambda^{1/2}-\eps^{1/3}\mu\right)\right\|_{L^\infty(D_1\cap
    D_2)}\right)\nonumber\\
&=& o\left(\eps^{-(2-3\beta)N/3-2\beta+1/3}\right),
\ee
and we also get similar estimates for $k_2$ and $k_3$:
\be\label{est-k2-k3}
k_2&=&o_{L^\infty((D_0\cap D_1)\cup (D_1\cap D_2))}(\eps^{-(2-3\beta)N/3-5\beta/2+1/3}),\\
k_3&=&o_{L^\infty((D_0\cap D_1)\cup (D_1\cap D_2))}(\eps^{-(2-3\beta)N/3-2\beta+1/3}).
\ee
Next, we estimate $l_1$. Clearly, $l_1$ is supported in $D_0\cap
D_1$. Moreover, Lemma \ref{lem-compar-om-nu-tau-l} implies 
$$\eps^{1/3}\nu=\omega+o_{L^\infty(D_0\cap D_1)}(\eps^{\beta(N-1/2)}),$$
and since $\eps^\beta\leq |z|\leq 2\eps^\beta$ for $x\in D_0\cap D_1$, it follows from the definition of $\omega$ given by
(\ref{troncat}), (\ref{eq-omega0}) and (\ref{asymp-omegan}) and from the asymptotics of the $\omega_m$'s as $z\to 0$ given in (\ref{eq-omega0}) and (\ref{asymp-omegan}) that
$$\|\omega\|_{L^\infty(D_0\cap D_1)}=O(\eps^{\beta/2}).$$
Thus, on $D_0\cap D_1$, we get
\be\label{est-l1}
l_1 &=&2\alpha_1
\eps^{-2N/3-5/3}\left[\Phi_\eps\omega^3+(1-\Phi_\eps)\left(\omega+o_{L^\infty(D_0\cap
  D_1)}(\eps^{\beta(N-1/2)})\right)^3\right.\nonumber\\
&&\left.-\left(\Phi_\eps\omega+(1-\Phi_\eps)(\omega+o_{L^\infty(D_0\cap
  D_1)}(\eps^{\beta(N-1/2)})\right)^3\right]\nonumber\\
&=&2\alpha_1
\eps^{-2N/3-5/3}o_{L^\infty(D_0\cap
  D_1)}(\eps^{\beta(N+1/2)})\quad=\quad o_{L^\infty(D_0\cap
  D_1)}\left(\eps^{-(2-3\beta)N/3+\beta/2-5/3}\right).
\ee
As for $l_2$, it is supported in $(D_0\cap D_1)\cup (D_1\cap
D_2)$. Taking into account Lemma \ref{lem-compar-om-nu-tau-l} and
Lemma \ref{lem-compar-D1D2}, $l_2$
can be rewritten as 
\be\label{est-l2}
l_2& =& \eps^{-2(N+1)/3}\eps^{-1}\left[\left\{\Phi_\eps\tau^2\omega+(1-\Phi_\eps)\left(\tau+o(\eps^{\beta
    N})\right)^2\left(\omega+o(\eps^{\beta
    (N-1/2)})\right)\right.\right.\nonumber\\
&&\left.\left.-\left(\Phi_\eps\omega+(1-\Phi_\eps)(\omega+o(\eps^{\beta
    (N-1/2)}))\right)\left(\Phi_\eps\tau+(1-\Phi_\eps)(\tau+o(\eps^{\beta
    N}))\right)^2\right\}\mathbf{1}_{D_0\cap
    D_1}\right.\nonumber\\
&&\left.+\eps\left\{\chi_\eps\lambda\nu-\chi_\eps\nu\left(\chi_\eps\lambda^{1/2}+(1-\chi_\eps)\mu\right)^2\right\}\mathbf{1}_{D_1\cap
    D_2}\right]\nonumber\\
&=&o_{L^\infty(D_0\cap D_1)}(\eps^{-(2-3\beta)N/3-5/3-\beta/2})\mathbf{1}_{D_0\cap D_1}+o_{L^\infty(D_1\cap D_2)}(\eps^\alpha)\mathbf{1}_{D_1\cap D_2},
\ee
where $\alpha$ is arbitrarily large. Similar calculations yield
\be\label{est-l3}
l_3&=&o_{L^\infty(D_0\cap
  D_1)}(\eps^{-(2-3\beta)N/3-5/3})\mathbf{1}_{D_0\cap D_1}+o_{L^\infty(D_1\cap
  D_2)}(\eps^{-(2-3\beta)N/3-5/3+\beta})\mathbf{1}_{D_1\cap D_2}
\ee
and
\be\label{est-l4}
l_4&=&o_{L^\infty(D_0\cap
  D_1)}(\eps^{-(2-3\beta)N/3-5/3})\mathbf{1}_{D_0\cap D_1}+o_{L^\infty(D_1\cap
  D_2)}(\eps^{\alpha})\mathbf{1}_{D_1\cap D_2}.
\ee
Combining all these inequalities and noting that the measure of $D_1$ is of the order of $\eps^\beta$,
 we deduce
\be\label{est-fe01}
\lefteqn{\left\|\left[f_\eps^0\right]_1\right\|_{L^2(\R^d)}}\nonumber\\
&\lesssim
 &\eps^{(2-3\beta)(M+1/2)-2N/3-2/3}+\eps^{-(2-3\beta)N/3+7/3-3\beta}+\eps^{-(2-3\beta)N/3+1/3-2\beta}+\eps^{-(2-3\beta)N/3-5/3}
\ee
and
\be\label{est-fe02}
\left\|\left[f_\eps^0\right]_2\right\|_{L^2(\R^d)}&\lesssim
  &\eps^{(2-3\beta)(M+1/2)-2N/3-2/3-\beta/2}+\eps^{-(2-3\beta)N/3+1/3-\beta/2}\nonumber\\
&&+\eps^{2(L-N)/3}+\eps^{-(2-3\beta)N/3+1/3-3\beta/2}+\eps^{-(2-3\beta)N/3-5/3+\beta/2}\nonumber\\
&\lesssim
  &\eps^{(2-3\beta)(M+1/2)-2N/3-2/3-\beta/2}+\eps^{2(L-N)/3}+\eps^{-(2-3\beta)N/3-5/3+\beta/2}.
\ee

\subsection{Estimate on the resolvent of $A_\eps$}
In order to solve equation (\ref{eq-reste}) with the choice of
$\nu,\mu,\lambda$ given in (\ref{troncat}), we have to invert $A_\eps$. For this
purpose, we prove that $A_\eps$ is a positive self-adjoint operator on
$L^2(\R^d)$. It will be convenient to have an idea about the size of
the functions $p_\eps$, $q_\eps$ and $r_\eps$ appearing in the
expression of $A_\eps$, depending on $x$. As a preliminary, let us
first simplify the expressions of $\rho_1^2$ and $\rho_2^2$, depending
on whether $x\in D_0$, $x\in D_1\string\(D_0\cup D_2)$ or $x\in
D_2$. Thanks to Lemma \ref{lem-compar-om-nu-tau-l}, (\ref{eq-omega0})
and (\ref{asymp-omegan}), we have, for $x\in D_0$,
\be\label{rho1d0}
\rho_1^2&=&\frac{1}{\eps^{2/3}}\left(\Phi_\eps\omega+\chi_\eps\eps^{1/3}\nu\right)^2\nonumber\\
&=&\frac{1}{\eps^{2/3}}\left(\Phi_\eps\omega+(1-\Phi_\eps)(\omega+o_{L^\infty(D_0\cap
  D_1)}(\eps^{\beta(N-1/2)}))\right)^2\nonumber\\
&=&\frac{1}{\eps^{2/3}}\left(\omega+o_{L^\infty(D_0\cap
  D_1)}(\eps^{\beta(N-1/2)})\mathbf{1}_{D_0\cap D_1}\right)^2\nonumber\\
&=&\frac{\omega^2}{\eps^{2/3}}+o_{L^\infty(D_0\cap
  D_1)}(\eps^{\beta N-2/3)})\mathbf{1}_{D_0\cap D_1}\nonumber\\
&=&\frac{\omega_0^2}{\eps^{2/3}}+\frac{1}{\eps^{2/3}}\sum_{m=1}^{2M}\eps^{2m}\sum_{\substack{m_1+m_2=m\\0\leq
    m_1,m_2\leq M}}\omega_{m_1}\omega_{m_2}+o_{L^\infty(D_0\cap
  D_1)}(\eps^{\beta N-2/3)})\mathbf{1}_{D_0\cap D_1}\nonumber\\
&=&\frac{\Gamma_2z}{2\alpha_1\Gamma_{12}\eps^{2/3}}+O_{L^\infty(D_0)}(\eps^{4/3-2\beta})+o_{L^\infty(D_0\cap
  D_1)}(\eps^{\beta N-2/3)})\mathbf{1}_{D_0\cap D_1},
\ee
where for the last equality, we have used (\ref{asymp-omegan}) to
infer that for $m\geq 1$,
$m_1+m_2=m$ and $x\in D_0$, $\omega_{m_1}\omega_{m_2}\lesssim
|z|^{1-3m}\leq\eps^{\beta-3\beta m}$, and that $2m+\beta-3\beta m\geq
2-2\beta$. The same kind of calculation yields, still for $x\in D_0$,
\be\label{rho2d0}
\rho_2^2&=&\frac{1}{\eps^{2/3}}\left(\Phi_\eps\tau+(1-\Phi_\eps)(\tau+o_{L^\infty(D_0\cap
  D_1)}(\eps^{\beta N}))\right)^2\nonumber\\
&=&\frac{\tau_0^2}{\eps^{2/3}}+\frac{1}{\eps^{2/3}}\sum_{m=1}^{2M}\eps^{2m}\sum_{\substack{m_1+m_2=m\\0\leq
    m_1,m_2\leq M}}\tau_{m_1}\tau_{m_2}+o_{L^\infty(D_0\cap
  D_1)}(\eps^{\beta N-2/3)})\mathbf{1}_{D_0\cap D_1}\nonumber\\
&=&\frac{R_2^2-R_1^2}{2\alpha_2\eps^{2/3}}+\frac{\Gamma_1z}{2\alpha_2\Gamma_{12}\eps^{2/3}}+O_{L^\infty(D_0)}(\eps^{4/3-2\beta})+o_{L^\infty(D_0\cap
  D_1)}(\eps^{\beta N-2/3)})\mathbf{1}_{D_0\cap D_1}.
\ee
Next, we deduce from (\ref{nu00}), (\ref{asymp-nun}),
and (\ref{asymp-nu0}) that for $x\in D_1\string\ (D_0\cup D_2)$,
\be\label{rho1d1}
\rho_1^2&=&\nu^2=\nu_0^2+(\nu^2-\nu_0^2)=\nu_0^2+\sum_{n=1}^{2N}\eps^{2n/3}\sum_{\substack{n_1+n_2=n\\0\leq
    n_1,n_2\leq
    N}}\nu_{n_1}\nu_{n_2}\nonumber\\
&=&\frac{R_1^{2/3}\Gamma_2^{2/3}}{2\alpha_1\Gamma_{12}}\gamma_0\left(\widetilde{y}_1\right)^2+O_{L^\infty(D_1\string\ (D_0\cup D_2))}(\eps^{2/3}),
\ee
and using (\ref{I-1}), (\ref{II0}), (\ref{asymp-lamn}) and
(\ref{asymp-lam0}), we get (again for $x\in D_1\string\ (D_0\cup D_2)$)
\be\label{rho2d1}
\rho_2^2&=&\lambda=\frac{R_2^2-R_1^2}{2\alpha_2\eps^{2/3}}+\frac{y_1}{2\alpha_2}-\frac{\alpha_0
  R_1^{2/3}\Gamma_2^{2/3}}{2\alpha_1\alpha_2\Gamma_{12}}\gamma_0\left(\widetilde{y}_1\right)^2
+O_{L^\infty(D_1\string\ (D_0\cup D_2))}(\eps^{2/3}),
\ee
where 
$$\widetilde{y}_1=\frac{\Gamma_2^{1/3}y_1}{R_1^{2/3}}.$$
For $x\in D_2$, we use Lemma \ref{lem-compar-D1D2} to
obtain, for $\alpha>0$ arbitrarily large,
\be\label{rho1d2}
\rho_1^2&=&\chi_\eps^2\nu^2=o_{L^\infty(D_1\cap D_2)}(\eps^\alpha)\mathbf{1}_{L^\infty(D_1\cap D_2)},
\ee
as well as, using also (\ref{mun}), (\ref{asympnu0+})
and Proposition \ref{proposition-mun},
\be\label{rho2d2}
\rho_2^2&=&\frac{1}{\eps^{2/3}}\left(\Psi_\eps\eps^{1/3}\mu+(1-\Psi_\eps)(\eps^{1/3}\mu+o_{L^\infty(D_1\cap D_2)}(\eps^{\beta(N+1)}))\right)^2\nonumber\\
&=&\frac{1}{\eps^{2/3}}\left(\eps^{1/3}\mu+o_{L^\infty(D_1\cap D_2)}(\eps^{\beta(N+1)})\mathbf{1}_{L^\infty(D_1\cap D_2)}\right)^2\nonumber\\
&=&\mu_0^2+\sum_{n=1}^{2L}\eps^{2n/3}\sum_{\substack{n_1+n_2=n\\0\leq
    n_1,n_2\leq L}}\mu_{n_1}\mu_{n_2}+o_{L^\infty(D_1\cap D_2)}(\eps^{\beta(N+1)-2/3})\mathbf{1}_{L^\infty(D_1\cap D_2)}\nonumber\\
&=&\frac{R_2^{2/3}}{2\alpha_2}\gamma_0(\widetilde{y}_2)^2+O_{L^\infty(D_2)}(\eps^{2/3})+o_{L^\infty(D_1\cap D_2)}(\eps^{\beta(N+1)-2/3})\mathbf{1}_{L^\infty(D_1\cap D_2)},
\ee
where
$$\widetilde{y}_2=\frac{y_2}{R_2^{2/3}}.$$
From (\ref{rho1d0}), (\ref{rho2d0}), (\ref{rho1d1}), (\ref{rho2d1}),
(\ref{rho1d2}), (\ref{rho2d2}) and the definitions of $p_\eps$ and
$q_\eps$, we deduce the following expressions of $p_\eps$ and
$q_\eps$, depending on whether $x\in D_0$, $x\in D_1\string\(D_0\cup
D_2)$ or $x\in D_2$. For each of this cases, we also calculate
$r_\eps^2$ and $-\Delta_\eps=p_\eps q_\eps-r_\eps^2$, a quantity which
will play a key role in the sequel. A large integer $N$ been fixed, We
assume that $\beta\in (0,2/3)$ satisfies $\beta
N-2/3>4/3-2\beta$ and $\beta(N+1)-2/3>2/3$ (which are equivalent to
$\beta >2/(N+2)$ if $N$ is large). For $x\in D_0$, we obtain
\be\label{ped0}
p_\eps(x)
&=&\frac{2\Gamma_2y_1}{\Gamma_{12}}+O_{L^\infty(D_0)}(\eps^{4/3-2\beta}),
\ee
\be\label{qed0}
q_\eps(x)&=\frac{2(R_2^2-R_1^2)}{\eps^{2/3}}+\frac{2\Gamma_1
  y_1}{\Gamma_{12}}+O_{L^\infty(D_0)}(\eps^{4/3-2\beta}),
\ee
\be\label{red0}
r_\eps(x)^2=\frac{4\alpha_0^2\Gamma_2
  y_1}{\alpha_1\alpha_2\Gamma_{12}\eps^{2/3}}\left(R_2^2-R_1^2+\frac{\Gamma_1}{\Gamma_{12}}z\right)+O_{L^\infty(D_0)}(\eps^{2/3-2\beta}),
\ee
\be\label{ded0}
-\Delta_\eps(x)&=&\frac{4\Gamma_2y_1}{\eps^{2/3}}\left(R_2^2-R_1^2+\frac{\Gamma_1}{\Gamma_{12}}z\right)+O_{L^\infty(D_0)}(\eps^{2/3-2\beta}).
\ee
For $x\in D_1\string\ (D_0\cup D_2)$, 
\be\label{ped1}
p_\eps(x)
&=&R_1^{2/3}\Gamma_2^{2/3}\widetilde{W}_0(\widetilde{y}_1)+O_{L^\infty(D_1\string\ (D_0\cup D_2))}(\eps^{2/3}),
\ee
where 
$$\widetilde{W}_0(y)=\left(1+\frac{2}{\Gamma_{12}}\right)\gamma_0(y)^2-y,$$
\be\label{qed1}
q_\eps(x)&=2y_2-\frac{2\alpha_0R_1^{2/3}\Gamma_{2}^{2/3}}{\alpha_1\Gamma_{12}}\gamma_0(\widetilde{y}_1)^2+O_{L^\infty(D_1\string\ (D_0\cup D_2))}(\eps^{2/3}),
\ee
\be\label{red1}
r_\eps(x)^2=\frac{4\alpha_0^2R_1^{2/3}\Gamma_{2}^{2/3}}{\alpha_1\alpha_2\Gamma_{12}}\gamma_0(\widetilde{y}_1)^2 \left(y_2-\frac{\alpha_0R_1^{2/3}\Gamma_{2}^{2/3}}{\alpha_1\Gamma_{12}}\gamma_0(\widetilde{y}_1)^2\right)+O_{L^\infty(D_1\string\ (D_0\cup D_2))}(1),
\ee
\be\label{ded1}
-\Delta_\eps&=&2R_1^{2/3}\Gamma_{2}^{2/3}\left(y_2-\frac{\alpha_0R_1^{2/3}\Gamma_{2}^{2/3}}{\alpha_1\Gamma_{12}}\gamma_0(\widetilde{y}_1)^2\right)W_0(\widetilde{y}_1)+O_{L^\infty(D_1\string\ (D_0\cup D_2))}(1).
\ee
Finally, for $x\in D_2$,
\be\label{ped2}
p_\eps(x)
&=&-\Gamma_2y_1+\frac{\alpha_0}{\alpha_2}R_2^{2/3}V_0(\widetilde{y}_2)+O_{L^\infty(D_2)}(\eps^{2/3}),
\ee
where 
$$V_0(y)=\gamma_0(y)^2-y,$$
\be\label{qed2}
q_\eps(x)&=R_2^{2/3}W_0(\widetilde{y}_2)+O_{L^\infty(D_2)}(\eps^{2/3}),
\ee
\be\label{red2}
r_\eps(x)^2=o_{L^\infty(D_2)}(\eps^{\alpha}),
\ee
for $\alpha$ arbitrarily large, and
\be\label{ded2}
-\Delta_\eps&=&\left(-\Gamma_2y_1+\frac{\alpha_0}{\alpha_2}R_2^{2/3}V_0(\widetilde{y}_2)\right)R_2^{2/3}W_0(\widetilde{y}_2)+O_{L^\infty(D_2)}(1).
\ee
Then, (\ref{ped0}), (\ref{ped1}) and (\ref{ped2}) will provide us
upper and lower bounds on $p_\eps$. For this purpose, since
the function $\widetilde{W}_0$ appears in (\ref{ped1}), we first
prove a lemma which gives informations about the size of this function.

\begin{lem}\label{W0} 
For $y\in \R$,
\begin{eqnarray}
W_0(y)\leq \widetilde{W}_0(y)\lesssim W_0(y),
\end{eqnarray}
and
\begin{eqnarray}\label{encW0}
\max(1,|y|)\lesssim W_0(y)\lesssim \max(1,|y|).
\end{eqnarray}
\end{lem}

\begin{Proof}
We write
$$\widetilde{W}_0(y)=W_0(y)+2\left(\frac{1}{\Gamma_{12}}-1\right)\gamma(y)^2,$$
where $1/\Gamma_{12}-1>0$, which directly provides the lower bound on
$\widetilde{W}_0$. Moreover, the analysis of the continuous functions $\gamma_0$ and
$W_0$ which was done in \cite{GP} ensures that $W_0(y)>0$ for every
$y\in\R$, $W_0(y)\underset{y\to +\infty}{\sim}2y$,
$W_0(y)\underset{y\to -\infty}{\sim}-y$, $\gamma_0(y)^2\underset{y\to
  +\infty}{\sim}y$ and $\gamma(y)^2\underset{y\to
  -\infty}{\longrightarrow}0$. We deduce (\ref{encW0}) and the
existence of $C_0>0$
such that $\gamma_0^2/W_0\leq C_0$. Then, we obtain the upper bound
$$\widetilde{W}_0\leq\left(1+2\left(\frac{1}{\Gamma_{12}}-1\right)C_0\right)W_0.$$
\end{Proof}
Then, we get lower and upper bounds on $p_\eps$ as stated in the lemma below.

\begin{lem}\label{ppos}
For $x\in \R^d$ and $\eps>0$ sufficiently small,
$$\max(1,|y_1|)\lesssim p_\eps(x)\lesssim \max(1,|y_1|).$$
\end{lem}

\begin{Proof}
On $D_0$, the two estimates directly follow from (\ref{ped0}), since
for $x\in D_0$, $y_1\geq
\eps^{\beta-2/3}\to+\infty$ as $\eps\to
  0$, whereas $\eps^{4/3-2\beta}\to 0$. On $D_1\string\ (D_0\cup D_2)$, they are
consequences from (\ref{ped1}) and Lemma \ref{W0}. On $D_2$, we deduce
them from (\ref{ped2}). Indeed, we know from the asymptotic expansion of
$\gamma_0$ (\ref{DAHM}) that $V_0(y)\underset{y\to
  +\infty}{=}O(y^{-2})$, and $V_0(y)\underset{y\to
  -\infty}{\sim}-y\to +\infty$, therefore $V_0$ is bounded from
below. Then, we have on the one side, for $\eps$ sufficiently small,
\begin{eqnarray}\label{pbelowcas2}
p_\eps(x) &\geq
&\frac{\alpha_0}{\alpha_2}R_2^{2/3}\underset{y\in\R}{\inf}V_0(y)+\Gamma_2|y_1|-1\gtrsim\max(1,|y_1|)
\gtrsim \eps^{\beta-2/3}\underset{\eps\to
  0}{\longrightarrow}+\infty.\ \ \ \ \ \ \ 
\end{eqnarray}
On the other side, the properties of $\gamma_0$ stated in
Proposition \ref{proposition-Painleve} imply 
$$\forall y\in\R,\quad V_0(y)\lesssim\max(1,-y).$$
Thus, using also Lemma \ref{W0} and the inequalities $y_1<0$ and $y_1<y_2$, we get
\begin{eqnarray}\label{pupcas2}
p_\eps(x)\lesssim\max(1,-y_2)+\max(1,-y_1)\lesssim\max(1,-y_1)=\max(1,|y_1|).
\end{eqnarray}
\end{Proof}
As for $q_\eps$, we infer similarly the next lemma from (\ref{qed0}),
(\ref{qed1}) and (\ref{qed2}).

\begin{lem}\label{qpos}
For $x\in\R^d$ and $\eps>0$ sufficiently small,
$$\max(1,|y_2|)\lesssim q_\eps(x)\lesssim \max(1,|y_2|).$$
\end{lem}

\begin{Proof}
In order to prove the two inequalities for $x\in D_0$, we rewrite
(\ref{qed0}) as 
\be
q_\eps(x)&=\frac{2}{\eps^{2/3}}\left(R_2^2-R_1^2+\frac{\Gamma_1
  }{\Gamma_{12}}z\right)+O_{L^\infty(D_0)}(\eps^{4/3-2\beta}).
\ee
As $x$ describes $D_0$, $z$ describes the interval
$[\eps^\beta,R_1^2]\subset[0,R_1^2]$. On this interval, $G(z)=
R_2^2-R_1^2+\frac{\Gamma_1}{\Gamma_{12}}z$ reaches its extrema at
$z=0$ and at $z=R_1^2$. Since $G(0)=R_2^2-R_1^2>0$ and
$G(R_1^2)=R_2^2-\frac{\alpha_0\Gamma_2}{\alpha_1\Gamma_{12}}R_1^2>0$
(thanks to (\ref{cond-disk})),
and because $-2/3<0<4/3-2\beta$, there exists a constant $c>1$ such
that for every $x\in D_0$,
$$\frac{1}{c}\eps^{-2/3}\leq q_\eps(x)\leq c\eps^{-2/3}.$$
The inequality follows for $x\in D_0$, since
$(R_2^2-R_1^2)\eps^{-2/3}\leq y_2\leq R_2^2\eps^{-2/3}$ on $D_0$.
On $D_1\string\ (D_1\cup D_2)$, the inequalities clearly follow from
(\ref{qed1}), since on that set, $y_2\gtrsim\eps^{-2/3}$. Finally, on
$D_2$, they are consequences of (\ref{qed2}) and Lemma \ref{W0}. 
\end{Proof}

We are now ready to prove positivity of the operator $A_\eps$

\begin{theo}\label{form-domain}
$A_\eps$ defines a positive self-adjoint operator on $L^2(\R^d)^2$, with form domain $H_w^1(\R^d)^2$, where
$$H_w^1(\R^d)=\left\{P\in H^1(\R^d)\big|
  \max(1,\min(|y_1|,|y_2|))^{1/2}P\in L^2(\R^d)\right\}.$$ 
Moreover, there exists $C>0$
  such that for every $(P,Q)$ in the domain of $A_\eps$,
$$\left<A_\eps \left[\begin{array}{c}
P\\Q
\end{array}\right], \left[\begin{array}{c}
P\\Q
\end{array}\right]\right>\geq \eps^{4/3}\int_{\R^d}\left(|\nabla
  P|^2+|\nabla Q|^2\right)dx+C\int_{\R^d}\max(1,\min(|y_1|,|y_2|))(|P|^2+|Q|^2)dx.$$
\end{theo}

\begin{Proof}
For $P,Q\in \mathcal{C}_c^\infty(\R^d)$, we have
\begin{eqnarray}\label{debproof}
\left<A_\eps \left[\begin{array}{c}
P\\Q
\end{array}\right], \left[\begin{array}{c}
P\\Q
\end{array}\right]\right>=\int_{\R^d}\left(\eps^{4/3}|\nabla
P|^2+\eps^{4/3}|\nabla Q|^2+p_\eps P^2+q_\eps Q^2+2r_\eps
PQ\right)dx.
\end{eqnarray}
Taking into account the positivity of $p_\eps$ and $q_\eps$ shown in
Lemmata \ref{ppos} and \ref{qpos}, 
$$p_\eps P^2+q_\eps Q^2+2r_\eps
PQ=\frac{1}{p_\eps}\left(p_\eps P+r_\eps Q\right)^2+\frac{p_\eps
  q_\eps-r_\eps^2}{p_\eps}Q^2\geq \frac{-\Delta_\eps}{p_\eps}Q^2,$$
where
$$\Delta_\eps=r_\eps^2-p_\eps q_\eps.$$
Symmetrically,
$$p_\eps P^2+q_\eps Q^2+2r_\eps PQ\geq\frac{-\Delta_\eps}{q_\eps}P^2,$$
and thus
\begin{eqnarray}\label{minquad}
p_\eps P^2+q_\eps Q^2+2r_\eps
PQ\geq
\frac{1}{2}\min\left(\frac{-\Delta_\eps}{p_\eps},\frac{-\Delta_\eps}{q_\eps}\right)(P^2+Q^2).
\end{eqnarray}
We shall see next that there exists  $c>0$ such that for every $x\in \R^d$,
\be\label{min-deltae}
-\Delta_\eps(x)\geq c p_\eps(x) q_\eps(x),
\ee
which is equivalent to
\be
-\Delta_\eps(x)\geq c \max(1,|y_1|)\max(1,|y_2|)
\ee
thanks to Lemmata \ref{ppos} and \ref{qpos} (up to a change of the
constant $c>0$), and which implies
\be\label{minmax}
\min\left(\frac{-\Delta_\eps}{p_\eps},\frac{-\Delta_\eps}{q_\eps}\right)\geq
c\min\left(p_\eps(x),q_\eps(x)\right)\gtrsim
\min\left(\max(1,|y_1|),\max(1,|y_2|)\right)=\max(1,\min(|y_1|,|y_2|)).\ 
\ee
For $x\in D_0$, (\ref{min-deltae}) comes from (\ref{ded0}), because
for such values of $x$, we have $y_1\gtrsim\eps^{\beta-2/3}\geq 1$ for
$\eps$ sufficiently small (and therefore $y_1=\max(1,|y_1|)$), we also
have $R_2^2\eps^{-2/3}\geq y_2\geq (R_2^2-R_1^2)\eps^{-2/3}$, and
therefore $\eps^{-2/3}\gtrsim \max(1,|y_2|)$, and finally, the remark
we have done to bound $q_\eps$ from below in the proof of Lemma
\ref{qpos} implies that $R_2^2-R_1^2+\frac{\Gamma_1}{\Gamma_{12}}z$ is bounded from below by a positive constant as
$x\in D_0$. For $x\in D_1\string\ (D_0\cup D_2)$ (\ref{min-deltae})
follows from (\ref{ded1}) and Lemma \ref{W0}, since
$y_2=\max(1,|y_2|)$ on that domain, and $\gamma_0(\widetilde{y_1})^2)\lesssim \eps^{\beta-2/3}$. For $x\in D_2$, note that
$W_0(\widetilde{y}_2)\gtrsim \max(1,|y_2|)$ thanks to Lemma
\ref{W0}. Then, using (\ref{ded2}) and the same arguments as to
obtain (\ref{pbelowcas2}), we complete the proof.
\end{Proof}

We deduce classicaly from Theorem \ref{form-domain} the following corollary.
\begin{cor}\label{est-res}
$A_\eps$ is invertible, and
  $$\|A_\eps^{-1}\|_{\mathcal{L}(L^2(\R^d)^2,H^1_w(\R^d)^2)}\lesssim\eps^{-4/3},$$
where $H^1_w(\R^d)^2$ is endowed by the
  norm 
$$\left\|(P,Q)\right\|_{H^1_w(\R^d)^2}=\left(\int_{\R^d}\left(|\nabla
  P|^2+|\nabla Q|^2\right)dx+\int_{\R^d}\max(1,\min(|y_1|,|y_2|))(|P|^2+|Q|^2)dx\right)^{1/2}.$$
\end{cor}

\begin{rem}\label{embH1w}
Note that the set $H^1_w(\R^d)^2$ does not depend on $\eps$, even though it's norm does. However, our choice of the $H^1_w(\R^d)^2$-norm ensures that the norm of the embedding of $H^1_w(\R^d)^2$ into $H^1(\R^d)^2$ is uniformly bounded in $\eps$.
\end{rem}

\subsection{The fixed point argument}\label{fp}
Let $1/3>\delta>0$, and $N$ a large integer. We fix $\beta\in (0,2/3)$ such that
$(2-3\beta)N/3<\delta$, and then $L$ and $M$ large enough such that 
$(2-3\beta)(M+1/2)-2N/3-2/3>-2$,
  $(2-3\beta)(M+1/2)-2N/3-2/3-\beta/2>-2$ and $2(L-N)/3>-2$, in such a
  way that (\ref{est-fe01}) and (\ref{est-fe02}) imply
\be\label{est-fe}
\left\|f_\eps^0\right\|_{L^2(\R^d)^2} &\lesssim & \eps^{-2}.
\ee
We are going to apply the Picard fix-point theorem to the map
$$\begin{array}{crcl}\Theta_\eps\ :&H^1_w(\R^d)^2&\longrightarrow
  &H^1_w(\R^d)^2\\ &(P,Q)&\longrightarrow & A_\eps^{-1}f_\eps^0+A_\eps^{-1}f_\eps^2(P,Q)+A_\eps^{-1}f_\eps^3(P,Q),\end{array}$$
in the ball $\mathcal{B}_R$ of $H^1_w(\R^d)^2$ centered at the origin, with radius
$R=2\left\|A_\eps^{-1}f_\eps^0\right\|_{H^1_w(\R^d)^2}$. Note that
it follows from Corollary \ref{est-res} and (\ref{est-fe}) that
\be\label{est-R}
R\lesssim \eps^{-10/3}.
\ee
From (\ref{rho1d0}), (\ref{rho2d0}), (\ref{rho1d1}), (\ref{rho2d1}),
(\ref{rho1d2}), (\ref{rho2d2}), it follows that for $x\in \R^d$,
$$|\rho_1|\lesssim\eps^{-1/3}\quad {\rm and}\quad |\rho_2|\lesssim\eps^{-1/3}.$$
Thus, the Sobolev embedding $H^1_w(\R^d)\subset
H^1(\R^d)\subset L^4(\R^d)$ ($d\leq 3$) implies that for every $(P,Q)\in
H^1_w(\R^d)^2$, we have $f_\eps^2\in L^2(\R^d)^2$, and
\be\label{est-fe2}
\|f_\eps^2(P,Q)\|_{L^2(\R^d)^2}\lesssim\eps^{2N/3+1/3}\|(P,Q)\|_{H^1_w(\R^d)^2}^2.
\ee
Then, Corollary \ref{est-res} yields
\be\label{sob2}
\left\|A_\eps^{-1}f_\eps^2(P,Q)\right\|_{H^1_w(\R^d)^2}\lesssim
\eps^{2N/3-1}\|(P,Q)\|_{H^1_w(\R^d)^2}^2.
\ee
Similarly, thanks to the Sobolev embedding $H^1_w(\R^d)\subset
H^1(\R^d)\subset L^6(\R^d)$ ($d\leq 3$), we get, for $(P,Q)\in
H^1_w(\R^d)^2$,
\be\label{est-fe3}
\left\|f_\eps^3(P,Q)\right\|_{H^1_w(\R^d)^2}\lesssim
\eps^{4N/3+4/3}\|(P,Q)\|_{H^1_w(\R^d)^2}^3
\ee
and
\be\label{sob3}
\left\|A_\eps^{-1}f_\eps^3(P,Q)\right\|_{H^1_w(\R^d)^2}\lesssim
\eps^{4N/3}\|(P,Q)\|_{H^1_w(\R^d)^2}^3.
\ee
From (\ref{sob2}) and (\ref{sob3}) we deduce that if
$(P,Q)\in\mathcal{B}_R$, for some positive constants $C_2$ and $C_3$, 
$$\left\|\Theta_\eps(P,Q)\right\|_{H^1_w(\R^d)^2}\leq\frac{R}{2}+C_2\eps^{2N/3-13/3}R+C_3\eps^{4N/3-20/3}R.
$$
Therefore, if $N\geq 7$ and $\eps$ is sufficiently small,
$\mathcal{B}_R$ is stable by $\Theta_\eps$. Similar arguments prove
that $\Theta_\eps$ is a contraction on that ball. As a result, $\Theta_\eps$  has a unique fixed point in $\mathcal{B}_R$.

\subsection{Positivity of $\eta_1$ and $\eta_2$.}\label{positivity}
This section is devoted to the proof of the positivity of the solution
$(\eta_1,\eta_2)$ to the system (\ref{sys}) given by
(\ref{ansatz-rho})-(\ref{ansatz-rho2})-(\ref{def-rho})-(\ref{def-rho2}), which has just been constructed in
the sections above. We proceed in three steps. First, we prove that
for $j=1,2$, $\rho_j$ (given by (\ref{def-rho}) or (\ref{def-rho2})) is bounded from below
by a positive constant on the set $S_j=\left\{x\in \R^d, |x|^2\leq
R_j^2+\eps^{2/3}\right\}$, provided $\eps$ is sufficiently
small. Second, we prove $L^\infty$ estimates on $P$ and $Q$, which
ensure that $\eta_1$ and $\eta_2$ are positive on $S_1$ and $S_2$
respectively. Finally, we prove positivity of $(\eta_1,\eta_2)$ on
$\R^d$ thanks to the maximum principle.

\paragraph{1$^{\rm st}$ step.} For some integers $N,M,L\geq 1$, let $\omega$, $\tau$, $\nu$,
$\lambda$ and $\mu$ be the functions given by (\ref{troncat}). Then, we decompose the functions $\rho_1$ and $\rho_2$ given by (\ref{def-rho})-(\ref{def-rho2})  as
\be\label{dec-rho1}
\rho_1&=& \eps^{-1/3}\omega\mathbf{1}_{D_0\string\ D_1}
+\left(\Phi_\eps\eps^{-1/3}\omega+\chi_\eps\nu\right)\mathbf{1}_{D_0\cap D_1}+\nu\mathbf{1}_{D_1\string\ D_0}
\ee
and
\be\label{dec-rho2}
\rho_2 &=&\eps^{-1/3}\tau\mathbf{1}_{D_0\string\ D_1}
+\left(\Phi_\eps\eps^{-1/3}\tau+\chi_\eps\lambda^{1/2}\right)\mathbf{1}_{D_0\cap D_1}+\lambda^{1/2}\mathbf{1}_{D_1\string\ (D_0\cup D_2)}\nonumber\\
&&+
\left(\chi_\eps\lambda^{1/2}+\Psi_\eps\mu\right)\mathbf{1}_{D_1\cap D_2}+\mu\mathbf{1}_{D_2\string\ D_1},
\ee
and we are going to bound from below $\omega$, $\tau$, $\nu$,
$\lambda$ and $\mu$ separately on the different sets appearing in the indicatrix functions above. According to Remark \ref{23}, $\omega$
and $\tau$ satisfy
\be
\omega=\omega_0+O_{L^\infty(D_0)}(\eps^{2-5\beta/2})\quad {\rm
  and}\quad \tau=\tau_0+O_{L^\infty(D_0)}(\eps^{2-2\beta}).
\ee
Moreover, thanks to the explicit expressions of $\omega_0$ and
$\tau_0$ (\ref{eq-omega0}) and (\ref{eq-tau0}), we deduce that for
$x\in D_0$, 
\be
\omega_0\geq\left(\frac{\Gamma_2}{2\alpha_1\Gamma_{12}}\right)^{1/2}\eps^{\beta/2}\quad
      {\rm and} \quad \tau_0\geq \left(\frac{R_2^2-R_1^2}{2\alpha_2}\right)^{1/2}+O(\eps^\beta).
\ee
Since $\beta<2/3$, we have $2-5\beta/2>\beta/2$ and $2-2\beta>\beta$,
so we conclude that for $x\in D_0$, 
\be\label{457}
\omega\geq\left(\frac{\Gamma_2}{2\alpha_1\Gamma_{12}}\right)^{1/2}\eps^{\beta/2}+O(\eps^{2-5\beta/2})\quad
      {\rm and} \quad \tau\geq \left(\frac{R_2^2-R_1^2}{2\alpha_2}\right)^{1/2}+O(\eps^\beta).
\ee
We have already seen in the proof of Lemma \ref{minor-lam} that for
$x\in D_1$,
\be\label{458}
\lambda\geq
\frac{R_2^2-R_1^2}{2\alpha_2\eps^{2/3}}+O(\eps^{\beta-2/3}).
\ee
Using similar arguments, we infer from Proposition \ref{summary} that
\be
\nu=\nu_0+O_{L^\infty(D_1)}(\eps^{2/3}).
\ee
Then, (\ref{nu00}), the fact that $\gamma_0$
is an increasing function and Proposition \ref{proposition-Painleve} imply
\be\label{460}
\nu\geq\frac{(R_1\Gamma_2)^{1/3}}{(2\alpha_1\Gamma_{12})^{1/2}}\gamma_0\left(\frac{\Gamma_2^{1/3}}{R_1^{2/3}}\eps^{\beta-2/3}\right)+O(\eps^{2/3})=\frac{\Gamma_2^{1/2}}{(2\alpha_1\Gamma_{12})^{1/2}}\eps^{\beta/2-1/3}+O(\eps^{2/3})\quad
            {\rm for\ } x\in D_0\cap D_1,
\ee
whereas
\be\label{461}
\nu\geq
\frac{(R_1\Gamma_2)^{1/3}}{(2\alpha_1\Gamma_{12})^{1/2}}\gamma_0\left(-\frac{\Gamma_2^{1/3}}{R_1^{2/3}}\right)+O(\eps^{2/3})\quad
{\rm for\ } x\in D_1\cap S_1.
\ee
From Proposition \ref{proposition-mun}, we get in the same way
\be
\mu=\mu_0+O_{L^\infty(D_2)}(\eps^{2/3}),
\ee
which implies thanks to (\ref{mun}) for $n=0$
\be\label{463}
\mu\geq \frac{R_2^{1/3}}{(2\alpha_2)^{1/2}}\gamma_0\left(-\frac{1}{R_2^{2/3}}\right)+O(\eps^{2/3})\quad {\rm for}\ x\in D_2\cap S_2.
\ee
Combining (\ref{dec-rho1}), (\ref{457}), (\ref{460}) and
(\ref{461}), we deduce that
\be\label{min-rho1}
\rho_1&\geq &
\left(\eps^{-1/3}\left(\frac{\Gamma_2}{2\alpha_1\Gamma_{12}}\right)^{1/2}\eps^{\beta/2}+O(\eps^{5/3-5\beta/2})\right)\mathbf{1}_{D_0\string\ D_1}\nonumber\\
&&+\left(\Phi_\eps\eps^{-1/3}\left(\frac{\Gamma_2}{2\alpha_1\Gamma_{12}}\right)^{1/2}\eps^{\beta/2}+\chi_\eps\left(\frac{\Gamma_2}{2\alpha_1\Gamma_{12}}\right)^{1/2}\eps^{\beta/2-1/3}+O(\eps^{5/3-5\beta/2})+O(\eps^{2/3})\right)\mathbf{1}_{D_0\cap
  D_1}\nonumber\\
&&+\left(\frac{(R_1\Gamma_2)^{1/3}}{(2\alpha_1\Gamma_{12})^{1/2}}\gamma_0\left(-\Gamma_2^{1/3}/R_1^{2/3}\right)+O(\eps^{2/3})\right)\mathbf{1}_{S_1\string\ D_0}.\nonumber\\
&\geq
&\left(\frac{\Gamma_2}{2\alpha_1\Gamma_{12}}\right)^{1/2}\eps^{\beta/2-1/3}\mathbf{1}_{D_0}+\frac{(R_1\Gamma_2)^{1/3}}{(2\alpha_1\Gamma_{12})^{1/2}}\gamma_0\left(-\Gamma_2^{1/3}/R_1^{2/3}\right)\mathbf{1}_{S_1\string\ D_0}\nonumber\\
&&+O(\eps^{5/3-5\beta/2})\mathbf{1}_{D_0}+O(\eps^{2/3})\mathbf{1}_{S_1\cap D_1}.\nonumber\\
&\geq &c_1\mathbf{1}_{S_1},
\ee
if $\eps\leq 1 $ is sufficiently small, where $c_1=\frac{1}{2}\min\left[\left(\frac{\Gamma_2}{2\alpha_1\Gamma_{12}}\right)^{1/2},\frac{(R_1\Gamma_2)^{1/3}}{(2\alpha_1\Gamma_{12})^{1/2}}\gamma_0\left(-\frac{\Gamma_2^{1/3}}{R_1^{2/3}}\right)\right]$ .
On the other side, using (\ref{dec-rho2}), (\ref{457}), (\ref{458}) and (\ref{463}), we have
\be\label{min-rho2}
\rho_2&\geq
&\left(\left(\frac{R_2^2-R_1^2}{2\alpha_2}\right)^{1/2}\eps^{-1/3}+O(\eps^{\beta-1/3})\right)\mathbf{1}_{(D_0\cup
  D_1)\string\ D_2}\nonumber\\
&&+\left(\min\left[  \left(\frac{R_2^2-R_1^2}{2\alpha_2}\right)^{1/2}\eps^{-1/3},\frac{R_2^{1/3}}{(2\alpha_2)^{1/2}}\gamma_0\left(\frac{-1}{R_2^{2/3}}\right)\right]+O(\eps^{\beta-1/3})\right)\mathbf{1}_{D_1\cap
  D_2}\nonumber\\
&&+\left(\frac{R_2^{1/3}}{(2\alpha_2)^{1/2}}\gamma_0\left(\frac{-1}{R_2^{2/3}}\right)+O(\eps^{2/3})\right)\mathbf{1}_{S_2\string\ D_1}\nonumber\\
&\geq & c_2\mathbf{1}_{S_2}
\ee
for $\eps\leq 1$ sufficiently small, where $c_2=\frac{1}{2}\min\left[
  \left(\frac{R_2^2-R_1^2}{2\alpha_2}\right)^{1/2},\frac{R_2^{1/3}}{(2\alpha_2)^{1/2}}\gamma_0\left(\frac{-1}{R_2^{2/3}}\right)\right]$.

\paragraph{2$^{\rm nd}$ step.} Let $N$ be a large integer, $R>0$
as in  Section \ref{fp} and $(P,Q)\in (H_w^1)^2$ the unique fixed
point of the map $\Theta_\eps$ constructed in that section. In order
to control the $L^\infty$ norm of $(P,Q)$, we will use the continuity
of the embedding of $H^2(\R^d)$ into $L^\infty(\R^d)$ (remember that
$d\leq 3$). Because of Remark \ref{embH1w} and (\ref{est-R}), we know that the $H^1$ norms of $P$ and $Q$ are controled by
$$\|(P,Q)\|_{(H^1)^2}\lesssim \eps^{-10/3},$$
such that in order to control the $H^2$ norms of $P$ and $Q$, we only
need to control the $L^2$ norms of $\Delta P$ and $\Delta Q$. For this purpose, let us introduce a $C^\infty$ function $\theta$ on $\R^d$, which is radial, positive, supported in $\{x\in\R^d, |x|\leq 2\}$ and such that $\theta(x)\equiv 1$ for $|x|\leq 1$. We also define, for integers $n\geq 1$, $\theta_n(x)=\theta(x/n)$. After integrations by parts, the $(L^2)^2$ scalar product of (\ref{eq-reste}) with $(\Delta P,\Delta Q)\theta_n$ yields
\be\label{est-Delta}
\lefteqn{\eps^{4/3}\int_{\R^d}(|\Delta P|^2+|\Delta Q|^2)\theta_n +\int_{\R^d}(p_\eps|\nabla P|^2+q_\eps|\nabla Q|^2+2r_\eps\nabla P\cdot \nabla Q)\theta_n}\nonumber\\
&=&-\int_{\R^d}(P\nabla p_\eps\cdot\nabla P+Q\nabla q_\eps\cdot\nabla Q)\theta_n-\int_{\R^d}(p_\eps P\nabla P+q_\eps Q \nabla Q)\cdot\nabla\theta_n\nonumber\\
&&-\int_{\R^d}\nabla r_\eps \cdot  \nabla (PQ)\theta_n
-\int_{\R^d}r_\eps  \nabla (PQ)\cdot \nabla \theta_n-\int_{\R^d}f_\eps\cdot (\Delta P,\Delta Q)\theta_n\nonumber\\
&=&\frac{1}{2}\int_{\R^d}(\Delta p_\eps P^2+\Delta q_\eps Q^2+2\Delta r_\eps PQ)\theta_n
+\frac{1}{2}\int_{\R^d}( p_\eps P^2+q_\eps Q^2+2 r_\eps PQ)\Delta\theta_n\nonumber\\
&&+\int_{\R^d}(\nabla p_\eps P^2+ \nabla q_\eps Q^2 +2\nabla r_\eps PQ)\cdot \nabla\theta_n
-\int_{\R^d}f_\eps\cdot (\Delta P,\Delta Q)\theta_n.
\ee
Thanks to Lemma \ref{ppos}, for $n\geq 1$ and $\eps\leq 1$
sufficiently small,
\be\label{mp1}
|p_\eps \Delta\theta_n| &\lesssim &\max(1,|y_1|)|\Delta\theta_n|\lesssim\frac{1}{\eps^{2/3}}\max(1,|x|^2)\frac{1}{n^2}\left|\Delta\theta\left(\frac{x}{n}\right)\right|\nonumber\\
&\lesssim &\frac{1}{\eps^{2/3}}\max\left(\frac{1}{n^2}\|\Delta\theta\|_{L^\infty},\||x|^2\Delta\theta\|_{L^\infty}\right)\lesssim\frac{1}{\eps^{2/3}}.
\ee
Similarly, Lemma \ref{qpos} yields
\be\label{mq1}
|q_\eps \Delta\theta_n| &\lesssim &\frac{1}{\eps^{2/3}},
\ee
and since $\Delta_\eps\leq 0$ thanks to (\ref{min-deltae}), (\ref{mp1}) and (\ref{mq1}) also imply
\be\label{mr1}
|r_\eps \Delta\theta_n| &\lesssim &\frac{1}{\eps^{2/3}}.
\ee
Next, we use the estimates
\be\label{m-der-pqreps}
\max(|\nabla p_\eps|,|\nabla q_\eps|,|\nabla r_\eps|)\lesssim \max(\eps^{-4/3},|x|/\eps^{2/3}),\quad \max(|\Delta p_\eps|,|\Delta q_\eps|,|\Delta r_\eps|)\lesssim \eps^{-2},
\ee
that will be proved later. Arguing like in (\ref{mp1}), it follows
from (\ref{m-der-pqreps}) that for $n\geq 1$,
\be\label{mpqr2}
|\nabla p_\eps\cdot\nabla\theta_n|\lesssim\eps^{-4/3},
\quad |\nabla q_\eps\cdot\nabla\theta_n|\lesssim\eps^{-4/3},
\quad |\nabla r_\eps\cdot\nabla\theta_n|\lesssim\eps^{-4/3}.
\ee
Letting $n\to \infty$, and using the positivity of the quadratic form
$a(P,Q)=p_\eps P^2+q_\eps Q^2+2r_\eps PQ$, shown in (\ref{minquad})
and (\ref{min-deltae}), we deduce from (\ref{est-Delta}), (\ref{mp1}),
(\ref{mq1}), (\ref{mr1}), (\ref{m-der-pqreps}), (\ref{mpqr2}), the
Young inequality and (\ref{est-R}) that
\be\label{est-Delta2}
\eps^{4/3}\int_{\R^d}(|\Delta P|^2+|\Delta
Q|^2)&\lesssim&\eps^{-2}\left(\|P\|_{L^2}^2+\|Q\|_{L^2}^2\right)+\eps^{-4/3}\|f_\eps\|_{(L^2)^2}^2\nonumber\\
&\lesssim&\eps^{-26/3}+\eps^{-4/3}\|f_\eps\|_{L^2(\R^d)^2}^2.
\ee
Thanks to (\ref{est-fe}), (\ref{est-fe2}),
(\ref{est-fe3}) and (\ref{est-R}), the $L^2$ norm of $f_\eps$ can
estimated as
$$\|f_\eps\|_{L^2(\R^d)^2}\lesssim \eps^{-2}+\eps^{2N/3-19/3}+\eps^{4N/3-26/3}\lesssim\eps^{-2},
$$
provided $N$ is large enough. Thus, (\ref{est-Delta2}) yields
$$\|(\Delta P,\Delta Q)\|_{L^2(\R^d)^2}\lesssim \eps^{-5},$$
which combined with (\ref{est-R}), implies
\be\label{est-PQ-infini}
\|(P,Q)\|_{L^\infty(\R^d)^2}\lesssim\|(P,Q)\|_{H^2(\R^d)^2}\lesssim\eps^{-5}.
\ee
In view of the ansatz (\ref{ansatz-rho})-(\ref{ansatz-rho2}) as well
asthe estimates (\ref{min-rho1}), (\ref{min-rho2}) and
(\ref{est-PQ-infini}), we conclude that if $N$ is sufficiently large
and if $\eps$ is small enough, $\eta_1$ and $\eta_2$ are strictly
positive respectively on  $S_1$ and $S_2$.
In order to complete the proof of this last statement, it remains to
prove estimates (\ref{m-der-pqreps}). This is the issue we address
now. First, we deduce from (\ref{eq-omega0}), (\ref{eq-tau0}) and
Lemma \ref{est-omega-tau} that
\be\label{estim-omega}
\|\omega\|_{L^\infty(D_0)}\lesssim 1,\quad
\|\nabla\omega\|_{L^\infty(D_0)}\lesssim \eps^{-\beta},\quad
\|\Delta\omega\|_{L^\infty(D_0)}\lesssim \eps^{-2\beta},
\ee
\be
\|\tau\|_{L^\infty(D_0)}\lesssim 1,\quad
\|\nabla\tau\|_{L^\infty(D_0)}\lesssim 1,\quad
\|\Delta\tau\|_{L^\infty(D_0)}\lesssim 1,\label{estim-tau}
\ee
where for the estimates on $\nabla\tau_0$ and $\Delta\tau_0$, we have
used assumption (\ref{cond-disk}). From (\ref{asymp-nu0}),
(\ref{asymp-nun}), (\ref{asymp-lam0}) and (\ref{asymp-lamn}), we infer
\be\label{estim-nu}
\|\nu\|_{L^\infty(D_1)}\lesssim \eps^{\beta/2-1/3},\quad
\|\nabla\nu\|_{L^\infty(D_1)}\lesssim \eps^{-2/3},\quad
\|\Delta\nu\|_{L^\infty(D_1)}\lesssim \eps^{-4/3},
\ee
\be
\underset{x\in D_1}{\inf}\lambda^{1/2}\gtrsim \eps^{-1/3},\quad\|\lambda^{1/2}\|_{L^\infty(D_1)}\lesssim \eps^{-1/3},\quad
\|\nabla(\lambda^{1/2})\|_{L^\infty(D_1)}\lesssim \eps^{-1/3},\quad
\|\Delta(\lambda^{1/2})\|_{L^\infty(D_1)}\lesssim \eps^{-1}.\label{estim-lambda}
\ee
Note that the first estimate in (\ref{estim-lambda}) has already been
proved in Lemma \ref{minor-lam}. (\ref{mun}) and Propositions
\ref{proposition-Painleve} and \ref{proposition-mun} imply
\be\label{estim-mu}
\|\mu\|_{L^\infty(D_2)}\lesssim\eps^{-1/3},\quad \|\nabla\mu\|_{L^\infty(D_2)}\lesssim\eps^{-2/3},\quad \|\Delta\mu\|_{L^\infty(D_2)}\lesssim\eps^{-4/3}.
\ee
Moreover, it follows from their definitions that the truncation
functions $\Phi_\eps$, $\chi_\eps$ and $\Psi_\eps$ satisfy the
estimates
\be\label{estim-phichipsi}
\|\nabla\Phi_\eps\|_{L^\infty},\ \|\nabla\chi_\eps\|_{L^\infty},
\|\nabla\Psi_\eps\|_{L^\infty}\lesssim\eps^{-\beta},\quad \|\Delta\Phi_\eps\|_{L^\infty},\ \|\Delta\chi_\eps\|_{L^\infty},
\|\Delta\Psi_\eps\|_{L^\infty}\lesssim\eps^{-2\beta}
\ee
Combining (\ref{estim-omega}), (\ref{estim-tau}), (\ref{estim-nu}),
(\ref{estim-lambda}), (\ref{estim-mu}) and (\ref{estim-phichipsi}) and
using Lemmata \ref{lem-compar-om-nu-tau-l} and \ref{lem-compar-D1D2}, we obtain
\be\label{estim-rho}
\|\rho_1\|_{L^\infty(\R^d)}\lesssim\eps^{-1/3},\quad \|\rho_2\|_{L^\infty(\R^d)}\lesssim\eps^{-1/3},
\ee
\be\label{estim-nabla-rho1}
\|\nabla\rho_1\|_{L^\infty(\R^d)}=\|\Phi_\eps\eps^{-1/3}\nabla\omega+\chi_\eps\nabla\nu+\nabla\Phi_\eps(\eps^{-1/3}\omega-\nu)\|_{L^\infty(\R^d)}\lesssim\eps^{\min(-1/3-\beta,-2/3)},
\ee
and
\be\label{estim-nabla-rho2}
\|\nabla\rho_2\|_{L^\infty(\R^d)}&=&\|\Phi_\eps\eps^{-1/3}\nabla\tau+\chi_\eps\nabla(\lambda^{1/2})+\Psi_\eps\nabla\mu+\nabla\Phi_\eps(\eps^{-1/3}\tau-\lambda^{1/2})+\nabla\Psi_\eps(\mu-\lambda^{1/2})\|_{L^\infty(\R^d)}\nonumber\\
&\lesssim&\eps^{-2/3}
\ee
provided $N$ is large enough, as well as
\be\label{estim-delta-rho1}
\|\Delta\rho_1\|_{L^\infty(\R^d)}&=&\|\Phi_\eps\eps^{-1/3}\Delta\omega+\chi_\eps\Delta\nu+2\nabla\Phi_\eps(\eps^{-1/3}\nabla\omega-\nabla\nu)+\Delta\Phi_\eps(\eps^{-1/3}\omega-\nu)\|_{L^\infty(\R^d)}\nonumber\\
&\lesssim&\eps^{\min(-1/3-2\beta,-4/3)},
\ee
and
\be\label{estim-delta-rho2}
\|\Delta\rho_2\|_{L^\infty(\R^d)}&=&\|\Phi_\eps\eps^{-1/3}\Delta\tau+\chi_\eps\Delta(\lambda^{1/2})+\Psi_\eps\Delta\mu+2\nabla\Phi_\eps\nabla(\eps^{-1/3}\tau-\lambda^{1/2})+2\nabla\Psi_\eps\nabla(\mu-\lambda^{1/2})\nonumber\\
&&+\Delta\Phi_\eps(\eps^{-1/3}\tau-\lambda^{1/2})+\Delta\Psi_\eps(\mu-\lambda^{1/2})\|_{L^\infty(\R^d)}\lesssim\eps^{-4/3},
\ee
where we assume again that $N$ is sufficiently
large. (\ref{m-der-pqreps}) follows by differentiation of the
definitions of $p_\eps$, $q_\eps$ and $r_\eps$ given in Section \ref{sec-deriv}.

\paragraph{3$^{\rm rd}$ step.} First, note that the functions $\eta_1$
and $\eta_2$ we have constructed are radial. Indeed, $\rho_1$ and
$\rho_2$ are functions of the variables $z$, $y_1$ and $y_2$, which
all depend on $x$ only through $|x|$. On the other side, the equation
(\ref{eq-reste}) is radially symmetric, such that the uniqueness
of its solution $(P,Q)$ in the ball $\mathcal{B}_R$, which was proved in
Section \ref{fp}, ensures that both $P$ and $Q$ are radial. For
convenience, we consider $\eta_1$ and $\eta_2$ as functions of
$r=|x|$. At this point, according to the conclusion of the second
step, we already know that for $j=1,2$, $\eta_j(r)>0$ for $r\in
[0,(R_j^2+\eps^{2/3})^{1/2}]$. So it remains to prove that
$\eta_j(r)>0$ for $r>(R_j^2+\eps^{2/3})^{1/2}$. We shall see that it
is a consequence of Hopf's lemma (see for instance \cite{E}). Indeed,
the system of equations (\ref{sys}) satisfied by $(\eta_1,\eta_2)$ can
be rewritten as 
$$\left(-\Delta+\frac{c_j}{\eps^2}\right)(-\eta_j)=0\quad {\rm for\ }j=1,2,$$ 
where
$$c_1=2\alpha_1\eta_1^2+2\alpha_0\eta_2^2-\frac{\alpha_0}{\alpha_2}(R_2^2-R_1^2)+(|x|^2-R_1^2)$$
and
$$c_2=2\alpha_0\eta_1^2+2\alpha_2\eta_2^2+|x|^2-R_2^2.$$
Let us fix $j\in\{1,2\}$. We shall see in Lemmata \ref{c1} and \ref{c2}
below that $c_j\geq 0$ for $|x|>(R_j^2+\eps^{2/3})^{1/2}$. Let us admit
provisionnaly this fact. We know that $-\eta_j(r)<0$ for
$r<(R_j^2+\eps^{2/3})^{1/2}$. Assume by contradiction that there
exists $r_0\geq (R_j^2+\eps^{2/3})^{1/2}$ such that
$-\eta_j(r_0)=0$. Then, Hopf's Lemma applied on the ball of $\R^d$
centered at the origin and with radius $r_0$ ensures that
$-\eta_j'(r_0)>0$. In particular, $r\mapsto -\eta_j(r)$ is strictly
increasing in a neighborhood of $r_0$, in such a way that we can
define $r_1\in (r_0, +\infty]$ by 
$$r_1=\sup\{r>r_0, -\eta_j {\rm\ is\ 
    stricly\ increasing\ on\ } (r_0,r_1)\}.$$ If $r_1$ is finite, we can
  apply again Hopf's Lemma on the ball centered at 0, with radius
  $r_1$, and conclude that $-\eta_j$ is increasing on a neighborhood
  of $r_1$, which is a contradiction with the definition of
  $r_1$. Thus, $r_1=+\infty$. Thus, $-\eta_j$ is strictly
  increasing on $[r_0, +\infty)$, with $-\eta_j(r_0)=0$. This is a
    contradiction with the fact that $\eta_j(r)\to 0$ as $r\to
    +\infty$ (which is itself a consequence of the decay of the
    $\mu_n(y_2)$'s as $y_2\to -\infty$ and of $(P,Q)\in
    H^2(\R^d)$). Therefore $-\eta_j(r)<0$ for every $r>0$.

\begin{lem}\label{c1}
For $\eps>0$ sufficiently small, $c_1(x)>0$ for every $x\in
\R^d\string\ S_1$.
\end{lem}
\begin{Proof} Note first that for $\eps\leq 1$, since $\beta<2/3$,
  $\R^d\string\ S_1$ is the disjoint union of the sets $D_1\backslash 
  S_1$ and $D_2\backslash D_1$. We first consider the case  where $x\in
  D_1\backslash  S_1$. Starting from (\ref{ansatz-rho}) and
  (\ref{def-rho}), we have
\be
\eta_1&=&\eps^{1/3}(\chi_\eps\nu+\eps^{2(N+1)/3}P)\nonumber\\
&=&\eps^{1/3}(\nu+o_{L^\infty(D_1\cap
  D_2)}(\eps^\alpha)+O_{L^\infty(D_1\backslash S_1)}(\eps^{2N/3-13/3}))\nonumber\\
&=&\eps^{1/3}(\nu+O_{L^\infty(D_1\backslash S_1)}(\eps^{2N/3-13/3}))\nonumber\\
&=&\eps^{1/3}(\nu_0+O_{L^\infty(D_1\backslash S_1)}(\eps^{2/3})),
\ee  
 where the first line holds because $x\not\in {\rm Supp}(\Phi_\eps)$, the second one because $\chi_\eps\equiv 1$ on $D_1\backslash D_2$ and thanks to Lemma \ref{lem-compar-D1D2} and 
 (\ref{est-PQ-infini}), the third line holds provided $\alpha$ is
 chosen large enough, and the last line is true for $N$ large enough,
 since  $D_1\backslash   S_1\subset \{x, y_1\leq -1\}$ and thanks to the asymptotics of the $\nu_n$'s as $y_1\to -\infty$ given in Proposition \ref{summary}. The same kind of arguments yield, still for $x\in  D_1\backslash   S_1$:
 \be
 \eta_2&=&\eps^{1/3}(\chi_\eps\lambda^{1/2}+\Psi_\eps\mu+\eps^{2(N+1)/3}Q)\nonumber\\
&=&\eps^{1/3}(\lambda^{1/2}+o_{L^\infty(D_1\cap
   D_2)}(\eps^{\beta(N+1)-1/3})+O_{L^\infty(D_1\backslash S_1)}(\eps^{2N/3-13/3})\\
&=&\left(\lambda_{-1}+\eps^{2/3}\lambda_0+O_{L^\infty(D_1\backslash S_1)}(\eps^{4/3})\right)^{1/2}+O_{L^\infty(D_1\backslash S_1)}(\eps^{4/3})=\left(\lambda_{-1}+\eps^{2/3}\lambda_0+O_{L^\infty(D_1\backslash S_1)}(\eps^{4/3})\right)^{1/2}.\nonumber
\ee   
As a result, thanks to (\ref{I-1}), (\ref{II0}) and (\ref{nu00}), for  $x\in  D_1\backslash   S_1$, we have
\be
c_1&=&2\alpha_1\eps^{2/3}(\nu_0+O_{L^\infty}(\eps^{2/3}))^2+2\alpha_0\left(\lambda_{-1}+\eps^{2/3}\lambda_0+O_{L^\infty}(\eps^{4/3})\right)-\frac{\alpha_0}{\alpha_2}(R_2^2-R_1^2)-\eps^{2/3}y_1\nonumber\\
&=& R_1^{2/3}\Gamma_2^{2/3}\eps^{2/3}\gamma_0\left(\frac{\Gamma_2^{1/3}y_1}{R_1^{2/3}}\right)^2-\Gamma_2\eps^{2/3}y_1+O_{L^\infty}(\eps^{4/3})\geq \Gamma_2\eps^{2/3}+O(\eps^{4/3})>0,
\ee
for $\eps$ sufficiently small, where the inequality holds because $x\not\in   S_1$, which implies $y_1\leq -1$. Let us now consider the case where $x\in D_2\backslash D_1$. Then, using again (\ref{est-PQ-infini}),
$$\eta_1=\eps^{2N/3+1}P=O_{L^\infty(D_2)}(\eps^{2N/3-4})=O_{L^\infty(D_2)}(\eps^{4/3})$$
and
\be
\eta_2=\eps^{1/3}(\mu+\eps^{2(N+1)/3}Q)=\eps^{1/3}(\mu_0+O_{L^\infty(D_2)}(\eps^{2/3})).
\ee
We infer that
\be
c_1 & =&2\alpha_0\eps^{2/3}\mu_0^2-\frac{\alpha_0}{\alpha_2}(R_2^2-R_1^2)+(|x|^2-R_1^2)+O_{L^\infty(D_2)}(\eps^{4/3})\nonumber\\
& =&\eps^{2/3}\left(\frac{\alpha_0}{\alpha_2}R_2^{2/3}\gamma_0\left(\frac{y_2}{R_2^{2/3}}\right)^2-\frac{\alpha_0}{\alpha_2}\frac{R_2^2-R_1^2}{\eps^{2/3}}-y_2+\frac{R_2^2-R_1^2}{\eps^{2/3}}\right)+O_{L^\infty(D_2)}(\eps^{4/3})\nonumber\\
&
=&\eps^{2/3}\left(\frac{\alpha_0}{\alpha_2}R_2^{2/3}\underbrace{\left(\gamma_0\left(\frac{y_2}{R_2^{2/3}}\right)^2-\frac{y_2}{R_2^{2/3}}\right)}_{\geq
  \inf_{y\in\R}\left[\gamma_0(y)^2-y\right]>-\infty}
+\Gamma_2
\underbrace{\left(\frac{R_2^2-R_1^2}{\eps^{2/3}}-y_2\right)}_{\geq
  2\eps^{\beta-2/3}\underset{\eps\to 0}{\longrightarrow}+\infty}\right)+O_{L^\infty(D_2)}(\eps^{4/3}),
\ee
and thus $c_1>0$ on $D_2\backslash D_1$ if $\eps$ is sufficiently small.
\end{Proof}

\begin{lem}\label{c2}
$c_2(x)>0$ for every $x\in
\R^d\string\ S_2$.
\end{lem}
\begin{Proof}
The lemma is a straightforward consequence of the definition of $c_2$,
since the assumption $x\in \R^d\string\ S_2$ can be rewritten as $|x|^2>R_2^2+\eps^{2/3}$.
\end{Proof}

\subsection{Uniqueness of the ground state}\label{sec-uniq}
In this section, we prove that the solution of (\ref{sys}) constructed
in the previous sections is the unique ground state of the system,
that is the unique solution of (\ref{sys}) with two positive
components. Uniqueness of the ground state of (\ref{sys}) was proved
in \cite{ANS}. We recall the arguments for the sake of
completeness. First, the next lemma gives an a priori upper bound on positive
solutions to (\ref{sys}).

\begin{lem}\label{est-unif}
Let $\eps>0$, and let $(\eta_1,\eta_2)$ be a positive solution of (\ref{sys}). 
Then, for every $\theta\in (0,1)$ and $x\in\R^d$, for $j=1,2$,
\be
\eta_j(x)\leq M_j\min\left[1,\exp\left(-\frac{\theta}{2\eps}\left(|x|^2-r_j^2\right)\right)\right],
\ee
where $a_1=\left(\frac{\alpha_0}{\alpha_2}(R_2^2-R_1^2)+R_1^2\right)^{1/2}$, $M_1=\frac{a_1}{(2\alpha_1)^{1/2}}$, $r_1=\frac{a_1}{(1-\theta^2)^{1/2}}$, $M_2=\frac{R_2}{(2\alpha_2)^{1/2}}$  and $r_2=\frac{R_2}{(1-\theta^2)^{1/2}}$.
\end{lem}

\begin{Proof}
We first prove that $\eta_1$ is uniformly bounded from above by the
constant $M_1$ defined in the statement of the lemma. The proof
follows an idea which is due to Farina \cite{F}, and which is also used in
\cite{IM} and \cite{ANS}. Let us define 
$$w_1=\frac{1}{\eps}\left((2\alpha_1)^{1/2}\eta_1-a_1\right),\quad
     {\rm and}\quad w_1^+=\max(0,w_1).$$
Then, Kato's inequality yields
\be
\Delta w_1^+ &\geq &\mathbf{1}_{\{w_1\geq 0\}}\Delta w_1\nonumber\\
&=&
\frac{\mathbf{1}_{\{w_1\geq 0\}}}{\eps^3}(2\alpha_1)^{1/2}\eta_1\left(2\alpha_1\eta_1^2+2\alpha_0\eta_2^2-\frac{\alpha_0}{\alpha_2}(R_2^2-R_1^2)-(R_1^2-|x|^2)\right)\nonumber\\
&\geq&\frac{\mathbf{1}_{\{w_1\geq 0\}}}{\eps^3}(\eps
w_1+a_1)\left((\eps
w_1+a_1)^2-(a_1)^2\right)\nonumber\\
&=& \frac{ \mathbf{1}_{\{w_1\geq 0\}}}{\eps^3}(\eps
w_1+a_1)\eps w_1(\eps
w_1+2a_1)
\nonumber\\
&\geq&(w_1^+)^3.
\ee
From Lemma 2 in \cite{B}, it follows that $w_1^+\leq 0$,
which means $\eta_1\leq M_1$. 

Next, like it was done in \cite{ANS}, we prove that $\eta_1$
 decays at least as fast as a gaussian as $|x|$ goes to
infinity. Easy calculations show that
\be\label{max-exp}
\left(-\Delta+\frac{\theta^2}{\eps^2}|x|^2\right)\exp\left(-\frac{\theta|x|^2}{2\eps}\right)=\frac{d\theta}{\eps}\exp\left(-\frac{\theta|x|^2}{2\eps}\right)>0,
\ee
whereas
\be\label{max-eta1}
\left(-\Delta+\frac{\theta^2}{\eps^2}|x|^2\right)\eta_1=\frac{1}{\eps^2}\left(a_1^2-(1-\theta^2)|x|^2\right)\eta_1-\frac{2\alpha_1}{\eps^2}\eta_1^3-\frac{2\alpha_0}{\eps^2}\eta_2^2\eta_1<0
\ee
for $x>r_1$. Then, we set
$W_1=M_1\exp\left(-\frac{\theta(|x|^2-r_1^2)}{2\eps}\right)-\eta_1$. We
know from the first part of the proof that $W_1(x)\geq 0$ for
$|x|=r_1$. Assume by contradiction that the inequality $W_1\geq 0$ does not
hold for every $x\in \R^d$ such that $|x|\geq r_1$. Then, since
$W_1(x)\to 0$ as $|x|\to \infty$, $W_1$ reaches a minimum at some
$x_0\in\R^d$ such that $|x_0|>r_1$. In particular, $\Delta
W_1(x_0)\geq 0$ and $W_1(x_0)<0$. This is in contradiction with the
difference between (\ref{max-exp}) multiplied by $M_1$ and (\ref{max-eta1}) evaluated
at $x_0$. The proof of the estimate on $\eta_2$ is
similar. 
\end{Proof}

The next lemma states the uniqueness of the ground state of
(\ref{sys}) and is also proved in \cite{ANS}. We give here a proof which is slightly simpler.

\begin{lem}
Let $\eps>0$, and let  $(\eta_1,\eta_2)$, $(\xi_1,\xi_2)$ be two positive solutions of (\ref{sys}). Then $\eta_1=\xi_1$ and $\eta_2=\xi_2$.
\end{lem}

\begin{Proof}
Let $v_1=\xi_1/\eta_1$ and $v_2=\xi_2/\eta_2$. Since $(\eta_1,\eta_2)$
and $(\xi_1,\xi_2)$ solve (\ref{sys}), it follows that for
$(i,j)=(1,2)$ or $(2,1)$, we have
\be\label{eq-vi}
\eps^2{\rm div}\left(\eta_i^2\nabla v_i\right) &=&2\alpha_i\eta_i^4v_i(v_i^2-1)+2\alpha_0\eta_i^2\eta_j^2v_i(v_j^2-1). 
\ee
Let $\zeta\in \mathcal{C}^\infty(\R^d)$ be a non-negative function
supported in $\{x\in\R^d,\ |x|\leq 2\}$ such that $\zeta(x)=1$ for
$|x|\leq 1$. For $n\geq 1$, we also define
$\zeta_n=\zeta(\cdot/n)$. Next, let us multiply (\ref{eq-vi}) by
$(v_i^2-1)\zeta_n^2/v_i$, sum over $\R^d$ and use integration by
parts. We obtain
\be\label{eq-vi-int}
\lefteqn{\int_{\R^d}\eta_i^2|\nabla
v_i|^2\left(1+\frac{1}{v_i^2}\right)\zeta_n^2dx
+\frac{2}{\eps^2}\int_{\R^d}\left[\alpha_i\left(\eta_i^2(v_i^2-1)\right)^2+\alpha_0\eta_i^2(v_i^2-1)\eta_j^2(v_j^2-1)\right]\zeta_n^2dx}\nonumber\\
&=&
-\int_{\R^d}\eta_i^2\nabla
v_i\left(v_i-\frac{1}{v_i}\right)\nabla\left(\zeta_n^2\right)dx\hspace{7cm}\nonumber\\
&=&
-2\int_{\R^d}\eta_i^2v_i\nabla
v_i\zeta_n\nabla\zeta_n dx+2\int_{\R^d}\eta_i^2\frac{\nabla
v_i}{v_i}\zeta_n\nabla\zeta_n dx.
\ee
Next, we estimate each integral in the right hand side of
(\ref{eq-vi-int}) thanks to the Cauchy-Schwarz inequality. For the
first one, we have
\be\label{int-1}
\left|\int_{\R^d}\eta_i^2v_i\nabla
v_i\zeta_n\nabla\zeta_n dx\right|&=&\left|\int_{\R^d}\eta_i\xi_i \nabla
v_i\zeta_n\nabla\zeta_n dx\right|\nonumber\\
&\leq &\left(\int_{\R^d}\eta_i^2|\nabla
v_i|^2\zeta_n^2dx\right)^{1/2}\left(\int_{\R^d}\xi_i^2|\nabla\zeta_n
|^2dx\right)^{1/2}\nonumber\\
&\leq &\frac{1}{4}\int_{\R^d}\eta_i^2|\nabla
v_i|^2\zeta_n^2dx+\int_{\R^d}\xi_i^2|\nabla\zeta_n
|^2dx,
\ee
whereas for the second one, we get
\be\label{int-2}
\left|\int_{\R^d}\eta_i^2\frac{\nabla
v_i}{v_i}\zeta_n\nabla\zeta_n dx\right|
&\leq&\left(\int_{\R^d}\eta_i^2\frac{|\nabla
v_i|^2}{v_i^2}\zeta_n^2dx\right)^{1/2}\left(\int_{\R^d}\eta_i^2|\nabla\zeta_n
|^2dx\right)^{1/2}\nonumber\\
&\leq &\frac{1}{4}\int_{\R^d}\eta_i^2\frac{|\nabla
v_i|^2}{v_i^2}\zeta_n^2dx+\int_{\R^d}\eta_i^2|\nabla\zeta_n
|^2dx.
\ee
Combining (\ref{eq-vi-int}), (\ref{int-1}) and (\ref{int-2}), we infer
\be\label{eq-vi-int-2}
\lefteqn{\frac{1}{2}\int_{\R^d}\eta_i^2|\nabla
v_i|^2\left(1+\frac{1}{v_i^2}\right)\zeta_n^2dx
+\frac{2}{\eps^2}\int_{\R^d}\left[\alpha_i\left(\eta_i^2(v_i^2-1)\right)^2+\alpha_0\eta_i^2(v_i^2-1)\eta_j^2(v_j^2-1)\right]\zeta_n^2dx}\nonumber\\
&\leq &2\int_{\R^d}(\xi_i^2+\eta_i^2)|\nabla\zeta_n
|^2dx.\hspace{9cm}
\ee
Finally, we sum the inequalities given by (\ref{eq-vi-int-2}) for
$(i,j)=(1,2)$ and for $(i,j)=(2,1)$. We deduce
\be\label{eq-unic}
\lefteqn{\frac{1}{2}\int_{\R^d}\eta_1^2|\nabla
v_1|^2\left(1+\frac{1}{v_1^2}\right)\zeta_n^2dx
+\frac{1}{2}\int_{\R^d}\eta_2^2|\nabla  
v_2|^2\left(1+\frac{1}{v_2^2}\right)\zeta_n^2dx+\frac{2}{\eps^2}
\int_{\R^d}q\left[\eta_1^2(v_1^2-1),
  \eta_2^2(v_2^2-1)\right]\zeta_n^2dx}\nonumber\\  
&\leq &2\int_{\R^d}(\xi_1^2+\eta_1^2+\xi_2^2+\eta_2^2)|\nabla\zeta_n
|^2dx,\hspace{9cm}
\ee
where $q[u_1,u_2]=\alpha_1 u_1^2+2\alpha_0 u_1 u_2+\alpha_2
u_2^2$. Note that the assumption $\Gamma_{12}>0$ can be rewritten as
$\alpha_0^2-\alpha_1\alpha_2<0$, which implies that there exists $c>0$
such that for every $u_1,u_2\in \R$, $q[u_1,u_2]\geq
c(u_1^2+u_2^2)$. As a result, in order to conclude that $v_1\equiv
v_2\equiv 1$, it is sufficient to prove that the right hand side of
(\ref{eq-unic}) converges to 0 as $n\to \infty$. It is the case thanks
to Lemma \ref{est-unif}. Indeed, for $n\geq \max(r_1,r_2)$, since
$\nabla\zeta_n$ is supported in $\{x\in\R^d, n\leq |x|\leq 2n\}$, we have
$$\int_{\R^d}(\xi_1^2+\eta_1^2+\xi_2^2+\eta_2^2)|\nabla\zeta_n
|^2dx\leq
2\left(M_1^2e^{\theta r_1^2/\eps}+M_2^2e^{\theta
  r_2^2/\eps}\right)\|\nabla\zeta\|_{L^\infty}^2\left|\{x\in\R^d,1\leq 
|x|\leq 2\}\right|n^{d-2}e^{-n^2\theta /\eps},$$
where the right hand side goes to 0 as $n\to \infty$.
\end{Proof}

\subsection{End of the proof of Theorem \ref{main}}
In section \ref{fp}, we have constructed a solution $(\eta_1,\eta_2)$
to (\ref{sys}) that converges to 0 at infinity. In section
\ref{positivity}, we have checked that this solution is positive. In
section \ref{sec-uniq}, we have seen that $(\eta_1,\eta_2)$ is in fact
the unique such solution of (\ref{sys}). So the first part of the
statement of Theorem \ref{main} has already been proved. Let us now
fix three integers $M_0,N_0$ and $L_0$, as well as $\beta\in
(0,\eps^{2/3})$. According to our construction of $(\eta_1,\eta_2)$
explained in sections \ref{sec-deriv} and \ref{fp}, provided $M> M_0$,
$N> \max(N_0,2/\beta-2)$ and $L> L_0$ are large integers that satisfy the conditions
listed at the beginning of sections \ref{sec-deriv} and \ref{fp},
$(\eta_1,\eta_2)$ can be written like in the ansatz
(\ref{ansatz})-(\ref{troncat}). Thus, defining $\eta_{1app}$ and
$\eta_{2app}$ as in the statement of Theorem \ref{main}, we have 
\be\label{N-N0-1}
\eta_1-\eta_{1app} 
=\Phi_\eps\!\!\!\!\sum_{m=M_0+1}^{M}\!\!\!\eps^{2m}\omega_m
+\eps^{1/3}\chi_\eps\!\!\!\sum_{n=N_0+1}^{N}\!\!\!\eps^{2n/3}\nu_n 
+\eps^{2N/3+1}P\ \ \ \ 
\ee
and
\be\label{N-N0-2}
\eta_2-\eta_{2app}
&=&\Phi_\eps\sum_{m=M_0+1}^{M}\eps^{2m}\tau_m
+\eps^{1/3}\chi_\eps\left(\left(\sum_{n=-1}^{N}\eps^{2n/3}\lambda_n\right)^{1/2}
 -\left(\sum_{n=-1}^{N_0}\eps^{2n/3}\lambda_n\right)^{1/2}\right)\nonumber\\  
&&+\eps^{1/3}\Psi_\eps\sum_{n=L_0+1}^{L}\eps^{2n/3}\mu_n  
 +\eps^{2N/3+1}Q.
\ee
The next step consists in evaluating the $L^p$ and $H^1$ norms of each
term in the right hand side of (\ref{N-N0-1}) and (\ref{N-N0-2}). Let
us start with the $L^p$ norm of $\Phi_\eps\omega_m$, for $m\geq 1$ and
$p\in [2,+\infty)$. Since
${\rm Supp}\Phi_\eps\subset D_0$, we have
\be\label{pre-est-omegam-lp}
\left\|\Phi_\eps\omega_m\right\|_{L^p(\R^d)}^p &\leq &\int_{|x|^2\leq
  R_1^2-\eps^\beta}|\omega_m(z)|^pdx=\int_{\mathbb{S}^{d-1}}
\int_0^{(R_1^2-\eps^\beta)^{1/2}}|\omega_m(R_1^2-r^2)|^pr^{d-1}drd\theta\nonumber\\  
&=&
|\mathbb{S}^{d-1}|\int_{\eps^\beta}^{ R_1^2}
|\omega_m(z)|^p(R_1^2-z)^{d/2-1}\frac{dz}{2}.  
\ee
Since $d/2-1\geq -1/2$, the integral converges at $z=R_1^2$. Moreover,
thanks to (\ref{asymp-omegan}), we deduce 
$$\int_{\eps^\beta}^{R_1^2}|\omega_m(z)|^p(R_1^2-z)^{d/2-1}dz\underset{\eps\to
  0}{\sim}R_1^{d-2}|\omega_{m0}|^p\int_{\eps^\beta}^{R_1^2}z^{p(1/2-3m)}dz\underset{\eps\to
  0}{\sim}\frac{R_1^{d-2}|\omega_{m0}|^p}{p(3m-1/2)-1}\eps^{\beta(
  -p(3 m-1/2)+1)}.$$
As a result,
\be\label{est-omegam-lp}
\left\|\Phi_\eps\omega_m\right\|_{L^p(\R^d)}=\mathcal{O}(\eps^{-3\beta
  m+\beta/2+\beta/p}). 
\ee
Similarly, (\ref{asymp-taun}) yields
\be\label{est-taum-lp}
\left\|\Phi_\eps\tau_m\right\|_{L^p(\R^d)}=\mathcal{O}(\eps^{-3\beta m+\beta+\beta/p}).
\ee
Note that (\ref{est-omegam-lp}) and (\ref{est-taum-lp}) also hold for
$p=+\infty$ thanks to (\ref{asymp-omegan}) and
(\ref{asymp-taun}). Note also that (\ref{est-omegam-lp}) and
(\ref{est-taum-lp}) are sharp. Indeed, since $\Phi_\eps\equiv 1$ for
$|x|^2\leq R_1^2-2\eps^\beta$, we deduce that
$\left\|\Phi_\eps\omega_m\right\|_{L^p(\R^d)}^p$ can be bounded from
below by an integral similar to the one that appears in the right hand
side of (\ref{pre-est-omegam-lp}).

Next, let us estimate the $L^p$ norm of
$\chi_\eps\nu_n$, for $n\geq 1$ and $p\in [2,+\infty)$. Since ${\rm
    Supp}\chi_\eps\subset D_1$, we get
\be\label{pre-est-nun-lp}
\left\|\chi_\eps\nu_n\right\|_{L^p(\R^d)}^p &\leq
&\int_{R_1^2-2\eps^\beta\leq |x|^2\leq
  R_1^2+2\eps^\beta}|\nu_n(y_1)|^pdx=\int_{\mathbb{S}^{d-1}}
\int_{(R_1^2-2\eps^\beta)^{1/2}}^{(R_1^2+2\eps^\beta)^{1/2}}\left|\nu_n\left(\frac{R_1^2-r^2}{\eps^{2/3}}\right)\right|^pr^{d-1}drd\theta\nonumber\\   
&=&
|\mathbb{S}^{d-1}|\int_{-2\eps^{\beta-2/3}}^{ 2\eps^{\beta-2/3}}
|\nu_n(y_1)|^p(R_1^2-\eps^{2/3}y_1)^{d/2-1}\frac{\eps^{2/3}dy_1}{2}.  
\ee
For $y_1\in [-2\eps^{\beta-2/3},2\eps^{\beta-2/3}]$, we have $1\lesssim R_1^2-\eps^{2/3}y_1\lesssim
1$, therefore according to the asymptotic behaviour of $\nu_n(y_1)$ as
$y_1\to\pm\infty$ given in Proposition \ref{summary}, we obtain
\be\label{est-nun-lp}
\left\|\chi_\eps\nu_n\right\|_{L^p(\R^d)} & =&
\left\{\begin{array}{ll}
\mathcal{O}(\eps^{\frac{2}{3p}}) & {\rm if\ } n=1\ {\rm or}\  (n=2\ {\rm and}\ p>2)\\
 \mathcal{O}(|\ln\eps|^{\frac{1}{2}}\eps^{\frac{1}{3}}) & {\rm if\ }
 n=2\ {\rm and}\ p=2\\
\mathcal{O}(\eps^{-\frac{2n}{3}+\beta(n-\frac{5}{2})+\frac{\beta}{p}+\frac{5}{3}})
  &  {\rm if\ } n\geq 3.
\end{array}\right.
\ee
Similarly,
\be\label{est-lambdan-lp}
\left\|\chi_\eps\lambda_n\right\|_{L^p(\R^d)} & =&
\left\{\begin{array}{ll}
\mathcal{O}(\eps^{\frac{2}{3p}}) & {\rm if\ } n=1\\
\mathcal{O}(\eps^{-\frac{2n}{3}+\beta(n-2)+\frac{\beta}{p}+\frac{4}{3}})
  &  {\rm if\ } n\geq 2.
\end{array}\right.
\ee
Again, it easily follows from Proposition  \ref{summary} that
(\ref{est-nun-lp}) and (\ref{est-lambdan-lp}) also hold for
$p=+\infty$, and the two estimates are sharp. Next, since ${\rm
    Supp}\Psi_\eps\subset D_2$, we infer
\be\label{pre-est-mun-lp}
\left\|\Psi_\eps\mu_n\right\|_{L^p(\R^d)}^p &\leq
&\int_{|x|^2\geq
  R_1^2+\eps^\beta}|\mu_n(y_2)|^pdx=\int_{\mathbb{S}^{d-1}}
\int_{(R_1^2+\eps^\beta)^{1/2}}^{+\infty}\left|\mu_n\left(\frac{R_2^2-r^2}{\eps^{2/3}}\right)\right|^pr^{d-1}drd\theta\nonumber\\   
&=&
|\mathbb{S}^{d-1}|\int_{-\infty}^{ \frac{R_2^2-R_1^2}{\eps^{2/3}}-\eps^{\beta-2/3}}
|\mu_n(y_2)|^p(R_2^2-\eps^{2/3}y_2)^{d/2-1}\frac{\eps^{2/3}dy_2}{2}.  
\ee
In order to estimate the integral in the right hand side of
(\ref{pre-est-mun-lp}), we split the integral into two pieces. First,
for $y_2\in
(-R_2^2/\eps^{2/3},(R_2^2-R_1^2)/\eps^{2/3}-\eps^{\beta-2/3})$, we
have
$1\lesssim R_1^2+\eps^\beta\leq R_2^2-\eps^{2/3}y_2\leq 2R_2^2\lesssim
1$. Therefore, according to (\ref{mun}) and Proposition
\ref{proposition-mun},
\be\label{est-mun-lp-bout1}
\int_{-R_2^2/\eps^{2/3}}^{ \frac{R_2^2-R_1^2}{\eps^{2/3}}-\eps^{\beta-2/3}}
|\mu_n(y_2)|^p(R_2^2-\eps^{2/3}y_2)^{d/2-1}dy_2\underset{\eps\to
  0}{=}\mathcal{O}(1).
\ee
If $d=1,2$ and $y_2\leq -R_2^2/\eps^{2/3}$, we still have
$(R_2^2-\eps^{2/3}y_2)^{d/2-1}\lesssim 1$, therefore
\be\label{est-mun-lp-bout2}
\int_{-\infty}^{-R_2^2/\eps^{2/3}}
|\mu_n(y_2)|^p(R_2^2-\eps^{2/3}y_2)^{d/2-1}dy_2\underset{\eps\to
  0}{=}\mathcal{O}(1),
\ee
whereas if $d=3$ and
$y_2\leq -R_2^2/\eps^{2/3}$, then $(R_2^2-\eps^{2/3}y_2)^{1/2}\leq
\sqrt{2}\eps^{1/3}|y_2|^{1/2}$ and since Proposition
\ref{proposition-mun} implies $\mu_n(y_2)\underset{y_2\to
  -\infty}{=}\mathcal{O}(|y_2|^{-5/2})$, we deduce that
(\ref{est-mun-lp-bout2}) also holds. Combining (\ref{pre-est-mun-lp}), 
(\ref{est-mun-lp-bout1}) and (\ref{est-mun-lp-bout2}), we deduce 
\be\label{est-mun-lp}
\left\|\Psi_\eps\mu_n\right\|_{L^p(\R^d)}=\mathcal{O}(\eps^{\frac{2}{3p}}).
\ee
Note again that (\ref{est-mun-lp}) is sharp and that Proposition \ref{proposition-mun} implies that it is also true for
$p=+\infty$. 

We are now ready to estimate $\eta_1-\eta_{1app}$ and
$\eta_2-\eta_{2app}$ in $L^p(\R^d)$. Remark first that since $\beta<2/3$,
(\ref{est-omegam-lp}) and (\ref{est-taum-lp}) imply that the larger is
$m\geq 1$, the smaller are $\eps^{2m}\Phi_\eps\omega_m$ and
$\eps^{2m}\Phi_\eps\tau_m$ in $L^p(\R^d)$, in the limit $\eps\to 0$. Similarly, since $\beta>0$,
it follows from (\ref{est-nun-lp}) and (\ref{est-lambdan-lp}) that the
larger is $n$, the smaller are $\eps^{2n/3}\chi_\eps\nu_n$ and
$\eps^{2n/3}\chi_\eps\lambda_n$ in $L^p(\R^d)$. Thus,
\be\label{pre-conclusion-1}
\left\|\eta_1-\eta_{1app}\right\|_{L^p(\R^d)}&=&\mathcal{O}(\eps^{(2-3\beta)(M_0+1)+\frac{\beta}{2}+\frac{\beta}{p}})+\left\{\begin{array}{ll}
\mathcal{O}(\eps^{1+\frac{2}{3p}}) & {\rm if\ } N_0=0\\
\mathcal{O}(\eps^{\frac{5}{3}+\frac{2}{3p}}) & {\rm if}\  N_0=1\ {\rm and}\ p>2\\
 \mathcal{O}(|\ln\eps|^{\frac{1}{2}}\eps^{2}) & {\rm if\ }
 N_0=1\ {\rm and}\ p=2\\
\mathcal{O}(\eps^{\beta(N_0-\frac{3}{2})+\frac{\beta}{p}+2})
  &  {\rm if\ } N_0\geq 2.
\end{array}\right\}\nonumber\\
&&+\|P\|_{L^p(\R^d)}\mathcal{O}(\eps^{2N/3+1}).
\ee
Now, remember that in (\ref{est-PQ-infini}), the $H^2(\R^d)$ norm of
$P$ is controled by some power of $\eps$ (namely, $\eps^{-5}$) wich is
independent
of $N$. Thus, thanks to Sobolev embeddings, for
fixed values of $M_0$, $N_0$ and $L_0$, if
$M,N$ and $L$ are chosen sufficiently large (and such that they
satisfy the conditions at the beginning of sections \ref{sec-deriv}
and \ref{fp}), for $\eps$ small, $\eps^{2N/3+1}
\|P\|_{L^p(\R^d)}$ becomes negligible in comparison with the other
terms in the right hand side of (\ref{pre-conclusion-1}). The estimate
on $\eta_1-\eta_{1app}$ in (\ref{est-main}) follows in the case
$E=L^p(\R^d)$.

As for the second component, using the same arguments, we infer from (\ref{est-taum-lp}), and (\ref{est-mun-lp}) that
\be\label{pre-pre-conclusion-2}
\left\|\eta_2-\eta_{2app}\right\|_{L^p(\R^d)}&=&\mathcal{O}(\eps^{(2-3\beta)(M_0+1)+\beta+\frac{\beta}{p}})\nonumber\\
&&+\chi_\eps\left(\left(\lambda_{-1}+\sum_{n=1}^{N+1}\eps^{2n/3}\lambda_{n-1}\right)^{1/2}
 -\left(\lambda_{-1}+\sum_{n=1}^{N_0+1}\eps^{2n/3}\lambda_{n-1}\right)^{1/2}\right)\nonumber\\
&&+\mathcal{O}(\eps^{\frac{1}{3}+\frac{2(L_0+1)}{3}+\frac{2}{3p}})+\|Q\|_{L^p(\R^d)}\mathcal{O}(\eps^{2N/3+1}).
\ee
In order to estimate the second term in the right hand side, note that thanks to the asymptotic behaviour of $\lambda_0$ given in Proposition \ref{summary} and (\ref{est-lambdan-lp}) for $p=+\infty$, we have
$$\sum_{n=1}^{N+1}\eps^{2n/3}\lambda_{n-1}=\mathcal{O}_{L^\infty(D_1)}(\eps^\beta),$$
and the same property holds for $N$ replaced by $N_0$. Thus, the mean value theorem
applied to the function square root close to $\lambda_{-1}$ and (\ref{est-lambdan-lp}) imply
\be\label{pre-conclusion-2}
\left\|\eta_2-\eta_{2app}\right\|_{L^p(\R^d)}
&=&\mathcal{O}(\eps^{(2-3\beta)(M_0+1)+\beta+\frac{\beta}{p}})+\left\{\begin{array}{ll}
\mathcal{O}(\eps^{\frac{4}{3}+\frac{2}{3p}}) & {\rm if\ } N_0=0\\
\mathcal{O}(\eps^{2+\beta(N_0-1)+\frac{\beta}{p}})
  &  {\rm if\ } N_0\geq 1.
\end{array}\right.\nonumber\\
&&+\mathcal{O}(\eps^{\frac{1}{3}+\frac{2(L_0+1)}{3}+\frac{2}{3p}}),
\ee
under the same condition of largeness on $M,N,L$ than for the estimate
on $\eta_1-\eta_{1app}$. We have proved the estimate on
$\eta_2-\eta_{2app}$ in (\ref{est-main}) for $E=L^p(\R^d)$.

Next, let us prove (\ref{est-main}) for $E=H^1(\R^d)$. For this
purpose, we have to estimate the $L^2(\R^d)$ norms of
$\nabla(\Phi_\eps\omega_m)$, $\nabla(\Phi_\eps\tau_m)$,
$\nabla(\chi_\eps\nu_n)$, $\nabla(\chi_\eps\lambda_n)$ and
$\nabla(\Psi_\eps\mu_n)$ for $m,n\geq 1$. In view of the definitions
of $\Phi_\eps$, $\chi_\eps$ and $\Psi_\eps$, it is clear that the
$L^\infty(\R^d)^d$ norms of their
gradients are all $\mathcal{O}(\eps^{-\beta})$. Thus, performing
calculations similar to the ones which were done to obtain
(\ref{est-omegam-lp}), (\ref{est-taum-lp}),
(\ref{est-nun-lp}),(\ref{est-lambdan-lp}) and (\ref{est-mun-lp}), we
obtain
\be
\label{est-nabla1-omegam-lp}
\left\|\nabla(\Phi_\eps)\omega_m\right\|_{L^2(\R^d)}&=&\mathcal{O}(\eps^{-3\beta
  m}), \\
\label{est-nabla1-taum-lp}
\left\|\nabla(\Phi_\eps)\tau_m\right\|_{L^2(\R^d)}&=&\mathcal{O}(\eps^{-3\beta
  m+\beta/2}),\\
\label{est-nabla1-nun-lp}
\left\|\nabla(\chi_\eps)\nu_n\right\|_{L^2(\R^d)} & =&
\left\{\begin{array}{ll}
\mathcal{O}(\eps^{\frac{1}{3}-\beta}) & {\rm if\ } n=1\\
 \mathcal{O}(|\ln\eps|^{\frac{1}{2}}\eps^{\frac{1}{3}-\beta}) & {\rm
   if\ } n=2\\ 
\mathcal{O}(\eps^{-\frac{2n}{3}+\beta(n-3)+\frac{5}{3}})
  &  {\rm if\ } n\geq 3
\end{array}\right.\\
\label{est-nabla1-lambdan-lp}
\left\|\nabla(\chi_\eps)\lambda_n\right\|_{L^2(\R^d)} & =&
\left\{\begin{array}{ll}
\mathcal{O}(\eps^{\frac{1}{3}-\beta}) & {\rm if\ } n=1\\
\mathcal{O}(\eps^{-\frac{2n}{3}+\beta(n-5/2)+\frac{4}{3}})
  &  {\rm if\ } n\geq 2
\end{array}\right.\\
\label{est-nabla1-mun-lp}
\left\|\nabla(\Psi_\eps)\mu_n\right\|_{L^2(\R^d)}&=&\mathcal{O}(\eps^{\frac{1}{3}-\beta}).
\ee
By differentiation of (\ref{asymp-omegan}) and (\ref{asymp-taun}),
since $\nabla=-2x\frac{d}{dz}$, similar calculations as the ones that
gave (\ref{est-omegam-lp}) and (\ref{est-taum-lp}) yield
\be
\label{est-nabla2-omegam-lp}
\left\|\Phi_\eps\nabla\omega_m\right\|_{L^2(\R^d)}&=&\mathcal{O}(\eps^{-3\beta
  m}), \\
\label{est-nabla2-taum-lp}
\left\|\Phi_\eps\nabla\tau_m\right\|_{L^2(\R^d)}&=&\mathcal{O}(\eps^{-3\beta
  m+\beta/2}).
\ee
Next, since $\nabla=\frac{-2x}{\eps^{2/3}}\frac{d}{dy_1}$, a calculation
similar to (\ref{pre-est-nun-lp}) yields
\be
\left\|\chi_\eps\nabla\nu_n\right\|_{L^2(\R^d)}^2 &\leq
&
\frac{2|\mathbb{S}^{d-1}|}{\eps^{2/3}}\int_{-2\eps^{\beta-2/3}}^{ 2\eps^{\beta-2/3}}
|\nu_n'(y_1)|^2(R_1^2-\eps^{2/3}y_1)^{d/2}dy_1.  
\ee
Then, after differentiation of (\ref{asymp-nun}), we deduce that
\be\label{est-nabla2-nun-lp}
\left\|\chi_\eps\nabla\nu_n\right\|_{L^2(\R^d)} & =&
\left\{\begin{array}{ll}
\mathcal{O}(\eps^{-1/3}) & {\rm if\ } n=1\ {\rm or}\ n=2\\
 \mathcal{O}(|\ln\eps|^{\frac{1}{2}}\eps^{-\frac{1}{3}}) & {\rm
   if\ } n=3\\ 
\mathcal{O}(\eps^{-2n/3+\beta(n-3)+5/3})
  &  {\rm if\ } n\geq 4.
\end{array}\right.
\ee
Similarly, differentiation of (\ref{asymp-lamn}) yields
\be\label{est-nabla2-lambdan-lp}
\left\|\chi_\eps\nabla\lambda_n\right\|_{L^2(\R^d)} & =&
\left\{\begin{array}{ll}
\mathcal{O}(\eps^{-1/3}) & {\rm if\ } n=1\ {\rm or}\ n=2\\
\mathcal{O}(\eps^{-2n/3+\beta(n-5/2)+4/3})
  &  {\rm if\ } n\geq 3.
\end{array}\right.
\ee
Using Proposition \ref{mun} like it was done to obtain
(\ref{est-mun-lp}), we infer 
\be\label{est-nabla2-mun-lp}
\left\|\Psi_\eps\nabla\mu_n\right\|_{L^2(\R^d)}&=&\mathcal{O}(\eps^{-1/3}).
\ee
Like in (\ref{pre-conclusion-1}), taking $M,N$ and $L$
large enough, we deduce from (\ref{est-nabla1-omegam-lp}),
(\ref{est-nabla1-nun-lp}), (\ref{est-nabla2-omegam-lp}) and
(\ref{est-nabla2-nun-lp}) that
\be\label{est-nabla-eta1}
\|\nabla(\eta_1-\eta_{1app})\|_{L^2(\R^d)}&=&\mathcal{O}(\eps^{(2-3\beta)(M_0+1)})+\left\{\begin{array}{ll}
\mathcal{O}(\eps^{2N_0/3+2/3}) & {\rm if\ } N_0=0\ {\rm or}\ N_0=1\\
 \mathcal{O}(|\ln\eps|^{\frac{1}{2}}\eps^{2}) & {\rm
   if\ } N_0=2\\
\mathcal{O}(\eps^{\beta(N_0-2)+2})
  &  {\rm if\ } N_0\geq 3.
\end{array}\right.
\ee
Next, we write
\be\label{nabla-N-N0-2}
\nabla(\eta_2-\eta_{2app})
&=&\underbrace{\nabla\left(\Phi_\eps\sum_{m=M_0+1}^{M}\eps^{2m}\tau_m\right)}_{=:T_1}
+\underbrace{\eps^{1/3}\nabla\chi_\eps\left(\left(\sum_{n=-1}^{N}\eps^{2n/3}\lambda_n\right)^{1/2}
 -\left(\sum_{n=-1}^{N_0}\eps^{2n/3}\lambda_n\right)^{1/2}\right)}_{=:T_2}\nonumber\\
&&  +\frac{1}{2}\underbrace{\eps^{1/3}\chi_\eps\left(\left(\sum_{n=-1}^{N}\eps^{2n/3}\lambda_n\right)^{-1/2}
 -\left(\sum_{n=-1}^{N_0}\eps^{2n/3}\lambda_n\right)^{-1/2}\right)\sum_{n=-1}^{N}\eps^{2n/3}\nabla\lambda_n}_{=:T_3}\nonumber\\
  &&  +\frac{1}{2}\underbrace{\eps^{1/3}\chi_\eps\left(\sum_{n=-1}^{N_0}\eps^{2n/3}\lambda_n\right)^{-1/2}\sum_{n=N_0+1}^{N}\eps^{2n/3}\nabla\lambda_n}_{=:T_4}\nonumber\\
&&+\underbrace{\eps^{1/3}\nabla\left(\Psi_\eps\sum_{n=L_0+1}^{L}\eps^{2n/3}\mu_n\right)  }_{=:T_5}
 +\underbrace{\eps^{2N/3+1}\nabla Q}_{=:T_6}.
\ee
Thanks to (\ref{est-nabla1-taum-lp}) and (\ref{est-nabla2-taum-lp}), we have 
\be\label{est-T1}
\|T_1\|_{L^2(\R^d)} &= & \mathcal{O}(\eps^{(2-3\beta)(M_0+1)+\beta/2}).
\ee
$T_2$ is estimated like the second term in the right hand side of (\ref{pre-pre-conclusion-2}) in (\ref{pre-conclusion-2}), using (\ref{est-nabla1-lambdan-lp}) instead of (\ref{est-lambdan-lp}). We obtain
\be\label{est-T2}
\|T_2\|_{L^2(\R^d)} &= & \left\{\begin{array}{ll}
\mathcal{O}(\eps^{5/3-\beta}) & {\rm if\ } N_0=0\\
\mathcal{O}(\eps^{\beta(N_0-3/2)+2})
  &  {\rm if\ } N_0\geq 1.
\end{array}\right.
\ee
In order to estimate $T_3$, note that thanks to (\ref{asymp-lam0}}) and  (\ref{asymp-lamn}}), $\lambda_n'$ is uniformly bounded on $\R$ for $n=0,1,2,3$, whereas for $n\geq 4$, $\lambda_n'=\mathcal{O}_{L^\infty(D_1)}(\eps^{-(2/3-\beta)(n-3)})$. Therefore, since $\lambda_{-1}$ is constant, 
$$\sum_{n=-1}^{N}\eps^{2n/3}\nabla\lambda_n
\underset{\eps\to 0}{=}-\frac{2x}{\eps^{2/3}}\lambda_0' +\mathcal{O}_{L^\infty(D_1)}(1).$$
Applying the mean value theorem to the inverse of the square root close to $\lambda_{-1}$, we use the same arguments as to obtain (\ref{pre-conclusion-2}) from (\ref{pre-pre-conclusion-2}) and we get thanks to (\ref{est-lambdan-lp})
\be\label{est-T3}
\|T_3\|_{L^2(\R^d)} &= & \left\{\begin{array}{ll}
\mathcal{O}(\eps^{5/3}) & {\rm if\ } N_0=0\\
\mathcal{O}(\eps^{\beta(N_0-1/2)+2})
  &  {\rm if\ } N_0\geq 1.
\end{array}\right.
\ee
Lemma \ref{minor-lam} and (\ref{est-nabla2-lambdan-lp}) yield
\be\label{est-T4}
\|T_4\|_{L^2(\R^d)} &= & \left\{\begin{array}{ll}
\mathcal{O}(\eps^{2N_0/3+1}) & {\rm if\ } N_0=0\ {\rm or}\ 1 \\
\mathcal{O}(\eps^{\beta(N_0-3/2)+2})
  &  {\rm if\ } N_0\geq 2.
\end{array}\right.
\ee
It follows from (\ref{est-nabla1-mun-lp}) and
(\ref{est-nabla2-mun-lp}) that
\be\label{est-T5}
\|T_5\|_{L^2(\R^d)} &= & \mathcal{O}(\eps^{2(L_0+1)/3}).
\ee
Finally, like in (\ref{pre-conclusion-1}), we deduce from (\ref{est-R}) that
 if $M$, $N$ and $L$ are chosen large enough, $T_6$ is neglectible in
 comparison with the sum of the five other terms. Therefore, combining
 (\ref{est-T1}), (\ref{est-T2}), (\ref{est-T3}), (\ref{est-T4}) and
 (\ref{est-T5}), we obtain
\be\label{est-nabla-eta2}
\|\nabla(\eta_2-\eta_{2app})\|_{L^2(\R^d)}&=&
\mathcal{O}(\eps^{(2-3\beta)(M_0+1)+\beta/2})
+\left\{\begin{array}{ll}
\mathcal{O}(\eps^{5/3-\beta}) & {\rm if\ } N_0=0\\
\mathcal{O}(\eps^{\beta(N_0-3/2)+2})
  &  {\rm if\ } N_0\geq 1.
\end{array}\right.\nonumber\\
&&+\left\{\begin{array}{ll}
\mathcal{O}(\eps^{5/3}) & {\rm if\ } N_0=0\\
\mathcal{O}(\eps^{\beta(N_0-1/2)+2})
  &  {\rm if\ } N_0\geq 1.
\end{array}\right.
+ \left\{\begin{array}{ll}
\mathcal{O}(\eps^{2N_0/3+1}) & {\rm if\ } N_0=0\ {\rm or}\ 1 \\
\mathcal{O}(\eps^{\beta(N_0-3/2)+2})
  &  {\rm if\ } N_0\geq 2.
\end{array}\right.\nonumber\\
&&+\mathcal{O}(\eps^{2(L_0+1)/3}).\nonumber\\
&=&\left\{\begin{array}{ll}
\mathcal{O}(\eps^{(2-3\beta)(M_0+1)+\beta/2})+\mathcal{O}(\eps)+\mathcal{O}(\eps^{2(L_0+1)/3}) & {\rm if\ } N_0=0\\
\mathcal{O}(\eps^{(2-3\beta)(M_0+1)+\beta/2})+\mathcal{O}(\eps^{5/3})+\mathcal{O}(\eps^{2(L_0+1)/3})
  &  {\rm if\ } N_0= 1\\
\mathcal{O}(\eps^{(2-3\beta)(M_0+1)+\beta/2})+\mathcal{O}(\eps^{\beta(N_0-3/2)+2})+\mathcal{O}(\eps^{2(L_0+1)/3})&  {\rm if\ } N_0\geq 2.
\end{array}\right.\ \ \ \ \ \ \ \ \ 
\ee
The estimate on $\eta_1-\eta_{1app}$ in (\ref{est-main}) for $E=H^1(\R^d)$
follows from (\ref{pre-conclusion-1}) and (\ref{est-nabla-eta1}), the estimate on
$\eta_2-\eta_{2app}$ comes from (\ref{pre-conclusion-2}) and (\ref{est-nabla-eta2}).

\paragraph{Acknowledgements.} I am grateful to A. Aftalion for
bringing my attention to this problem. I also thank her and B. Noris
for helpful discussions during the work. This work is supported by the
project ANR-12-MONU-0007 BECASIM.

\end{document}